\documentclass[twoside,11pt]{article}
\usepackage{a4wide,cite,latexsym,amsfonts,amssymb,exscale,epsfig,hyperref}
\usepackage[centertags,sumlimits,intlimits,namelimits,reqno]{amsmath}

\usepackage{pstricks}
\usepackage{pst-grad}
\usepackage{pst-plot}
\usepackage[tiling]{pst-fill}
\psset{unit=1em}

\usepackage{amsthm}
\theoremstyle{definition}
\newtheorem{theorem}{Theorem}[section]
\newtheorem{lemma}[theorem]{Lemma}
\newtheorem{corollary}[theorem]{Corollary}
\newtheorem{proposition}[theorem]{Proposition}
\newtheorem{definition}[theorem]{Definition}
\newtheorem{remark}[theorem]{Remark}

\usepackage{bbm}
\def\N{{\mathbbm N}}
\def\C{{\mathbbm C}}
\def\Q{{\mathbbm Q}}
\def\Z{{\mathbbm Z}}

\def\ssl{{\mathfrak{sl}}}

\def\ie{{\sl i.e.\/}}

\def\etc{{\sl etc.\/}}
\def\cf{{\sl c.f.\/}}

\let\phi=\varphi
\let\theta=\vartheta
\let\epsilon=\varepsilon

\def\tr{\mathop{\rm tr}\nolimits}
\def\ev{\mathop{\rm ev}\nolimits}
\def\coev{\mathop{\rm coev}\nolimits}
\def\dim{\mathop{\rm dim}\nolimits}
\def\id{\mathop{\rm id}\nolimits}
\def\Span{\mathop{\rm span}\nolimits}
\def\End{\mathop{\rm End}\nolimits}
\def\Nat{\mathop{\rm Nat}\nolimits}
\def\Hom{\mathop{\rm Hom}\nolimits}
\def\im{\mathop{\rm im}\nolimits}
\def\coim{\mathop{\rm coim}\limits}

\let\hat=\widehat
\let\tilde=\widetilde
\def\pprime{{\prime\prime}}
\def\ppprime{{\prime\prime\prime}}
\def\hotimes{{\hat\otimes}}

\numberwithin{equation}{section}

\makeatletter

\newfont{\@aidxte}{cmsy10}
\newfont{\@aidxel}{cmsy10 scaled 1095}
\newfont{\@aidxtw}{cmsy10 scaled 1200}
\newlength\@aidxtexvi
\newlength\@aidxtexvii
\newlength\@aidxelxvi
\newlength\@aidxelxvii
\newlength\@aidxtwxvi
\newlength\@aidxtwxvii
\newcommand{\alignidx}[1]{%
\@aidxtexvi=\fontdimen16\@aidxte
\@aidxtexvii=\fontdimen17\@aidxte
\@aidxelxvi=\fontdimen16\@aidxel
\@aidxelxvii=\fontdimen17\@aidxel
\@aidxtwxvi=\fontdimen16\@aidxtw
\@aidxtwxvii=\fontdimen17\@aidxtw
{\mbox{$%
\fontdimen16\@aidxte=2.9pt
\fontdimen17\@aidxte=2.9pt
\fontdimen16\@aidxel=3.1pt
\fontdimen17\@aidxel=3.1pt
\fontdimen16\@aidxtw=3.3pt
\fontdimen17\@aidxtw=3.3pt
#1$}}%
\fontdimen16\@aidxte=\@aidxtexvi
\fontdimen17\@aidxte=\@aidxtexvii
\fontdimen16\@aidxel=\@aidxelxvi
\fontdimen17\@aidxel=\@aidxelxvii
\fontdimen16\@aidxtw=\@aidxtwxvi
\fontdimen17\@aidxtw=\@aidxtwxvii}

\makeatother

\newenvironment{myenumerate}{%

\begin{enumerate}
\setlength{\partopsep}{0pt}
\setlength{\parskip}{0pt}}{\end{enumerate}}

\def\nn{\notag}

\def\emph#1{{\sl #1\/}}
\def\sym#1{{\mathcal #1}}
\def\one{\mathbbm{1}}%
\def\bar#1{\overline{#1}}%

\def\op{\mathrm{op}}
\def\cop{\mathrm{cop}}

\def\Set{\mathbf{Set}}
\def\Vect{\mathbf{Vect}}
\def\fdVect{\mathbf{fdVect}}
\def\CoAlg{\mathbf{CoAlg}}
\def\coend{\mathbf{coend}}

\def\msc#1{\noindent{\small Mathematics Subject Classification (2000):
#1\par}}%
\def\keywords#1{\noindent {\small keywords: #1\par}}%

\newcommand{\luturn}{
\psbezier(-0.5,2)(-0.5,1)(-1,1.5)(-1,0)
\psline(-1,0)(-0.5,0)
\psbezier(-0.5,0)(-0.5,0.9)(-0.4,1)(0,1)
\psline(0,1)(0,2)
\psline(0,2)(-0.5,2)
}

\newcommand{\lucurve}{
\psbezier(-0.5,2)(-0.5,1)(-1,1.5)(-1,0)
\psline(-1,0)(-0.5,0)
\psbezier(-0.5,0)(-0.5,1.5)(0,1)(0,2)
\psline(0,2)(-0.5,2)
}

\newcommand{\ruturn}{
\psbezier(0.5,2)(0.5,1)(1,1.5)(1,0)
\psline(1,0)(0.5,0)
\psbezier(0.5,0)(0.5,0.9)(0.4,1)(0,1)
\psline(0,1)(0,2)
\psline(0,2)(0.5,2)
}

\newcommand{\rucurve}{
\psbezier(0.5,2)(0.5,1)(1,1.5)(1,0)
\psline(1,0)(0.5,0)
\psbezier(0.5,0)(0.5,1.5)(0,1)(0,2)
\psline(0,2)(0.5,2)
}

\newcommand{\llturn}{
\psbezier(-0.5,-2)(-0.5,-1)(-1,-1.5)(-1,0)
\psline(-1,0)(-0.5,0)
\psbezier(-0.5,0)(-0.5,-0.9)(-0.4,-1)(0,-1)
\psline(0,-1)(0,-2)
\psline(0,-2)(-0.5,-2)
}

\newcommand{\llcurve}{
\psbezier(-0.5,-2)(-0.5,-1)(-1,-1.5)(-1,0)
\psline(-1,0)(-0.5,0)
\psbezier(-0.5,0)(-0.5,-1.5)(0,-1)(0,-2)
\psline(0,-2)(-0.5,-2)
}

\newcommand{\rlturn}{
\psbezier(0.5,-2)(0.5,-1)(1,-1.5)(1,0)
\psline(1,0)(0.5,0)
\psbezier(0.5,0)(0.5,-0.9)(0.4,-1)(0,-1)
\psline(0,-1)(0,-2)
\psline(0,-2)(0.5,-2)
}

\newcommand{\rlcurve}{
\psbezier(0.5,-2)(0.5,-1)(1,-1.5)(1,0)
\psline(1,0)(0.5,0)
\psbezier(0.5,0)(0.5,-1.5)(0,-1)(0,-2)
\psline(0,-2)(0.5,-2)
}

\newcommand{\straight}{
\psline(0,0)(0,2)
\psline(0,2)(0.5,2)
\psline(0.5,2)(0.5,0)
\psline(0.5,0)(0,0)
}

\newcommand{\trivalent}[3]{
\psline(0,-1)(0,0)
\psline(0,0)(-1,1)
\psline(0,0)(1,1)
\pscircle[fillstyle=solid,fillcolor=black](0,0){0.15}
\rput(-1.5,0.8){$\scriptstyle #1$}
\rput(1.5,0.8){$\scriptstyle #2$}
\rput(0.5,-1){$\scriptstyle #3$}
}

\newcommand{\trivalentdown}[3]{
\psline(0,1)(0,0)
\psline(0,0)(-1,-1)
\psline(0,0)(1,-1)
\pscircle[fillstyle=solid,fillcolor=black](0,0){0.15}
\rput(0.5,1.0){$\scriptstyle #1$}
\rput(-1.5,-0.8){$\scriptstyle #2$}
\rput(1.5,-0.8){$\scriptstyle #3$}
}

\newcommand{\boxbasis}[1]{
\psframe(-1,0.5)(1,-0.5)
\rput(0,0){$\scriptstyle #1$}
\psline(0,1)(0,0.5)
\psline(-0.5,-0.5)(-0.5,-1)
\psline(0.5,-0.5)(0.5,-1)
}

\newcommand{\boxbasisdual}[1]{
\psframe(-1,0.5)(1,-0.5)
\rput(0,0){$\scriptstyle #1$}
\psline(0,-0.5)(0,-1)
\psline(-0.5,0.5)(-0.5,1)
\psline(0.5,0.5)(0.5,1)
}

\newcommand{\boxbasiswide}[1]{
\psframe(-1,0.5)(1.5,-0.5)
\rput(0.25,0){$\scriptstyle #1$}
\psline(0.25,1)(0.25,0.5)
\psline(-0.5,-0.5)(-0.5,-1)
\psline(1,-0.5)(1,-1)
}

\newcommand{\boxbasisdualwide}[1]{
\psframe(-1,0.5)(1.5,-0.5)
\rput(0.25,0){$\scriptstyle #1$}
\psline(0.25,-0.5)(0.25,-1)
\psline(-0.5,0.5)(-0.5,1)
\psline(1,0.5)(1,1)
}

\newcommand{\circlebox}[1]{
\pscircle(0,0){0.75}
\rput(0,0){$\scriptstyle #1$}
\psline(0,0.75)(0,1)
\psline(0,-0.75)(0,-1)
}

\newcommand{\circlebig}[1]{
\pscircle(0,0){0.9}
\rput(0,0){$\scriptstyle #1$}
\psline(0,0.9)(0,1)
\psline(0,-0.9)(0,-1)
}

\usepackage[ps,matrix,arrow,curve]{xy}
\xyoption{dvips}

\begin{document}

\title{
Tannaka--Kre\v\i n reconstruction and a characterization\\
of modular tensor categories}
\author{Hendryk Pfeiffer\thanks{E-mail: \texttt{pfeiffer@math.ubc.ca}}}
\date{\small{Department of Mathematics, University of British Columbia,\\
1984 Mathematics Road, Vancouver, BC, V2T 1Z2, Canada}\\[1ex]
January 30, 2008}

\maketitle

\begin{abstract}

We show that every modular category is equivalent as an additive
ribbon category to the category of finite-dimensional comodules of a
Weak Hopf Algebra. This Weak Hopf Algebra is finite-dimensional, split
cosemisimple, weakly cofactorizable, coribbon and has trivially
intersecting base algebras. In order to arrive at this
characterization of modular categories, we develop a generalization of
Tannaka--Kre\v\i n reconstruction to the long version of the canonical
forgetful functor which is lax and oplax monoidal, but not in general
strong monoidal, thereby avoiding all the difficulties related to
non-integral Frobenius--Perron dimensions. In the more general case of a
finitely semisimple additive ribbon category, not necessarily modular, the
reconstructed Weak Hopf Algebra is finite-dimensional, split cosemisimple,
coribbon and has trivially intersecting base algebras.

\end{abstract}

\msc{%
16W30,
18D10
}

\section{Introduction}

A modular category~\cite{ReTu91,Tu94} is a finitely semisimple
additive ribbon category that satisfies a non-degeneracy
condition. For the precise definition, see Section~\ref{sect_modular}
below. Modular categories are of interest in a variety of areas from
low-dimensional topology~\cite{ReTu91,Tu94} to Conformal Field
Theory~\cite{MoSe89,FuRu02}, subfactor theory~\cite{Mu03} and
$3$-dimensional quantum gravity~\cite{Ba03}.

It is well known that some modular categories are equivalent to
categories of the form ${}_A\sym{M}$, the category of left $A$-modules
of some suitable ring or algebra $A$. Since algebras are often easier
to deal with than categories, it is an interesting problem to
understand whether all modular categories are of this form. We show
that this is indeed the case. We restrict ourselves to modular
categories $\sym{C}$ for which the commutative ring $k=\End(\one)$,
\ie\ the endomorphisms of the monoidal unit object, is a field.

In many cases, it is outright obvious that a modular category
$\sym{C}$ is equivalent to the category ${}_A\sym{M}$ for some
$k$-algebra $A$. This is the case, for example, if all simple objects
$V_j$ of $\sym{C}$ are finite-dimensional vector spaces over some
field $k$. Since $\sym{C}$ is by definition finitely
semisimple\footnote{The assumption of finite semisimplicity includes
that $\End(V_j)\cong k$ for each simple object $V_j$ of $\sym{C}$.},
$A$ is just a finite direct sum of the appropriate $n_j\times
n_j$-matrix algebras with coefficients in $k$ where $n_j=\dim_k
V_j$. The equivalence $\sym{C}\simeq{}_A\sym{M}$ is just an
equivalence of (ordinary) categories, and one still needs to determine
which additional structure and properties of $A$ give rise to the
monoidal structure, braiding, ribbon structure and special properties
of $\sym{C}$.

Which sort of additional structure on a $k$-algebra $A$ would be
sufficient in order to equip the category ${}_A\sym{M}$ of left
$A$-modules with the structure of a monoidal category? The most widely
known answer to this question is that one can employ the structure of
a Hopf algebra or a bialgebra. It is further known, for example, that
ribbon Hopf algebras $H$~\cite{ReTu90}, a special sort of
quasitriangular Hopf algebras, have categories ${}_H\sym{M}$ of left
$H$-modules that carry the structure of ribbon categories. In order to
answer the converse question, \ie\ which ribbon categories are of the
form ${}_H\sym{M}$ for some $k$-algebra $H$, Tannaka--Kre\v\i n
reconstruction~\cite{Sa72,DeMi82} was generalized from (the coordinate
rings of) groups to Hopf algebras, ribbon Hopf algebras and even quasi
Hopf algebras, see, for
example~\cite{Ul90,Ma91b,JoSt91b,Sc92,Ma95}. These constructions
successfully deal with the additional structure such as duality,
braiding and the ribbon structure, but quite a basic problem with the
monoidal structure is left unsolved.

The problem is that not every rigid monoidal category is monoidally equivalent
to the category of modules over a Hopf algebra or a quasi Hopf algebra. For
example, a modular category $\sym{C}$ with $\End(\one)=\C$ is the category of
$H$-modules for some finite-dimensional quasi Hopf algebra $H$ if and only if
each simple object of $\sym{C}$ has an integer Frobenius--Perron
dimension~\cite[Theorem~8.33]{EtNi05}. But there exist interesting examples of
modular categories that contain objects of non-integer Frobenius--Perron
dimension\footnote{Important examples of modular categories are related to (a
finite-dimensional version of) the Hopf algebras $U_q(\mathfrak{g})$, $q$ a
root of unity, see, for example~\cite{ReTu91,Sa06}. These are, however, not
${}_{U_q(\mathfrak{g})}\sym{M}$. One rather has to take first the full
subcategory of tilting modules of ${}_{U_q(\mathfrak{g})}\sym{M}$ and then a
quotient of that subcategory, dividing out the negligible morphisms. The
resulting modular category is in general no longer of the form ${}_H\sym{M}$
for any Hopf algebra $H$.}.

In order to deal with non-integer Frobenius--Perron dimensions,
B{\"o}hm, Nill and Szlach{\'a}nyi have invented the concept of Weak
Bialgebras (WBAs) and Weak Hopf Algebras
(WHAs)~\cite{BoSz96,Ni98,BoNi99,BoSz00,Ni02,Ni04}. B{\"o}hm's
thesis~\cite{Bo97} contains the first examples of modular categories
which have objects of non-integer Frobenius--Perron dimension and
which are shown to be the categories of modules of some
finite-dimensional WBA. The definitions of a WBA and of a WHA are
summarized in detail in Section~\ref{sect_wha} below.

Is the concept of a WBA general enough in order to show that every
modular category $\sym{C}$ is equivalent (first as a monoidal and then
as a ribbon category) to the category of modules of some WBA $H$? It
is useful to subdivide this question into the following three steps:
\begin{myenumerate}
\item
Can every object $X\in|\sym{C}|$ be viewed as a $k$-vector space for
some $k$?
\item
Does the monoidal structure of $\sym{C}$ arise from the WBA
structure of $H$?
\item
Which additional structure and properties of $H$ are required in
order to obtain duality, braiding and ribbon structure of $\sym{C}$
and in order to satisfy the non-degeneracy condition?
\end{myenumerate}

Question~(1) was answered by Hayashi~\cite{Ha99} who showed that there
is a canonical forgetful functor
$\hat\omega\colon\sym{C}\to{}_R\sym{M}_R$ into the category
${}_R\sym{M}_R$ of $(R,R)$-bimodules. Here $R=\End(\hat V)$ is the
commutative $k$-algebra, $k=\End(\one)$, of endomorphisms of the
universal object,
\begin{equation}
\hat V=\bigoplus_{j\in I} V_j,
\end{equation}
the direct sum over one representative $V_j$ for each isomorphism
class of simple objects.

This functor $\hat\omega\colon\sym{C}\to{}_R\sym{M}_R$ is now known as
the \emph{short forgetful functor}. In order to solve question~(1)
above, one composes it with the forgetful functor
${}_R\sym{M}_R\to\Vect_k$ that assigns to each $(R,R)$-bimodule the
underlying $k$-vector space, and thereby obtains the \emph{long
forgetful functor} $\omega\colon\sym{C}\to\Vect_k$. While the short
forgetful functor is strong monoidal, the long forgetful functor is in
general not strong monoidal and therefore not a \emph{fibre functor}
in the usual technical sense. Szlach{\'a}nyi~\cite{Sz05} has
characterized those long forgetful functors that originate from the
categories of modules of WBAs.

Hayashi~\cite{Ha99} and H{\'a}i~\cite{Ha02} have studied the generalization of
Tannaka--Kre\v\i n reconstruction to the case of the short forgetful functor,
\ie\ to a strong monoidal functor into the bimodule category
${}_R\sym{M}_R$. It is known that the reconstructed algebraic structure is a
bialgebroid over $R$ and, furthermore, since $R$ is a finite-dimensional
separable commutative $k$-algebra, one actually gets a
WBA~\cite{Sz05}. Therefore, Ostrik~\cite{Os03} concludes from these abstract
considerations that the answer to question~(2) above is `yes'.

Tannaka--Kre\v\i n reconstruction using the short forgetful functor
$\hat\omega\colon\sym{C}\to{}_R\sym{M}_R$ alone, however, uses the
language of bialgebroids, and it is thus not transparent how duality,
braiding and ribbon structure carry over from the modular category to
the reconstructed WBA.

It is the purpose of the present article to complete the programme of
Tannaka--Kre\v\i n reconstruction including question~(3) above, and to
prove the following

\begin{theorem}
\label{thm_main}
Every modular category for which $k=\End(\one)$ is a field, is
equivalent as a $k$-linear ribbon category to the category of
finite-dimensional comodules of a finite-dimensional split
cosemisimple weakly cofactorizable coribbon WHA over $k$ whose base
algebras intersect trivially.
\end{theorem}

This theorem also holds, more generally, without the non-degeneracy
condition on the $S$-matrix:
\begin{theorem}
Every finitely semisimple additive ribbon category for which $k=\End(\one)$ is
a field, is equivalent as a $k$-linear ribbon category to the category of
finite-dimensional comodules of a finite-dimensional split cosemisimple
coribbon WHA over $k$ whose base algebras intersect trivially.
\end{theorem}

We reconstruct this WHA, characterize all its operations by the
universal property of the appropriate coend, \ie\ the universal
coacting coalgebra, and also write down the operations in terms of a
convenient basis.

Several authors have given sufficient conditions for the category of
modules ${}_A\sym{M}$ of some $k$-algebra $A$ to be modular, see, for
example~\cite[Lemma 1.1]{EtGe98} for Drinfel'd doubles of Hopf
algebras and~\cite[Lemma 8.2]{NiTu03} for WHAs. As far as we know,
Theorem~\ref{thm_main} is the first one to establish the precise form
of the converse implication, \ie\ that every modular category can
indeed be obtained from a WHA with the properties stated.

In order to prove Theorem~\ref{thm_main}, we generalize
Tannaka--Kre\v\i n reconstruction to the long forgetful functor
$\omega\colon\sym{C}\to\Vect_k$. Since this functor has the category
$\Vect_k$ of vector spaces over $k$ as its codomain, reconstruction
immediately yields a coalgebra object\footnote{For the usual technical
reasons, \ie\ because we want to exploit that the category $\Vect_k$
is (small) cocomplete and its tensor product preserves colimits in
both arguments, we prefer to reconstruct a coalgebra rather than an
algebra. For more details, we refer to Section~\ref{sect_tkreview}.}
$H$ in $\Vect_k$. This is substantially more transparent than a
functor into the bimodule category and allows us to recover all
additional operations of $H$ by exploiting the universal property of
the coend. We emphasize that the long forgetful functor
$\omega\colon\sym{C}\to\Vect_k$ is in general not strong monoidal, but
nevertheless both lax and oplax monoidal~\cite{Sz05}, and we have to
generalize Tannaka--Kre\v\i n reconstruction to this case. It is this
property of being lax and oplax rather than strong monoidal that
enables us to deal with non-integer Frobenius--Perron dimensions.

The present article is structured as follows. In
Section~\ref{sect_prelim}, we review the definitions and some key
results on modular categories, WHAs, comodules, and on
Tannaka--Kre\v\i n reconstruction. In Section~\ref{sect_long}, we
study the properties of the long forgetful functor. We reconstruct a
coribbon WHA from each modular category in
Section~\ref{sect_reconstruct}. In Section~\ref{sect_corep}, we study
the category of finite-dimensional comodules of a coribbon WHA, and in
Section~\ref{sect_equiv}, we show that the original modular category
is equivalent to the category of finite-dimensional comodules of the
reconstructed coribbon WHA. For convenient reference, we
compile the relevant definitions and results about monoidal categories
in Appendix~\ref{app_tensor}.

The reader who just wants to get a quick overview of how the
reconstructed WHA looks like, without going through all the technical
details, is invited to go straight to Section~\ref{sect_example} where
we present the reconstructed WHA for the modular category associated
with $U_q(\ssl_2)$, $q$ a root of unity, in term of the familiar
diagrams.

\section{Preliminaries}
\label{sect_prelim}

\subsection{Modular categories}
\label{sect_modular}

In this section, we summarize the definition and some basic properties
of modular categories. For more details, we refer to the
book~\cite{Tu94}.

Our notation is as follows. If $\sym{C}$ is a category, we write
$X\in|\sym{C}|$ for the objects $X$ of $\sym{C}$, $\Hom(X,Y)$ for the
collection of all morphisms $f\colon X\to Y$ and
$\End(X)=\Hom(X,X)$. By $\id_X\colon X\to X$ we denote the identity
morphism of $X$ and by $g\circ f\colon X\to Z$ the composition of
morphisms $f\colon X\to Y$ and $g\colon Y\to Z$. If two objects
$X,Y\in|\sym{C}|$ are isomorphic, we write $X\cong Y$. If two
categories are equivalent, we write $\sym{C}\simeq\sym{D}$. The
identity functor on $\sym{C}$ is denoted by $1_{\sym{C}}$, and
$\sym{C}^\op$ is the opposite category of $\sym{C}$. The category of
vector spaces over a field $k$ is denoted by $\Vect_k$ and its full
subcategory of finite-dimensional vector spaces by $\fdVect_k$.

We assume that the reader is familiar with the notions of
$\mathbf{Ab}$-enriched, additive, abelian, monoidal, braided monoidal,
autonomous and ribbon categories. For convenience, we have compiled
the relevant definitions in Appendix~\ref{app_tensor}.

\begin{definition}
\label{def_modular}
A \emph{modular category}
$(\sym{C},\otimes,\one,\alpha,\lambda,\rho,{(-)}^\ast,\ev,\coev,\sigma,\nu)$
is an additive ribbon category (\cf\ Definitions~\ref{def_ribboncat}
and~\ref{def_preadditivecat}) that satisfies the following conditions:
\begin{myenumerate}
\item
$k=\End(\one)$ is a field.
\item
There is a finite family ${\{V_j\}}_{j\in I}$ of objects
$V_j\in|\sym{C}|$ where $I$ denotes some finite index set, that
satisfies the following conditions:
\begin{myenumerate}
\item
Each $V_j$, $j\in I$, satisfies $\End(V_j)\cong k$, \ie\ it is
\emph{simple}.
\item
There is an element $0\in I$ such that $V_0\cong\one$.
\item
For each $j\in I$, there is some $j^\ast\in I$ such that $V_{j^\ast}\cong{(V_j)}^\ast$.
\item
For each object $X\in|\sym{C}|$, there is a finite sequence
$(j^X_1,\ldots,j^X_{n^X})\in I^{n^X}$, $n^X\in\N_0$, and
morphisms $\imath^X_\ell\colon V_{j^X_\ell}\to X$ and
$\pi_\ell^X\colon X\to V_{j^X_\ell}$ for $1\leq\ell\leq n$ such
that
\begin{equation}
\id_X = \sum_{\ell=1}^{n^X} \imath^X_\ell\circ\pi^X_\ell.
\end{equation}
\end{myenumerate}
\item
The matrix ${(S_{ij})}_{i,j\in I}$ (\emph{$S$-matrix}) whose
coefficients are
\begin{equation}
\label{eq_smatrix}
S_{ij} = \tr_{V_i\otimes V_j}(\sigma_{V_j,V_i}\circ\sigma_{V_i,V_j})\in k,
\end{equation}
is invertible.
\end{myenumerate}
\end{definition}

Compared with the definition of Turaev~\cite{Tu94}, we have added in
our Definition~\ref{def_modular} the conditions that $k=\End(\one)$ be
a field and that $\sym{C}$ be additive rather than just
$\mathbf{Ab}$-enriched. The former is related to the fact that we
reconstruct a WHA over $k$ and we only deal with the case in which
this is a field. The latter makes sure that $\sym{C}$ has all finite
biproducts ('direct sums'). Otherwise, one could remove some of the
objects of $\sym{C}$ that are biproducts of simple objects, without
violating any condition of the definition. We disallow this because we
want to compare $\sym{C}$ to the category of comodules of the
reconstructed WHA which automatically has all finite biproducts.

Note that in the definition of a modular category, one usually
requires $\End(V_j)\cong k$ for the simple objects although one does
not impose any restriction on the field $k$ such as algebraic
closure. Many algebraic examples of modular categories, see, for
example~\cite{Sa06}, have $k=\Q(\epsilon)$, a cyclotomic extension of
the rationals, far from algebraically closed, and nevertheless
$\End(V_j)\cong k$ for all simple objects.

\begin{proposition}
\label{prop_modular}
Let $\sym{C}$ be a modular category, $k=\End(\one)$, and
${\{V_j\}}_{j\in I}$ be a family of objects as in
Definition~\ref{def_modular}(2).
\begin{myenumerate}
\item
$\sym{C}$ is $k$-linear as a monoidal category (\cf\
Definition~\ref{def_preadditivecat})~\cite[Section I.1.5]{Tu94}.
\item
For all objects $X,Y\in|\sym{C}|$, the abelian group $\Hom(X,Y)$ is
a finite-dimensional vector space over $k$~\cite[Lemma II.4.2.1]{Tu94}.
\item
$\sym{C}$ is non-degenerate, \ie\ its traces define non-degenerate
bilinear forms (\cf\ Definition~\ref{def_nondegenerate})~\cite[Lemma
II.4.2.3]{Tu94}.
\item
The morphisms $\imath_\ell^X$ and $\pi_\ell^X$ of
Definition~\ref{def_modular}(2d) can be chosen in such a way that
\begin{equation}
\pi^X_\ell\circ\imath^X_m = \left\{
\begin{matrix}
\id_{V_{j^X_\ell}},&\mbox{if}\quad \ell=m,\\
0,                 &\mbox{else}
\end{matrix}
\right.
\end{equation}
for all $\ell,m\in I$ (Proposition~\ref{prop_nondegenerate}).
\item
If $j,\ell\in I$ and $j\neq\ell$, then $\Hom(V_j,V_\ell)=\{0\}$~\cite[Lemma II.1.5]{Tu94}.
\item
$\sym{C}$ is finitely semisimple according to
Definition~\ref{def_semisimple}(3).
\item
If $X\in|\sym{C}|$ is simple, then $0\neq \dim(X)\in k$~\cite[Lemma II.4.2.4]{Tu94}.
\item
If $X\in|\sym{C}|$ is simple, then there exists some $j\in I$ such
that $X\cong V_j$ (Corollary~\ref{cor_allsimplethere}).
\end{myenumerate}
\end{proposition}

Most results of this article already hold without the non-degeneracy condition
on the $S$-matrix~\eqref{eq_smatrix}, \ie\ for finitely semisimple additive
ribbon categories for which $k=\End(\one)$ is a field.

\subsection{Weak Hopf algebras}
\label{sect_wha}

In this section, we summarize the definitions of a Weak Bialgebra
(WBA) and of a Weak Hopf Algebra (WHA). For more details, we refer
to~\cite{Ni98,BoNi99,BoSz00,Ni02,Ni04}.

\begin{definition}
\label{def_wba}
A \emph{Weak Bialgebra} $(H,\mu,\eta,\Delta,\epsilon)$ over a field
$k$ is a $k$-vector space $H$ with linear maps $\mu\colon H\otimes
H\to H$ (\emph{multiplication}), $\eta\colon k\to H$ (\emph{unit}),
$\Delta\colon H\to H\otimes H$ (\emph{comultiplication}), and
$\epsilon\colon H\to k$ (\emph{counit}) such that the following
conditions hold:
\begin{myenumerate}
\item
$(H,\mu,\eta)$ is an associative unital algebra, \ie\
$\mu\circ(\mu\otimes\id_H)=\mu\circ(\id_H\otimes\mu)$ and
$\mu\circ(\eta\otimes\id_H)=\id_H=\mu\circ(\id_H\otimes\eta)$.
\item
$(H,\Delta,\epsilon)$ is a coassociative counital coalgebra, \ie\
$(\Delta\otimes\id_H)\circ\Delta=(\id_H\otimes\Delta)\circ\Delta$
and
$(\epsilon\otimes\id_H)\circ\Delta=\id_H=(\id_H\otimes\epsilon)\circ\Delta$.
\item
The following compatibility conditions hold:
\begin{eqnarray}
\label{eq_wba1}
\Delta\circ\mu
&=& (\mu\otimes\mu)\circ(\id_H\otimes\sigma_{H,H}\otimes\id_H)\circ(\Delta\otimes\Delta),\\
\label{eq_wba2}
\epsilon\circ\mu\circ(\mu\otimes\id_H)
&=& (\epsilon\otimes\epsilon)\circ(\mu\otimes\mu)\circ(\id_H\otimes\Delta\otimes\id_H)\nn\\
&=& (\epsilon\otimes\epsilon)\circ(\mu\otimes\mu)\circ(\id_H\otimes\Delta^\op\otimes\id_H),\\
\label{eq_wba3}
(\Delta\otimes\id_H)\circ\Delta\circ\eta
&=& (\id_H\otimes\mu\otimes\id_H)\circ(\Delta\otimes\Delta)\circ(\eta\otimes\eta)\nn\\
&=& (\id_H\otimes\mu^\op\otimes\id_H)\circ(\Delta\otimes\Delta)\circ(\eta\otimes\eta).
\end{eqnarray}
\end{myenumerate}
Here $\sigma_{V,W}\colon V\otimes W\to W\otimes V$, $v\otimes w\mapsto
w\otimes v$ is the transposition of the tensor factors in $\Vect_k$, and
by $\Delta^\op=\sigma_{H,H}\circ\Delta$ and $\mu^\op=\mu\circ\sigma_{H,H}$, we
denote the \emph{opposite comultiplication} and \emph{opposite
multiplication}, respectively. We tacitly identify the vector spaces
$(V\otimes W)\otimes U\cong V\otimes(W\otimes U)$ and $V\otimes k\cong V\cong
k\otimes V$, exploiting the coherence theorem for the monoidal category
$\Vect_k$.
\end{definition}

We use the term \emph{comultiplication} for the operation $\Delta$ in
a coalgebra whereas \emph{coproduct} always refers to a colimit in a
category.

\begin{definition}
A \emph{homomorphism} $\phi\colon H\to H^\prime$ of WBAs
$(H,\mu,\eta,\Delta,\epsilon)$ and
$(H^\prime,\mu^\prime,\Delta^\prime,\epsilon^\prime)$ over the same
field $k$ is a $k$-linear map that is a homomorphism of unital
algebras, \ie\ $\phi\circ\eta=\eta^\prime$ and
$\phi\circ\mu=\mu^\prime\circ(\phi\otimes\phi)$, as well as a
homomorphism of counital coalgebras, \ie\
$\epsilon^\prime\circ\phi=\epsilon$ and
$\Delta^\prime\circ\phi=(\phi\otimes\phi)\circ\Delta$.
\end{definition}

\begin{definition}
\label{def_sourcetarget}
Let $(H,\mu,\eta,\Delta,\epsilon)$ be a WBA. The linear maps
$\epsilon_t\colon H\to H$ (\emph{target counital map}) and
$\epsilon_s\colon H\to H$ (\emph{source counital map}) are defined by
\begin{eqnarray}
\label{eq_epsilont}
\epsilon_t&:=&(\epsilon\otimes\id_H)\circ(\mu\otimes\id_H)\circ(\id_H\otimes\sigma_{H,H})
\circ(\Delta\otimes\id_H)\circ(\eta\otimes\id_H),\\
\label{eq_epsilons}
\epsilon_s&:=&(\id_H\otimes\epsilon)\circ(\id_H\otimes\mu)\circ(\sigma_{H,H}\otimes\id_H)
\circ(\id_H\otimes\Delta)\circ(\id_H\otimes\eta).
\end{eqnarray}
\end{definition}

Both $\epsilon_t$ and $\epsilon_s$ are idempotents. A WBA
$(H,\mu,\eta,\Delta,\epsilon)$ is a bialgebra if and only if
$\Delta\circ\eta=\eta\otimes\eta$, if and only if
$\epsilon\circ\mu=\epsilon\otimes\epsilon$, if and only if
$\epsilon_s=\eta\circ\epsilon$ and if and only if
$\epsilon_t=\eta\circ\epsilon$.

\begin{proposition}
\label{prop_hsht}
Let $(H,\mu,\eta,\Delta,\epsilon)$ be a WBA.
\begin{myenumerate}
\item
The subspace $H_t:=\epsilon_t(H)$ (\emph{target base algebra}) forms
a subalgebra with unit and a left coideal, \ie\
\begin{equation}
\Delta(H_t)\subseteq H\otimes H_t.
\end{equation}
\item
The subspace $H_s:=\epsilon_s(H)$ (\emph{source base algebra}) forms
a subalgebra with unit and a right coideal, \ie\
\begin{equation}
\label{eq_rightcoideal}
\Delta(H_s)\subseteq H_s\otimes H.
\end{equation}
\item
The subalgebras $H_s$ and $H_s$ commute, \ie\ $xy=yx$ for all $x\in
H_t$ and $y\in H_s$.
\end{myenumerate}
\end{proposition}

\begin{definition}
A \emph{Weak Hopf Algebra} $(H,\mu,\eta,\Delta,\epsilon,S)$ is a Weak
Bialgebra $(H,\mu,\eta,\Delta,\epsilon)$ with a linear map $S\colon
H\to H$ (\emph{antipode}) that satisfies the following conditions:
\begin{eqnarray}
\label{eq_wha1}
\mu\circ(\id_H\otimes S)\circ\Delta &=& \epsilon_t,\\
\label{eq_wha2}
\mu\circ(S\otimes\id_H)\circ\Delta &=& \epsilon_s,\\
\label{eq_wha3}
\mu\circ(\mu\otimes\id_H)\circ(S\otimes\id_H\otimes S)
\circ(\Delta\otimes\id_H)\circ\Delta&=&S.
\end{eqnarray}
\end{definition}

For convenience, we write $1=\eta(1)$ and omit parentheses in products,
exploiting associativity. We also use Sweedler's notation and write
$\Delta(x)=x^\prime\otimes x^\pprime$ for the comultiplication of $x\in
H$ as an abbreviation of the expression $\Delta(x)=\sum_k a_k\otimes b_k$ with
some $a_k,b_k\in H$. Similarly, we write
$((\Delta\otimes\id_H)\circ\Delta)(x)=x^\prime\otimes x^\pprime\otimes
x^\ppprime$, exploiting coassociativity. Then, for example,
equation~\eqref{eq_epsilont} reads $\epsilon_t(x) = \epsilon(1^\prime
x)1^\pprime$ for all $x\in H$.

The concepts of a WBA and of a WHA are formally self-dual, \ie\ if $H$
is a [WBA, WHA] that is finite-dimensional as a vector space, then its
dual space $H^\ast$ is a [WBA, WHA] as well.

\begin{definition}
A \emph{homomorphism} $\phi\colon H\to H^\prime$ of WHAs is a
homomorphism of WBAs for which $\phi\circ S=S^\prime\circ\phi$.
\end{definition}

If $H$ is a WHA, we denote by $H^\op$ its opposite WHA
$(H,\mu^\op,\eta,\Delta,\epsilon,S^{-1})$, by $H^\cop$ its coopposite
WHA $(H,\mu,\eta,\Delta^\op,\epsilon,S^{-1})$ and by $H^{\op,\cop}$
the WHA $(H,\mu^\op,\eta,\Delta^\op,\epsilon,S)$. The antipode of a
WHA is an algebra antihomomorphism, \ie\
$S\circ\mu=\mu^\op\circ(S\otimes S)$ and $S\circ\eta=\eta$, as well as
a coalgebra antihomomorphism, \ie\ $(S\otimes
S)\circ\Delta=\Delta^\op\circ S$ and $\epsilon\circ S=\epsilon$.

\begin{definition}[see~\cite{Ni02,Ni04}]
Let $(H,\mu,\eta,\Delta,\epsilon,S)$ be a WHA.
\begin{myenumerate}
\item
The \emph{minimal Weak Hopf Algebra} $H_{\mathrm{min}}$ of $H$ is
the smallest sub WHA of $H$ that contains the unit $\eta(1)\in H$.
\item
$H$ is called \emph{regular} if
$S^2|_{H_{\mathrm{min}}}=\id_{H_{\mathrm{min}}}$.
\end{myenumerate}
\end{definition}

\subsection{Coalgebras and comodules}

The following definitions and results can be found, for example,
in~\cite{Gr76}.

\begin{definition}
Let $(C,\Delta,\epsilon)$ be a coalgebra over some field $k$. A right
$C$-comodule $(V,\beta_V)$ is a $k$-vector space $V$ with a $k$-linear
map $\beta_V\colon V\to V\otimes C$ that satisfies
\begin{eqnarray}
\label{eq_comodule1}
(\id_V\otimes\epsilon)\circ\beta_V &=& \id_V,\\
\label{eq_comodule2}
(\beta_V\otimes\id_H)\circ\beta_V  &=& (\id_V\otimes\Delta)\circ\beta_V.
\end{eqnarray}
\end{definition}

\begin{definition}
Let $(C,\Delta,\epsilon)$ be a coalgebra and $(V,\beta_V)$ and
$(W,\beta_W)$ be right $C$-comodules. A \emph{morphism of coalgebras}
$f\colon V\to W$ is a $k$-linear map that satisfies
\begin{equation}
(f\otimes\id_C)\circ\beta_V = \beta_W\circ f.
\end{equation}
\end{definition}

We extend Sweedler's notation to comodules and write
$\beta(v)=v_V\otimes v_C$. The conditions~\eqref{eq_comodule1}
and~\eqref{eq_comodule2} then read $v_V\epsilon(v_C)=v$ and
${(v_V)}_V\otimes{(v_V)}_C\otimes
v_C=v_V\otimes{(v_C)}^\prime\otimes{(v_C)}^\pprime$.

\begin{proposition}
Let $(C,\Delta,\epsilon)$ be a coalgebra over some field $k$ and
$\sym{M}^C$ be the category whose objects are the right $C$-comodules
that are finite-dimensional as vector spaces over $k$, and whose
morphisms are morphisms of right $C$-comodules. Then $\sym{M}^C$ is
$k$-linear and abelian, and $\Hom(V,W)$ is finite-dimensional for all
$V,W\in|\sym{M}^C|$.
\end{proposition}

If $V$ is a finite-dimensional right $C$-comodule with basis
${(v_j)}_j$, then there are elements $c_{\ell j}\in C$ uniquely
determined by the condition that $\beta_V(v_j)=\sum_\ell v_\ell\otimes
c_{\ell j}$. They are called the \emph{coefficients of} $V$ with
respect to that basis. They span the \emph{coefficient coalgebra}
$C(V)=\Span_k\{c_{\ell j}\}$, a sub coalgebra of $C$.

Let $W$ be a finite-dimensional vector space over $k$ with dual space
$W^\ast$ and a pair of dual bases ${(e_j)}_j$ and ${(e^j)}_j$ of $W$
and $W^\ast$, respectively. We abbreviate $c_{jk}=e^j\otimes e_k\in
W^\ast\otimes W$. The coalgebra $(W^\ast\otimes W,\Delta,\epsilon)$
with $\Delta(c_{jk})=\sum_\ell c_{j\ell}\otimes c_{\ell k}$ and
$\epsilon(c_{jk})=\delta_{jk}$ is called the \emph{matrix coalgebra}
associated with $W$. In this case, $W$ is a right $W^\ast\otimes
W$-comodule, and $W^\ast\otimes W$ is its coefficient coalgebra.

\begin{definition}
A coalgebra $(C,\Delta,\epsilon)$ over a field $k$ is called
\emph{cosimple} if $C$ has no sub coalgebras other than $C$ and
$\{0\}$. The coalgebra $C$ is called \emph{cosemisimple} if it is a
coproduct in $\Vect_k$ of cosimple coalgebras. The coalgebra $C$ is
called \emph{split cosemisimple} if it is cosemisimple and every
cosimple sub coalgebra is a matrix coalgebra.
\end{definition}

We prefer the term \emph{cosemisimple} rather than the more common
\emph{semisimple} because it indicates that this is a property of a
coalgebra. In the following, \emph{semisimple WHA} therefore means
that the underlying algebra of the WHA is semisimple whereas
\emph{cosemisimple WHA} means that its underlying coalgebra has the
property just defined above.

\begin{definition}
Let $(C,\Delta,\epsilon)$ be a coalgebra over a field $k$. A right
$C$-comodule $(V,\beta_V)$ is called \emph{irreducible} if
$V\neq\{0\}$ and $V$ has no sub comodules other than $V$ and $\{0\}$.
\end{definition}

We here use the term \emph{irreducible} as opposed to \emph{simple} in
order to distinguish it from the property that an object $X$ of a
$k$-linear category satisfies $\End(X)\cong k$.

\begin{lemma}
\label{lem_corep}
Let $(C,\Delta,\epsilon)$ be a coalgebra over a field $k$.
\begin{myenumerate}
\item
Every irreducible right $C$-comodule is finite-dimensional as a
vector space over $k$.
\item
If $V$ and $W$ are irreducible right $C$-comodules and $V\not\cong
W$, then $\Hom(V,W)=\{0\}$.
\item
If $C$ is split cosemisimple and $V$ an irreducible right
$C$-comodule, then $\End(V)\cong k$.
\item
If $C$ is cosemisimple and $V$ a finite-dimensional right
$C$-comodule, then
\begin{equation}
\label{eq_cosemisimplebiprod}
V\cong \bigoplus_{i=1}^n V_i
\end{equation}
for some irreducible right $C$-comodules $V_i$ and $n\in\N_0$.
\end{myenumerate}
\end{lemma}

\subsection{Tannaka--Kre\v\i n reconstruction}
\label{sect_tkreview}

In this section, we summarize the main results on Tannaka--Kre\v\i n
reconstruction of a coalgebra from a category $\sym{C}$ with a functor
$\sym{C}\to\Vect_k$, following~\cite{Sc92}.

Let $\sym{C}$ be a small category and $\omega\colon\sym{C}\to\Vect_k$
be a functor taking values in $\fdVect_k$. Then the \emph{coend}
\begin{equation}
\coend(\sym{C},\omega) = \int^{X\in|\sym{C}|}{\omega(X)}^\ast\otimes\omega(X)
\end{equation}
exists.

In the following, we ignore all set theoretic issues and no longer
mention the requirement that $\sym{C}$ be small. In fact, all examples
relevant to topology and mathematical physics that we are aware of,
can already be obtained with essentially small $\sym{C}$, and whenever
a coend appears, we can therefore replace $\sym{C}$ by an equivalent
small category.

By $\Nat(\omega,\omega\otimes-)\colon\Vect_k\to\Set$ we denote the
functor that sends each vector space $M$ to the set
$\Nat(\omega,\omega\otimes M)$ of natural transformations
$\omega\Rightarrow\omega\otimes M$ and each linear map $\phi\colon
M\to N$ to the map of sets
$(\id_\omega\otimes\phi)\circ-\colon\Nat(\omega,\omega\otimes
M)\to\Nat(\omega,\omega\otimes N)$.

\begin{theorem}
\label{thm_coend}
Let $\sym{C}$ be a category and $\omega\colon\sym{C}\to\Vect_k$ be a
functor taking values in $\fdVect_k$. For a any vector space $C$, the
following are equivalent:
\begin{myenumerate}
\item
$C\cong\coend(\sym{C},\omega)$.
\item
The functor $\Nat(\omega,\omega\otimes-)\colon\Vect_k\to\Set$ is
representable with representing object $C$.
\item
There is a natural transformation
$\delta^\omega\colon\omega\Rightarrow\omega\otimes C$ such that for
each vector space $M$ and each natural transformation
$\phi\colon\omega\Rightarrow\omega\otimes M$, there is a unique
linear map $f\colon C\to M$ such that the diagram
\begin{equation}
\label{eq_universal}
\begin{aligned}
\xymatrix{
\omega\ar[rr]^{\delta^\omega}\ar[ddrr]_\phi&&
\omega\otimes C\ar[dd]^{\id_\omega\otimes f}\\
\\
&&\omega\otimes M
}
\end{aligned}
\end{equation}
of natural transformations between functors $\sym{C}\to\Vect_k$
commutes.
\end{myenumerate}
\end{theorem}

\begin{proposition}
Let $\sym{C}$ be a category and $\omega\colon\sym{C}\to\Vect_k$ be a
functor taking values in $\fdVect_k$. The vector space
$C=\coend(\sym{C},\omega)$ forms a coassociative counital coalgebra
$(C,\Delta,\epsilon)$. The operations $\Delta\colon C\to C\otimes C$
and $\epsilon\colon C\to k$ are determined from the universal property
of the coend by commutativity of the following diagrams of natural
transformations between functors $\sym{C}\to\Vect_k$:
\begin{equation}
\begin{aligned}
\xymatrix{
\omega\ar[rr]^{d^\omega}\ar[d]_{\delta^\omega}\ar[rrdd]&&
\omega\otimes C\ar[dd]^{\id_\omega\otimes\Delta}\\
\omega\otimes C\ar[d]_{\delta^\omega\otimes\id_C}\\
(\omega\otimes C)\otimes C\ar[rr]_{\alpha_{\omega(-),C,C}}&&
\omega\otimes(C\otimes C)
}
\end{aligned}
\end{equation}
and
\begin{equation}
\begin{aligned}
\xymatrix{
\omega\ar[rr]^{\delta^\omega}\ar[ddrr]_{\rho_{\omega(-)}^{-1}}&&
\omega\otimes C\ar[dd]^{\id_\omega\otimes\epsilon}\\
\\
&&\omega\otimes k
}
\end{aligned}
\end{equation}
Here, $\alpha$ and $\rho$ denote the associator and the right unit
constraint of $\Vect_k$. We always draw the diagonal in these diagrams
in order to remind the reader of~\eqref{eq_universal}.
\end{proposition}

Part~(3) of Theorem~\ref{thm_coend} thus states that the coend
$C=\coend(\sym{C},\omega)$ is the universal coalgebra that coacts on
all objects of $\sym{C}$. The coaction of $C$ on the vector space
$\omega(X)$ associated with an object $X\in|\sym{C}|$ is given by
$\delta^\omega_X\colon\omega(X)\to\omega(X)\otimes C$. The following
proposition describes $\coend(\sym{C},\omega)$ as a vector space in
terms of generators and relations.

\begin{proposition}
Let $\sym{C}$ be a category and $\omega\colon\sym{C}\to\Vect_k$ be a
functor taking values in $\fdVect_k$. The coend is the vector space,
\begin{equation}
\label{eq_coendvect}
\coend(\sym{C},\omega)\cong\Biggl(\coprod_{X\in|\sym{C}|}{\omega(X)}^\ast\otimes\omega(X)\Biggr)/N,
\end{equation}
where $\coprod$ denotes the coproduct in the category $\Vect_k$ and
\begin{equation}
N=\{\,({\omega(f)}^\ast\theta)\otimes v-\theta\otimes(\omega(f)v)\,\mid\quad
\theta\in{\omega(Y)}^\ast; v\in\omega(X); f\colon X\to Y; X,Y\in|\sym{C}|\,\}.
\end{equation}
\end{proposition}

The coalgebra structure of the coend is a quotient modulo $N$ of a
coproduct of matrix coalgebras. Let
$({\omega(X)}^\ast,\ev_{\omega(X)},\coev_{\omega(X)})$ be a left-dual
of $\omega(X)$ in $\Vect_k$. Such a left-dual exists because
$\omega(X)$ is by assumption finite-dimensional. Then the structure of
the coalgebra $\coend(\sym{C},\omega)$ is given on the homogeneous
elements of~\eqref{eq_coendvect} by
\begin{eqnarray}
\Delta\colon{\omega(X)}^\ast\otimes\omega(X)&\to&
({\omega(X)}^\ast\otimes\omega(X))\otimes({\omega(X)}^\ast\otimes\omega(X)),\nn\\
\theta\otimes v&\mapsto&\sum_j\theta\otimes e_j^{(X)}\otimes e^j_{(X)}\otimes v,\\
\epsilon\colon{\omega(X)}^\ast\otimes\omega(X)&\to&k,\nn\\
\theta\otimes v&\mapsto&\ev_{\omega(X)}(\theta\otimes v).
\end{eqnarray}
Here we have written $\coev_{\omega(X)}(1)=\sum_je_j^{(X)}\otimes
e^j_{(X)}$. The universal coaction of $\coend(\sym{C},\omega)$ on
$\omega(X)$ is given by
\begin{equation}
\label{eq_coaction}
\delta^\omega_X\colon\omega(X)\to\omega(X)\otimes({\omega(X)}^\ast\otimes\omega(X)),\quad
v\mapsto\sum_je_j^{(X)}\otimes e^j_{(X)}\otimes v.
\end{equation}
Below, we make use of this reconstruction of the coalgebra
$\coend(\sym{C},\omega)$ in the context in which $\sym{C}$ is a
modular category and $\omega$ the long forgetful functor.

In this section, we have used the fact that the category $\Vect_k$ is
small cocomplete and that the tensor product $\otimes$ preserves
colimits in both arguments.

\section{The long forgetful functor}
\label{sect_long}

Let us now define the long forgetful functor by composing the canonical
functor $\hat\omega\colon\sym{C}\to{}_R\sym{M}_R$ of~\cite{Ha99} with the
forgetful functor ${}_R\sym{M}_R\to\Vect_k$ and show that this functor
satisfies Szlach{\'a}nyi's conditions~\cite{Sz05}, \ie\ that the functor is
equipped with a \emph{separable Frobenius structure}. Before we can show this,
we need to establish some facts about modular categories and their
non-degenerate traces. In this section, unless specified otherwise,
$(\sym{C},\otimes,\one,\alpha,\lambda,\rho,{(-)}^\ast,\ev,\coev,\sigma,\nu)$
is a finitely semisimple additive ribbon category for which $k=\End(\one)$
is a field. ${\{V_j\}}_{j\in I}$ denotes a family of objects as in
Definition~\ref{def_semisimple}(3).

\subsection{Traces and convenient bases}

The traces of the ribbon category $\sym{C}$ can be used in order to
relate ${\Hom(X,Y)}^\ast$ with $\Hom(Y,X)$, $X,Y\in|\sym{C}|$. For the
reconstruction, it turns out to be convenient if one `rescales' the traces by
the following isomorphisms.

\begin{proposition}
There is a natural equivalence $D\colon 1_{\sym{C}}\Rightarrow
1_{\sym{C}}$ of the identity functor, given by
\begin{equation}
D_X\colon X\to X,\qquad D_X:=\sum_{\ell=1}^{n^X}\imath^X_\ell\circ\pi^X_\ell{(\dim V_{j^X_\ell})}^{-1}
\end{equation}
for all objects $X\in|\sym{C}|$. Here $j^X_\ell$, $\imath^X_\ell$ and
$\pi^X_\ell$ are as in Definition~\ref{def_modular}(2d) or
Definition~\ref{def_semisimple}(3c).
\end{proposition}

\begin{corollary}
\label{cor_trace}
For any two objects $X,Y\in|\sym{C}|$, the map
\begin{equation}
\phi_{X,Y}\colon\Hom(Y,X)\otimes\Hom(X,Y)\to k,\qquad
f\otimes g\mapsto\tr_X(D_X\circ f\circ g).
\end{equation}
is a non-degenerate symmetric and associative $k$-bilinear form, \ie\
it is a non-degenerate $k$-bilinear form and satisfies
\begin{myenumerate}
\item
Symmetry: $\phi_{X,Y}(f\otimes g) = \phi_{Y,X} (g\otimes f)$ for all
morphisms $f\colon Y\to X$ and $g\colon X\to Y$ of $\sym{C}$, and
\item
Associativity: $\phi_{X,Z}(f\otimes (g\circ h)) = \phi_{X,Y}((f\circ
g)\otimes h)$ for all $f\colon Z\to X$, $g\colon Y\to Z$ and
$h\colon X\to Y$.
\end{myenumerate}
\end{corollary}

\begin{proof}
This follows from the cyclic property of the trace, from
non-degeneracy of $\sym{C}$ (Definition~\ref{def_nondegenerate}) and
from the fact that $D_X$ is invertible for all $X\in|\sym{C}|$.
\end{proof}

It is then possible to write down a pair of dual bases of $\Hom(X,Y)$
and $\Hom(Y,X)$ with respect to $\phi_{X,Y}$.

\begin{proposition}
\label{prop_dualbasis}
Let $X,Y\in|\sym{C}|$. Then
\begin{eqnarray}
\{\,e_{\alpha\beta}=\imath^Y_\alpha\circ\pi^X_\beta\colon X\to Y
\mid\quad 1\leq\alpha\leq n^Y, 1\leq\beta\leq n^X, j^Y_\alpha=j^X_\beta\,\},\\
\{\,e^{\gamma\delta}=\imath^X_\delta\circ\pi^Y_\gamma\colon X\to Y
\mid\quad 1\leq\gamma\leq n^Y, 1\leq\delta\leq n^X, j^Y_\gamma=j^X_\delta\,\}
\end{eqnarray}
form a pair of dual basis of $\Hom(X,Y)$ and $\Hom(Y,X)$ with respect
to $\phi_{X,Y}$, \ie\
\begin{equation}
\phi_{X,Y} (e^{\gamma\delta}\otimes e_{\alpha\beta}) = \delta_{\alpha\gamma}\delta_{\beta\delta}.
\end{equation}
Here we have used Definition~\ref{def_modular}(2d)
(Definition~\ref{def_semisimple}(3c)) for both $X$ and $Y$.
\end{proposition}

\begin{proof}
Use Proposition~\ref{prop_modular}(4) for both $X$ and $Y$.
\end{proof}

\subsection{The long forgetful functor}

\begin{definition}
The \emph{universal object} of $\sym{C}$ is defined as
\begin{equation}
\hat V:=\bigoplus_{j\in I} V_j
\end{equation}
\end{definition}

Note that the universal object is determined up to isomorphism by the
category $\sym{C}$ and that it is determined fully as soon as a family
${\{V_j\}}_{j\in I}$ and a total order of $I$ have been fixed. We assume from
now on that such a choice has been made.

\begin{proposition}
\label{prop_verlinde}
\begin{myenumerate}
\item
The $k$-vector space $R=\End(\hat V)$ forms a commutative separable
$k$-algebra with respect to composition.
\item
A basis ${(\lambda_j)}_{j\in I}$ of orthogonal idempotents for $R$
is given by $\lambda_j(v)=0$ if $v\in V_{\ell}$, $\ell\neq j$, and
$\lambda_j(v)=v$ if $v\in V_j$.
\item
Every morphism $f\colon\hat V\to\hat V$ is of the form $f=\sum_{j\in
I}f_j\lambda_j$ with some $f_j\in k$.
\end{myenumerate}
\end{proposition}

\begin{proof}
Definition~\ref{def_modular}(2a) and Proposition~\ref{prop_modular}(4).
\end{proof}

The following definition is the canonical functor
$\hat\omega\colon\sym{C}\to{}_R\sym{M}_R$ of~\cite{Ha99} composed with
the forgetful functor ${}_R\sym{M}_R\to\Vect_k$.

\begin{definition}
The \emph{long forgetful functor}
is the functor
\begin{eqnarray}
\omega\colon\sym{C}\to\Vect_k,\quad X &\mapsto&\Hom(\hat V,\hat V\otimes X),\\
f &\mapsto& (\id_{\hat V}\otimes f)\circ-.\nn
\end{eqnarray}
\end{definition}

\noindent
Note that the long forgetful functor is $k$-linear and takes values in
$\fdVect_k$.

\begin{proposition}
\label{prop_hayashidual}
Let $\omega\colon\sym{C}\to\Vect_k$ be the long forgetful functor. Then
$\omega(X)$, $X\in|\sym{C}|$, has a left-dual
$({\omega(X)}^\ast,\ev_{\omega(X)},\coev_{\omega(X)})$ where
${\omega(X)}^\ast=\Hom(\hat V\otimes X,\hat V)$,
\begin{alignat}{2}
\ev_{\omega(X)}&\colon{\omega(X)}^\ast\otimes\omega(X)\to k,&&\quad
\theta\otimes v\mapsto\tr_{\hat V}(D_{\hat V}\circ\theta\circ v),\\
\label{eq_hayashicoev}
\coev_{\omega(X)}&\colon k\to\omega(X)\otimes{\omega(X)}^\ast,&&\quad
1\mapsto \sum_je_j^{(X)}\otimes e^j_{(X)}.
\end{alignat}
Here, ${\{e_j^{(X)}\}}_j$ and ${\{e^j_{(X)}\}}_j$ denote a pair of
dual bases of $\omega(X)=\Hom(\hat V,\hat V\otimes X)$ and
${\omega(X)}^\ast=\Hom(\hat V\otimes X,\hat V)$ with respect to
$\phi_{\hat V,\hat V\otimes X}$. Given any morphism $f\colon X\to Y$
of $\sym{C}$, the morphism dual to $\omega(f)=(\id_{\hat V}\otimes
f)\circ-$ is given by
\begin{equation}
\label{eq_hayashimor}
\omega(f)^\ast = -\circ(\id_{\hat V}\otimes f).
\end{equation}
\end{proposition}

\begin{proof}
Corollary~\ref{cor_trace} and Proposition~\ref{prop_dualbasis} imply
the triangle identities.
\end{proof}

\begin{proposition}
\label{prop_faithful}
The long forgetful functor $\omega\colon\sym{C}\to\Vect_k$ is faithful.
\end{proposition}

\begin{proof}
Let $X,Y\in|\sym{C}|$ and $f,g\colon X\to Y$ be arbitrary morphisms of
$\sym{C}$. We have to show that $\omega(f)=\omega(g)$ implies $f=g$.

Choose some arbitrary $\ell\in I$, $q\colon V_\ell\to X$ and $p\colon
Y\to V_\ell$. For every $\ell\in I$, denote by $\imath_\ell^{(\hat
V)}\colon V_\ell\to \hat V$ and $\pi_\ell^{(\hat V)}\colon\hat V\to V_\ell$
the morphisms of Definition~\ref{def_modular}(2d)
{(Definition~\ref{def_semisimple}(3c))} associated with
$V_\ell$ for $X=\hat V$. Then define $v:=(\imath_0^{(\hat V)}\otimes
q)\circ\lambda_{V_\ell}^{-1}\circ\pi_\ell^{(\hat V)}\colon \hat
V\to\hat V\otimes X$ and $\eta:=\imath_\ell^{(\hat
V)}\circ\lambda_{V_\ell}\circ(\pi_0^{(\hat V)}\otimes p)\colon\hat
V\otimes Y\to\hat V$. We compute that
\begin{equation}
0 = \eta\circ(\omega(f)-\omega(g))(v)
= \eta\circ(\id_{\hat V}\otimes (f-g))\circ v
\end{equation}
and
\begin{equation}
0 = \tr_{\hat V}(\eta\circ(\id_{\hat V}\otimes (f-g))\circ v)
= \tr_{V_\ell}(p\circ (f-g)\circ q).
\end{equation}
This holds for any $\ell\in I$ and any $p$ and $q$. If we insert all
$q=\imath_m^{(X)}$ and $p=\pi_n^{(Y)}$ with $j_m=j_n=\ell$, we
conclude that $0=f-g$.
\end{proof}

\begin{remark}
Our definition of a modular category does not assume the existence of
all finite limits (preabelian category) nor that all monomorphisms and
all epimorphisms are normal (abelian category). In
Corollary~\ref{cor_abelian} below, we nevertheless see that all
finitely semisimple additive ribbon categories with $k=\End(\one)$ a
field and therefore all modular categories are in fact
abelian and that the long forgetful functor is exact.
\end{remark}

The remainder of the present subsection can be skipped on first
reading. The results are, however, needed in several proofs below.

\begin{lemma}
\label{lem_isoanti}
Let $X\in|\sym{C}|$. Then there are natural isomorphisms
\begin{eqnarray}
\Phi_X \colon \omega(X)&\to& {\omega(X^\ast)}^\ast,\nn\\
v&\mapsto& D_{\hat V}^{-1}\circ\rho_{\hat V}\circ(\id_{\hat V}\otimes\bar\ev_X)
\circ\alpha_{\hat V,X,X^\ast}\circ (v\otimes\id_{X^\ast})\circ (D_{\hat V}\otimes\id_{X^\ast}),\\
\Psi_X \colon {\omega(X)}^\ast&\to& \omega(X^\ast),\nn\\
\theta&\mapsto& (\theta\otimes\id_{X^\ast})\circ\alpha^{-1}_{\hat V,X,X^\ast}\circ(\id_{\hat V}\otimes\coev_X)
\circ\rho^{-1}_{\hat V}.
\end{eqnarray}
Their composites are given by
\begin{eqnarray}
\Xi_X:=\Psi_{X^\ast}\circ\Phi_X \colon\omega(X)&\to&\omega({X^\ast}^\ast),\nn\\
v&\mapsto& (D_{\hat V}^{-1}\otimes\tau_X)\circ v\circ D_{\hat V},\\
\Theta_X:=\Phi_{X^\ast}\circ\Psi_X \colon{\omega(X)}^\ast&\to&{\omega({X^\ast}^\ast)}^\ast,\nn\\
\theta&\mapsto& D_{\hat V}^{-1}\circ\theta\circ(D_{\hat V}\otimes\tau_X^{-1}).
\end{eqnarray}
Here $\tau_X\colon X\to{X^\ast}^\ast$ denotes the isomorphism
of~\eqref{eq_pivotal}.
\end{lemma}

\begin{proof}
Naturality follows from the properties of dual morphisms. The
morphisms $\Phi_X$ and $\Psi_X$ are invertible because $D_{\hat V}$ is
and because of the triangle identities~\eqref{eq_zigzag1}
to~\eqref{eq_zigzag4}. In order to determine their composites, one
needs~\eqref{eq_pivotal}, \eqref{eq_rightdual1}
and~\eqref{eq_rightdual2}.
\end{proof}

\begin{proposition}
Let $X\in|\sym{C}|$, and let
${\{e^{(X)}_j\}}_j$ and ${\{e^j_{(X)}\}}_j$ form a pair of dual bases
of $\omega(X)$ and ${\omega(X)}^\ast$ with respect to $\phi_{\hat
V,\hat V\otimes X}$. Then
\begin{equation}
\label{eq_complete1}
\sum_j e^{(X)}_j\circ e^j_{(X)} = \id_{\hat V\otimes X}
\end{equation}
and
\begin{equation}
\label{eq_complete2}
\sum_j \Psi(e^j_{(X)})\circ \Phi(e^{(X)}_j) = \id_{\hat V\otimes X^\ast}.
\end{equation}
\end{proposition}

\begin{proof}
Show~\eqref{eq_complete1} first for the pair of dual bases of
Proposition~\ref{prop_dualbasis} with $\hat V$ instead of $X$ and
$\hat V\otimes X$ instead of $Y$. Then the claim follows for any other
pair of dual bases. In order to verify~\eqref{eq_complete2}, show that
${(\Psi(e^j_{(X)}))}_j$ and ${(\Phi(e^{(X)}_j))}_j$ form a pair of
dual bases of $\omega(X^\ast)$ and ${\omega(X^\ast)}^\ast$,
respectively. Then the claim follows from~\eqref{eq_complete1}.
\end{proof}

\begin{proposition}
\label{prop_coevprod}
Let $X,Y\in|\sym{C}|$, and let
${\{e^{(X)}_j\}}_j$ and ${\{e^j_{(X)}\}}_j$ as well as
${\{e^{(Y)}_\ell\}}_\ell$ and ${\{e^\ell_{(Y)}\}}_\ell$ be pairs of
dual bases of $\omega(X)$ and ${\omega(X)}^\ast$ as well as of
$\omega(Y)$ and ${\omega(Y)}^\ast$, respectively. Define
\begin{eqnarray}
e_{j\ell}^{(X\otimes Y)}
&:=& \alpha_{\hat V,X,Y}\circ(e_j^{(X)}\otimes\id_Y)\circ e_\ell^{(Y)}\in\omega(X\otimes Y),\\
e^{j\ell}_{(X\otimes Y)}
&:=& e^\ell_{(Y)}\circ(e^j_{(X)}\otimes\id_Y)\circ\alpha_{\hat V,X,Y}^{-1}\in{\omega(X\otimes Y)}^\ast.
\end{eqnarray}
Then $({\omega(X\otimes Y)}^\ast,\ev_{\omega(X\otimes
Y)},\coev_{\omega(X\otimes Y)})$ is a left-dual of $\omega(X\otimes
Y)$ with
\begin{alignat}{3}
\ev_{\omega(X\otimes Y)}&\colon{\omega(X\otimes Y)}^\ast\otimes\omega(X\otimes Y)\to k,&&\qquad
\theta\otimes v\mapsto \phi_{\hat V,\hat V\otimes(X\otimes Y)}(\theta\otimes v),\\
\label{eq_coev}
\coev_{\omega(X\otimes Y)}&\colon k\to\omega(X\otimes Y)\otimes{\omega(X\otimes Y)}^\ast,&&\qquad
1\mapsto\sum_{j,\ell}e_{j\ell}^{(X\otimes Y)}\otimes e^{j\ell}_{(X\otimes Y)}.
\end{alignat}
\end{proposition}

\begin{proof}
The triangle identities for $\ev_{\omega(X\otimes Y)}$ and
$\coev_{\omega(X\otimes Y)}$ follow from the triangle identities for
$\ev_{\omega(X)}$ and $\coev_{\omega(X)}$ as well as for $\ev_{\omega(Y)}$
and $\coev_{\omega(Y)}$, \cf\ Proposition~\ref{prop_hayashidual}, and
from~\eqref{eq_complete1}. Note that the $e^{(X\otimes Y)}_{j\ell}$ are in
general linearly dependent. The sum in~\eqref{eq_coev} nevertheless yields a
perfectly acceptable coevaluation map.
\end{proof}

\subsection{Functors with separable Frobenius structure}

We can now show that the long forgetful functor
$\omega\colon\sym{C}\to\Vect_k$ satisfies the following conditions due
to Szlach{\'a}nyi~\cite{Sz05}.

\begin{definition}
\label{def_specialfrob}
Let $\sym{D}$ and $\sym{D}^\prime$ be monoidal categories. A
\emph{functor with separable Frobenius structure}
$(F,F_{X,Y},F_0,F^{X,Y},F^0)\colon\sym{D}\to\sym{D}^\prime$ is a
functor $F\colon\sym{D}\to\sym{D}^\prime$ which is lax monoidal as
$(F,F_{X,Y},F_0)$ and oplax monoidal as $(F,F^{X,Y},F^0)$ (\cf\
Definition~\ref{def_lax}) and which satisfies the following
compatibility conditions.
\begin{equation}
\label{eq_splitepic}
F_{X,Y}\circ F^{X,Y} = \id_{F(X\otimes Y)},
\end{equation}
\begin{equation}
\label{eq_specialfrob1}
\begin{aligned}
\xymatrix{
F(X\otimes Y)\otimes^\prime FZ\ar[rr]^{F_{X\otimes Y,Z}}\ar[dd]_{F^{X,Y}\otimes^\prime\id_{FZ}}&&
F((X\otimes Y)\otimes Z)\ar[rr]^{F\alpha_{X,Y,Z}}&&
F(X\otimes(Y\otimes Z))\ar[dd]^{F^{X,Y\otimes Z}}\\
\\
(FX\otimes^\prime FY)\otimes^\prime FZ\ar[rr]_{\alpha^\prime_{FX,FY,FZ}}&&
FX\otimes^\prime(FY\otimes^\prime FZ)\ar[rr]_{\id_{FX}\otimes^\prime F_{Y,Z}}&&
FX\otimes^\prime F(Y\otimes Z),
}
\end{aligned}
\end{equation}
\begin{equation}
\label{eq_specialfrob2}
\begin{aligned}
\xymatrix{
FX\otimes^\prime F(Y\otimes Z)\ar[rr]^{F_{X,Y\otimes Z}}\ar[dd]_{\id_{FX}\otimes^\prime F^{Y,Z}}&&
F(X\otimes (Y\otimes Z))\ar[rr]^{F\alpha^{-1}_{X,Y,Z}}&&
F((X\otimes Y)\otimes Z)\ar[dd]^{F^{X\otimes Y,Z}}\\
\\
FX\otimes^\prime(FY\otimes^\prime FZ)\ar[rr]_{{\alpha^\prime}^{-1}_{FX,FY,FZ}}&&
(FX\otimes^\prime FY)\otimes^\prime FZ\ar[rr]_{F_{X,Y}\otimes^\prime\id_{FZ}}&&
F(X\otimes Y)\otimes^\prime FZ,
}
\end{aligned}
\end{equation}
for all $X,Y,Z\in|\sym{D}|$.
\end{definition}

The reason for choosing the term \emph{Frobenius structure} becomes
obvious if one visualizes the compatibility conditions by the
following diagrams. For more details on these diagrams, we refer
to~\cite{La05,LP1,LP2}.
\begin{equation}
\begin{aligned}
\begin{pspicture}(4,5)
\rput(2,2){\luturn}
\rput(2,2){\ruturn}
\rput(2,2){\llturn}
\rput(2,2){\rlturn}
\rput(2,4.5){$\scriptstyle X\otimes Y$}
\rput(2,-0.5){$\scriptstyle X\otimes Y$}
\end{pspicture}
\end{aligned}
=
\begin{aligned}
\begin{pspicture}(2.5,5)
\rput(1,0){\straight}
\rput(1.5,0){\straight}
\rput(1,2){\straight}
\rput(1.5,2){\straight}
\rput(1.5,4.5){$\scriptstyle X\otimes Y$}
\rput(1.5,-0.5){$\scriptstyle X\otimes Y$}
\end{pspicture}
\end{aligned}
\qquad
\begin{aligned}
\begin{pspicture}(3.5,5)
\rput(2,4){\llturn}
\rput(2,4){\rlturn}
\rput(2.5,4){\rlcurve}
\rput(2,0){\lucurve}
\rput(2.5,0){\luturn}
\rput(2.5,0){\ruturn}
\rput(1.25,4.5){$\scriptstyle X$}
\rput(3,4.5){$\scriptstyle Y\otimes Z$}
\rput(1.5,-0.5){$\scriptstyle X\otimes Y$}
\rput(3.25,-0.5){$\scriptstyle Z$}
\end{pspicture}
\end{aligned}
=
\begin{aligned}
\begin{pspicture}(3.5,5)
\rput(1,2){\llturn}
\rput(1,2){\rlturn}
\rput(0,2){\straight}
\rput(2.5,2){\luturn}
\rput(2.5,2){\ruturn}
\rput(3,0){\straight}
\rput(0.25,4.5){$\scriptstyle X$}
\rput(2.5,4.5){$\scriptstyle Y\otimes Z$}
\rput(1,-0.5){$\scriptstyle X\otimes Y$}
\rput(3.25,-0.5){$\scriptstyle Z$}
\end{pspicture}
\end{aligned}
\qquad
\begin{aligned}
\begin{pspicture}(3.5,5)
\rput(2,4){\llcurve}
\rput(2.5,4){\llturn}
\rput(2.5,4){\rlturn}
\rput(2,0){\luturn}
\rput(2,0){\ruturn}
\rput(2.5,0){\rucurve}
\rput(1.5,4.5){$\scriptstyle X\otimes Y$}
\rput(3.25,4.5){$\scriptstyle Z$}
\rput(1.25,-0.5){$\scriptstyle X$}
\rput(3,-0.5){$\scriptstyle Y\otimes Z$}
\end{pspicture}
\end{aligned}
=
\begin{aligned}
\begin{pspicture}(3.5,5)
\rput(2.5,2){\llturn}
\rput(2.5,2){\rlturn}
\rput(3,2){\straight}
\rput(1,2){\luturn}
\rput(1,2){\ruturn}
\rput(0,0){\straight}
\rput(1,4.5){$\scriptstyle X\otimes Y$}
\rput(3.25,4.5){$\scriptstyle Z$}
\rput(0.25,-0.5){$\scriptstyle X$}
\rput(2.5,-0.5){$\scriptstyle Y\otimes Z$}
\end{pspicture}
\end{aligned}
\end{equation}

\begin{theorem}
\label{thm_specialfrob}
The long forgetful functor $\omega\colon\sym{C}\to\Vect_k$ has a separable
Frobenius structure
$(\omega,\omega_{X,Y},\omega_0,\omega^{X,Y},\omega^0)$ with
\begin{alignat}{2}
\label{eq_longfrob1}
\omega_{X,Y}&\colon\omega(X)\otimes\omega(Y)\to\omega(X\otimes Y),&&\quad
f\otimes g\mapsto\alpha_{\hat V,X,Y}\circ(f\otimes\id_Y)\circ g,\\
\omega_0&\colon k\to\omega(\one),&&\quad
1\mapsto\rho_{\hat V}^{-1},
\end{alignat}
and
\begin{alignat}{2}
\omega^{X,Y}&\colon\omega(X\otimes Y)\to\omega(X)\otimes\omega(Y),&&\quad
h\mapsto\sum_{j,\ell}\ev_{\omega(X\otimes Y)}(e^{j\ell}_{(X\otimes Y)}\otimes h)\,e_j^{(X)}\otimes e_{\ell}^{(Y)},\\
\omega^0&\colon\omega(\one)\to k,&&\quad
v\mapsto \ev_{\omega(\one)}(\rho_{\hat V}\otimes v),
\end{alignat}
using the $e^{j\ell}_{(X\otimes Y)}$ of
Proposition~\ref{prop_coevprod}.
\end{theorem}

\begin{proof}
We need to verify the following.
\begin{myenumerate}
\item
$\omega$ is indeed a functor.
\item
$\omega_{X,Y}$ and $\omega^{X,Y}$ are natural transformations.
\item
$(\omega,\omega_{X,Y},\omega_0)$ is lax monoidal. The hexagon axiom
for the lax monoidal functor follows from the pentagon axiom in
$\sym{C}$, and the two squares follow from the triangle axiom.
\item
$(\omega,\omega^{X,Y},\omega^0)$ is oplax monoidal. Again, the
hexagon axiom follows from the pentagon, and the two squares from
the triangle.
\item
In order to show the compatibility
conditions~\eqref{eq_specialfrob1} and~\eqref{eq_specialfrob2}, we
verify for each $h\in\omega(X\otimes Y)$ and $w\in\omega(Z)$ that
\begin{eqnarray}
\omega^{X,Y\otimes Z}\circ\omega(\alpha_{X,Y,Z})\circ\omega_{X\otimes Y,Z} (h\otimes w)\nn\\
= (\id_{\omega(X)}\otimes\omega_{Y,Z})\circ\alpha_{\omega(X),\omega(Y),\omega(Z)}
\circ(\omega^{X,Y}\otimes\id_{\omega(Z)})(h\otimes w),
\end{eqnarray}
which follows from the definitions, using the left-duals of
Proposition~\ref{prop_hayashidual}, the basis of
Proposition~\ref{prop_coevprod} and the pentagon axiom of
$\sym{C}$. Similarly, for $v\in\omega(X)$ and $f\in\omega(Y\otimes
Z)$, we verify that
\begin{eqnarray}
\omega^{X\otimes Y,Z}\circ\omega(\alpha^{-1}_{X,Y,Z})\circ\omega_{X,Y\otimes Z}(v\otimes f)\nn\\
= (\omega_{X,Y}\otimes\id_{\omega(Z)})\circ\alpha^{-1}_{\omega(X),\omega(Y),\omega(Z)}
\circ(\id_{\omega(X)}\otimes\omega^{Y,Z}) (v\otimes f).
\end{eqnarray}
The condition~\eqref{eq_splitepic} follows from
Proposition~\ref{prop_coevprod}.
\end{myenumerate}
\end{proof}

\begin{remark}
\label{rem_splitmonic}
Under the conditions of Theorem~\ref{thm_specialfrob},
\begin{equation}
\label{eq_splitmonic}
\omega^0\circ\omega_0 = |I|,
\end{equation}
where $|I|$ is the number of (isomorphism classes of) simple objects as in
Definition~\ref{def_modular}(2) or Definition~\ref{def_semisimple}(3).

If the characteristic of $k$ divides $|I|$, then this is
zero. Otherwise, $\omega_0$ is a split monomorphism with left-inverse
$\omega^0/|I|$. Recall that because of~\eqref{eq_splitepic}, the
$\omega_{X,Y}$ are split epimorphisms with right-inverse
$\omega^{X,Y}$. If the characteristic of $k$ does not divide $|I|$,
one may call the Frobenius structure \emph{special},
\cf~\cite[Definition 2.3]{FuSc03}. It is interesting to note that the
question of whether the right hand side of~\eqref{eq_splitmonic} is
zero or not, does not play any role in the following.
\end{remark}

\section{Tannaka--Kre\v\i n reconstruction}
\label{sect_reconstruct}

In this section, unless specified otherwise, $\sym{C}$ denotes a finitely
semisimple additive ribbon category for which $k=\End(\one)$ is a field.

\subsection{Coalgebra structure}

If $\omega\colon\sym{C}\to\Vect_k$ is the long forgetful functor, then the
coend $H=\coend(\sym{C},\omega)$ has the structure of a coassociative counital
coalgebra $(H,\Delta,\epsilon)$ as in Section~\ref{sect_tkreview}.

For the homogeneous elements of the coend~\eqref{eq_coendvect}, we write
${[\theta|v]}_X\in{\omega(X)}^\ast\otimes\omega(X)$ with
$\theta\in{\omega(X)}^\ast$ and $v\in\omega(X)$, \ie\ $\theta\colon\hat
V\otimes X\to\hat V$ and $v\colon\hat V\to\hat V\otimes X$, using the
left-duals of Proposition~\ref{prop_hayashidual}. The relations in the
quotient~\eqref{eq_coendvect} then read,
\begin{equation}
\label{eq_relations}
{[\zeta\circ(\id_{\hat V}\otimes f)|v]}_X
= {[\zeta|(\id_{\hat V}\otimes f)\circ v]}_Y,
\end{equation}
where $v\colon\hat V\to\hat V\otimes X$, $f\colon X\to Y$ and
$\zeta\colon\hat V\otimes Y\to\hat V$. The coalgebra operations of $H$
can be written as
\begin{eqnarray}
\label{eq_coalg1}
\Delta  ({[\theta|v]}_X) &=& \sum_j{[\theta|e^{(X)}_j]}_X\otimes{[e_{(X)}^j|v]}_X,\\
\label{eq_coalg2}
\epsilon({[\theta|v]}_X) &=& \ev_{\omega(X)}(\theta\otimes v),
\end{eqnarray}
and the universal coaction as
\begin{equation}
\delta^\omega_X(v) = \sum_je_j^{(X)}\otimes{[e^j_{(X)}|v]}_X
\end{equation}
for all $v\colon\hat V\to\hat V\otimes X$ and $\theta\colon\hat
V\otimes X\to\hat V$, $X\in|\sym{C}|$.

\subsection{Semisimplicity}

With the explicit description of the coend of Section~\ref{sect_tkreview}, we
can show that the coend $H=\coend(\sym{C},\omega)$ is a
finite-dimensional split cosemisimple coalgebra.

\begin{proposition}
Let $\sym{D}$ be an $\mathbf{Ab}$-enriched ribbon category,
$k=\End(\one)$ and $\omega\colon\sym{D}\to\Vect_k$ be a $k$-linear
functor taking values in $\fdVect_k$.
\begin{myenumerate}
\item
If $\sym{D}$ is semisimple with a family ${\{V_j\}}_j$, $j\in I$, of
simple objects as in Definition~\ref{def_semisimple}(3), then
\begin{equation}
\coend(\sym{D},\omega)\cong\coprod_{j\in I}{\omega(V_j)}^\ast\otimes\omega(V_j).
\end{equation}
With the operations~\eqref{eq_coalg1} and~\eqref{eq_coalg2}, the
coend therefore forms a split cosemisimple coalgebra.
\item
If $\sym{D}$ is finitely semisimple, then the coend is
finite-dimensional, \ie\
\begin{equation}
\label{eq_coendsum}
\coend(\sym{D},\omega)\cong\bigoplus_{j\in I}{\omega(V_j)}^\ast\otimes\omega(V_j),
\end{equation}
where the coproduct has turned into a finite biproduct.
\end{myenumerate}
\end{proposition}

\begin{proof}
We show that the composition of the inclusion
\begin{equation}
\bigoplus_{j\in I}{\omega(V_j)}^\ast\otimes\omega(V_j)\to
\coprod_{X\in|\sym{D}|}{\omega(X)}^\ast\otimes\omega(X)
\end{equation}
with the canonical projection
\begin{equation}
\coprod_{X\in|\sym{D}|}{\omega(V_j)}^\ast\otimes\omega(V_j)\to
\biggl(\coprod_{X\in|\sym{D}|}{\omega(X)}^\ast\otimes\omega(X)\biggr)/N
\end{equation}
is a bijection. In order to see that it is surjective, we use
Definition~\ref{def_modular}(2d) or Definition~\ref{def_semisimple}(3c),
and in order to see that it is injective, Proposition~\ref{prop_modular}(5).
\end{proof}

\subsection{Algebra structure}

In this section, we use the monoidal structure of $\sym{C}$ in order to
equip the coend with the structure of an associative unital algebra. In the
remainder of Section~\ref{sect_reconstruct}, all commutative diagrams are in
$\Vect_k$. We write $k$ for the monoidal unit object of $\Vect_k$.

In order to reconstruct an algebra structure on
$H=\coend(\sym{C},\omega)$ from the monoidal structure of $\sym{C}$,
we consider the category $\sym{C}$ with the functor
$\omega\otimes\omega\colon\sym{C}\to\Vect_k$, $X\mapsto
\omega(X)\otimes\omega(X)$, $f\mapsto\omega(f)\otimes\omega(f)$. The
corresponding coend and the universal coaction are given as follows.

\begin{proposition}[see, for example~\cite{Sc92}]
Let $\sym{D}$ be a monoidal category, $\omega\colon\sym{D}\to\Vect_k$
be a functor taking values in $\fdVect_k$ and
$H=\coend(\sym{D},\omega)$.
\begin{myenumerate}
\item
The coend of $\omega\otimes\omega\colon\sym{D}\times\sym{D}\to\Vect_k$ is the
tensor product coalgebra,
\begin{equation}
H\otimes H\cong\coend(\sym{D}\times\sym{D},\omega\otimes\omega),
\end{equation}
with the operations $\Delta_{H\otimes
H}=(\id_H\otimes\sigma_{H,H}\otimes\id_H)\circ(\Delta_H\otimes\Delta_H)$
and $\epsilon_{H\otimes H}=\epsilon_H\otimes\epsilon_H$.
\item
The corresponding universal coaction is given by
$\delta^{\omega\otimes\omega}_{X,Y}\colon\omega(X)\otimes\omega(Y)\to(\omega(X)\otimes\omega(Y))\otimes
(H\otimes H)$ where
\begin{equation}
\delta^{\omega\otimes\omega}_{X,Y}=(\id_{\omega(X)}\otimes\sigma_{H,\omega(Y)}\otimes\id_{H})
\circ(\delta^\omega_X\otimes\delta^\omega_Y).
\end{equation}
\end{myenumerate}
\end{proposition}

In addition, let $\mathbf{1}$ denote the monoidal category whose only morphism
is the identity of the monoidal unit. Then $\mathbf{1}$ with the functor
$\omega^{\otimes 0}\colon\mathbf{1}\to\Vect_k$, $\ast\mapsto k$,
$\id_\ast\mapsto\id_k$, has the trivial coalgebra as the coend,
$\coend(\mathbf{1},\omega^{\otimes 0})\cong k$, and the universal coaction
$\delta^{\omega^{\otimes 0}}\colon k\to k\otimes k$, $\delta^{\omega^{\otimes
0}}=\rho_k^{-1}$.

\begin{theorem}
\label{thm_algebra}
Let $\omega\colon\sym{C}\to\Vect_k$ be the long forgetful functor. Then
the coend $H=\coend(\sym{C},\omega)$ is equipped with the structure
$(H,\mu,\eta)$ of an associative unital algebra. Its operations are
determined from the universal property of the coend by commutativity
of
\begin{equation}
\label{eq_alguniv1}
\begin{aligned}
\xymatrix{
\omega(X)\otimes\omega(Y)\ar[rr]^{\delta^{\omega\otimes\omega}_{X,Y}}\ar[d]_{\omega_{X,Y}}
\ar[ddrr]&&
(\omega(X)\otimes\omega(Y))\otimes(H\otimes H)\ar[dd]^{\id_{\omega(X)\otimes\omega(Y)}\otimes\mu}\\
\omega(X\otimes Y)\ar[d]_{\delta^\omega_{X\otimes Y}}\\
\omega(X\otimes Y)\otimes H\ar[rr]_{\omega^{X,Y}\otimes\id_H}&&
(\omega(X)\otimes\omega(Y))\otimes H
}
\end{aligned}
\end{equation}
and of
\begin{equation}
\label{eq_alguniv2}
\begin{aligned}
\xymatrix{
k\ar[rr]^{\delta^{\omega^{\otimes 0}}}\ar[d]_{\omega_0}\ar[ddrr]&&
k\otimes k\ar[dd]^{\id_k\otimes\eta}\\
\omega(\one)\ar[d]_{\delta^\omega_{\one}}\\
\omega(\one)\otimes H\ar[rr]_{\omega^0\otimes\id_H}&&
k\otimes H.
}
\end{aligned}
\end{equation}
\end{theorem}

In order to prove the theorem, it is convenient to compute the
operations $\mu$ and $\eta$ in a basis of $H=\coend(\sym{C},\omega)$
that is adapted to the matrix coalgebra structure~\eqref{eq_coendsum}
as follows.

\begin{lemma}
\label{lem_algebra}
Under the conditions of Theorem~\ref{thm_algebra}, the operations are
given by
\begin{eqnarray}
\label{eq_alg1}
\mu({[\theta|v]}_X\otimes{[\zeta|w]}_Y)
&=& {[\zeta\circ(\theta\otimes\id_Y)\circ\alpha^{-1}_{\hat V,X,Y}|
\alpha_{\hat V,X,Y}\circ(v\otimes\id_Y)\circ w]}_{X\otimes Y},\\
\label{eq_alg2}
\eta(1) &=& {[\rho_{\hat V}|\rho_{\hat V}^{-1}]}_{\one}
\end{eqnarray}
for $\theta\in{\omega(X)}^\ast$, $\zeta\in{\omega(Y)}^\ast$,
$v\in\omega(X)$ and $w\in\omega(Y)$.
\end{lemma}

\begin{proof}
In order to show that the operations~\eqref{eq_alg1}
and~\eqref{eq_alg2} make the diagrams~\eqref{eq_alguniv1}
and~\eqref{eq_alguniv2} commute, one needs the definitions and
Proposition~\ref{prop_coevprod}.
\end{proof}

\begin{proof}[Proof of Theorem~\ref{thm_algebra}]
We use the operations $\mu$ and $\eta$ as given in
Lemma~\ref{lem_algebra}. In order to show that $\mu$ is well defined
on the quotient modulo $N$ of~\eqref{eq_coendvect}, one needs the
relations~\eqref{eq_relations}. For associativity of $\mu$, one needs
the pentagon axiom for the associator of $\sym{C}$, and for the unit
laws the triangle axiom.
\end{proof}

In this section, we have not only used that $\Vect_k$ is small
cocomplete and that $\otimes$ preserves colimits in both arguments,
but also that $\Vect_k$ has a symmetric braiding.

\subsection{Weak Hopf Algebra structure}

In this section, we show that the algebra and coalgebra structure of
the coend satisfy the compatibility conditions of a WBA, and we use the
left-duals in $\sym{C}$ in order to construct an antipode that turns
it into a WHA.

\begin{theorem}
Let $\omega\colon\sym{C}\to\Vect_k$ be the long forgetful functor. Then
the coend $H=\coend(\sym{C},\omega)$ has the structure of a WBA
$(H,\mu,\eta,\Delta,\epsilon)$.
\end{theorem}

\begin{proof}
In order to verify the conditions of Definition~\ref{def_wba}, we
express the operations of $H$ in the form~\eqref{eq_alg1},
\eqref{eq_alg2}, \eqref{eq_coalg1} and~\eqref{eq_coalg2}. In order to
verify~\eqref{eq_wba1}, we use
Proposition~\ref{prop_coevprod}. For~\eqref{eq_wba2}
and~\eqref{eq_wba3}, we need the triangle equations for the evaluation
and coevaluation maps of Proposition~\ref{prop_hayashidual} as well as
the associativity property of Corollary~\ref{cor_trace} for the
traces involved.
\end{proof}

In order to reconstruct an antipode from duality in $\sym{C}$, we
consider the opposite category $\sym{C}^\op$ with the functor
$\omega^\ast\colon\sym{C}^\op\to\Vect_k$, $X\mapsto{\omega(X)}^\ast$,
$f\mapsto{\omega(f)}^\ast$. The corresponding coend and the universal
coaction on the duals are characterized by the following result.

\begin{proposition}[see, for example~\cite{Sc92}]
Let $\sym{D}$ be a left-autonomous monoidal category,
$\omega\colon\sym{D}\to\Vect_k$ be a functor taking values in
$\fdVect_k$ and $H=\coend(\sym{D},\omega)$.
\begin{myenumerate}
\item
The coend of $\omega^\ast\colon\sym{D}^\op\to\Vect_k$ is the
coopposite coalgebra,
\begin{equation}
H^\cop\cong\coend(\sym{D}^\op,\omega^\ast).
\end{equation}
\item
The corresponding universal coaction is given by
$\delta^{\omega^\ast}_X\colon{\omega(X)}^\ast\to{\omega(X)}^\ast\otimes
H^\cop$ where
\begin{eqnarray}
\delta^{\omega^\ast}_X
&=& \sigma_{H,{\omega(X)}^\ast}\circ(\lambda_H\otimes\id_{{\omega(X)}^\ast})
\circ((\ev_{\omega(X)}\otimes\id_H)\otimes\id_{{\omega(X)}^\ast})\nn\\
&&  \circ(\alpha^{-1}_{{\omega(X)}^\ast,\omega(X),H}\otimes\id_{{\omega(X)}^\ast})
\circ((\id_{{\omega(X)}^\ast}\otimes\delta^\omega_X)\otimes\id_{{\omega(X)}^\ast})\nn\\
&&  \circ\alpha^{-1}_{{\omega(X)}^\ast,\omega(X),{\omega(X)}^\ast}
\circ(\id_{{\omega(X)}^\ast}\otimes\coev_{\omega(X)})\circ\rho^{-1}_{{\omega(X)}^\ast}.
\end{eqnarray}
\end{myenumerate}
\end{proposition}

\begin{theorem}
\label{thm_wha}
Let $\omega\colon\sym{C}\to\Vect_k$ be the long forgetful functor. Then
the coend $H=\coend(\sym{C},\omega)$ has the structure of a WHA
$(H,\mu,\eta,\Delta,\epsilon,S)$. Here the antipode $S\colon H\to H$
is determined from the universal property of the coend by
commutativity of
\begin{equation}
\label{eq_antipodeuniv}
\begin{aligned}
\xymatrix{
{\omega(X)}^\ast\ar[rr]^{\delta^{\omega^\ast}_X}\ar[d]_{\Psi_X}\ar[ddrr]&&
{\omega(X)}^\ast\otimes H\ar[dd]^{\id_{{\omega(X)}^\ast}\otimes S}\\
\omega(X^\ast)\ar[d]_{\delta^\omega_{X^\ast}}\\
\omega(X^\ast)\otimes H\ar[rr]_{\Psi_X^{-1}\otimes\id_H}&&{\omega(X)}^\ast\otimes H.
}
\end{aligned}
\end{equation}
with $\Psi_X$ as in Lemma~\ref{lem_isoanti}.
\end{theorem}

\noindent
Again, we first express the antipode in our preferred basis.
\begin{lemma}
\label{lem_antipode}
Under the conditions of Theorem~\ref{thm_wha}, the antipode is given by
\begin{equation}
\label{eq_antipodebasis}
S({[\theta|v]}_X) = {[\Phi_X(v)|\Psi_X(\theta)]}_{X^\ast}
\end{equation}
with $\Phi$ and $\Psi$ as defined in Lemma~\ref{lem_isoanti}.
\end{lemma}

\begin{proof}
We verify in a direct computation that the linear map $S$
of~\eqref{eq_antipodebasis} makes the diagram~\eqref{eq_antipodeuniv}
commute. The top right of that diagram can be computed from the
definitions using the left-duals of Proposition~\ref{prop_hayashidual}
whereas the bottom left can be obtained from the explicit expression
for $\Psi$ in Lemma~\ref{lem_isoanti}.
\end{proof}

\begin{proof}[Proof of Theorem~\ref{thm_wha}]
Using the expression~\eqref{eq_antipodebasis}, we can employ the
relations~\eqref{eq_relations} to show that $S$ is well defined on the
quotient~\eqref{eq_coendvect}. Before we verify~\eqref{eq_wha1}
and~\eqref{eq_wha2}, we first compute $\epsilon_t$ and $\epsilon_s$,
\begin{eqnarray}
\epsilon_t({[\theta|v]}_X)
&=& {[\Phi_X(v)\circ \Psi_X(\theta)\circ\rho_{\hat V}|\rho_{\hat V}^{-1}]}_{\one},\\
\epsilon_s({[\theta|v]}_X)
&=& {[\rho_{\hat V}|\rho_{\hat V}^{-1}\circ\theta\circ v]}_{\one},
\end{eqnarray}
for all $\theta\in{\omega(X)}^\ast$ and $v\in\omega(X)$. In order to
verify~\eqref{eq_wha1}, one needs~\eqref{eq_complete1}, and
for~\eqref{eq_wha2}, one needs the triangle and pentagon axioms in
$\sym{C}$ as well as~\eqref{eq_complete2}. Finally, \eqref{eq_wha3}
can be obtained from~\eqref{eq_wha2} and~\eqref{eq_complete1}.
\end{proof}

\begin{proposition}
Let $\omega\colon\sym{C}\to\Vect_k$
be the long forgetful functor and $H=\coend(\sym{C},\omega)$ the
reconstructed WHA. Then
\begin{equation}
\label{eq_ssquare}
S^2({[\theta|v]}_X)
= {[D_{\hat V}^{-1}\circ\theta\circ(D_{\hat V}\otimes\id_X)|(D_{\hat V}^{-1}\otimes\id_X)\circ v\circ D_{\hat V}]}_X,
\end{equation}
for all $v\in\omega(X)$ and $\theta\in{\omega(X)}^\ast$,
$X\in|\sym{C}|$.
\end{proposition}

\begin{proof}
Using Lemma~\ref{lem_antipode} and Lemma~\ref{lem_isoanti}, we get
\begin{equation}
S^2({[\theta|v]}_X) = {[\Theta_X(\theta)|\Xi_X(v)]}_{{X^\ast}^\ast}
\end{equation}
which implies the claim upon using the relations~\eqref{eq_relations}.
\end{proof}

\begin{remark}
In the reconstructed WHA, we have
\begin{equation}
\epsilon\circ\eta = \omega^0\circ\omega_0 = |I|\in k,
\end{equation}
\cf\ Remark~\ref{rem_splitmonic}.
\end{remark}

\subsection{Coribbon structure}

In this section, we define the notion of a coribbon WHA and show that the WHA
reconstructed from $\sym{C}$ has this structure.

\begin{definition}
Let $(H,\mu,\eta,\Delta,\epsilon,S)$ be a WHA. A linear form $f\colon
H\to k$ is called
\begin{myenumerate}
\item
\emph{convolution invertible} if there exists some linear $\bar
f\colon H\to k$ such that $f(x^\prime)\bar
f(x^\pprime)=\epsilon(x)=\bar f(x^\prime)f(x^\pprime)$ for all $x\in H$,
\item
\emph{dual central} if $f(x^\prime)x^\pprime = x^\prime f(x^\pprime)$
for all $x\in H$,
\item
\emph{dual group-like} if it is convolution invertible and
\begin{equation}
\label{eq_dualgrouplike}
f(xy)=\epsilon(x^\prime y^\prime)f(x^\pprime)f(y^\pprime) =
f(x^\prime)f(y^\prime)\epsilon(x^\pprime y^\pprime)
\end{equation}
for all $x,y\in H$.
\end{myenumerate}
\end{definition}

Note that $\bar f$ in~(1) is uniquely determined by $f$. Every dual
group-like linear form also satisfies
$f(\epsilon_t(x))=\epsilon(x)=f(\epsilon_s(x))$ and $f(S(x))=\bar
f(x)$ for all $x\in H$.

\begin{definition}
\label{def_coquasi}
A \emph{coquasitriangular} WHA $(H,\mu,\eta,\Delta,\epsilon,S,r)$ over
a field $k$ is a WHA $(H,\mu,\eta,\Delta,\epsilon,S)$ over $k$ with a
linear form $r\colon H\otimes H\to k$ (\emph{universal $r$-form}) that
satisfies the following conditions:
\begin{myenumerate}
\item
For all $x,y\in H$,
\begin{equation}
\label{eq_coquasidef}
r(x\otimes y)=\epsilon(x^\prime y^\prime)r(x^\pprime\otimes y^\pprime)
=r(x^\prime\otimes y^\prime)\epsilon(y^\pprime x^\pprime)
\end{equation}
\item
There exists some linear $\bar r\colon H\otimes H\to k$ such that
\begin{eqnarray}
\label{eq_coquasiinv1}
\bar r(x^\prime\otimes y^\prime)r(x^\pprime\otimes y^\pprime)&=&\epsilon(yx)\\
\label{eq_coquasiinv2}
r(x^\prime\otimes y^\prime)\bar r(x^\pprime\otimes y^\pprime)&=&\epsilon(xy)
\end{eqnarray}
\item
For all $x,y,z\in H$,
\begin{eqnarray}
\label{eq_almostcomm}
x^\prime y^\prime r(x^\pprime\otimes y^\pprime)
&=&r(x^\prime\otimes y^\prime)y^\pprime x^\pprime\\
\label{eq_coquasitensor1}
r((xy)\otimes z)&=&r(y\otimes z^\prime) r(x\otimes z^\pprime)\\
\label{eq_coquasitensor2}
r(x\otimes (yz))&=&r(x^\prime\otimes y) r(x^\pprime\otimes z)
\end{eqnarray}
\end{myenumerate}
The WHA $H$ is called \emph{cotriangular} if in addition
\begin{equation}
r(x^\prime\otimes y^\prime) r(y^\pprime\otimes x^\pprime) = \epsilon(xy)
\end{equation}
for all $x,y\in H$.
\end{definition}

Note that $\bar r$ in~(2) is uniquely determined by $r$ if one
imposes~\eqref{eq_coquasidef}, \eqref{eq_coquasiinv1}
and~\eqref{eq_coquasiinv2}. Condition~\eqref{eq_coquasidef} says that
$r$ is well defined on the tensor product $H\hotimes H$ and on its
opposite if $H$ is viewed as the right-regular $H$-comodule. The
conditions~\eqref{eq_coquasiinv1} and~\eqref{eq_coquasiinv2} express
\emph{weak convolution invertibility}, \eqref{eq_almostcomm}
\emph{almost commutativity} and \eqref{eq_coquasitensor1}
and~\eqref{eq_coquasitensor2} compatibility with the tensor product.

\begin{theorem}
\label{thm_coquasi}
Let $\omega\colon\sym{C}\to\Vect_k$ be the long forgetful functor. Then
$H=\coend(\sym{C},\omega)$ is a coquasitriangular WHA
$(H,\mu,\eta,\Delta,\epsilon,S,r)$ whose universal $r$-form $r\colon
H\otimes H\to k$ is determined from the universal property of the
coend by commutativity of
\begin{equation}
\label{eq_coquasiuniv1}
\begin{aligned}
\xymatrix{
\omega(X)\otimes\omega(Y)\ar[rr]^{\delta^{\omega\otimes\omega}_{X,Y}}
\ar[dddrr]\ar[d]_{\sigma_{\omega(X),\omega(Y)}}&&
(\omega(X)\otimes\omega(Y))\otimes(H\otimes H)\ar[ddd]^{\id_{\omega(X)}\otimes\id_{\omega(Y)}\otimes r}\\
\omega(Y)\otimes\omega(X)\ar[d]_{\omega_{X,Y}}\\
\omega(Y\otimes X)\ar[d]_{\omega(\sigma_{Y,X})}\\
\omega(X\otimes Y)\ar[r]_{\omega^{X,Y}}&
\omega(X)\otimes\omega(Y)\ar[r]_{\rho^{-1}_{\omega(X)\otimes\omega(Y)}}&
(\omega(X)\otimes\omega(Y))\otimes k
}
\end{aligned}
\end{equation}
for all $X,Y\in|\sym{C}|$. Here, $\sigma_{Y,X}$ denotes the braiding of
$\sym{C}$, $\sigma_{\omega(X),\omega(Y)}$ the braiding of $\Vect_k$
and $\rho_{\omega(X)\otimes\omega(Y)}$ the right unit constraint of
$\Vect_k$. Its weak convolution inverse $\bar r\colon H\otimes H\to k$
is determined by commutativity of
\begin{equation}
\label{eq_coquasiuniv2}
\begin{aligned}
\xymatrix{
\omega(X)\otimes\omega(Y)\ar[rr]^{\delta^{\omega\otimes\omega}_{X,Y}}
\ar[dddrr]\ar[d]_{\omega_{X,Y}}&&
(\omega(X)\otimes\omega(Y))\otimes(H\otimes H)\ar[ddd]^{\id_{\omega(X)}\otimes\id_{\omega(Y)}\otimes\bar r}\\
\omega(X\otimes Y)\ar[d]_{\omega(\sigma^{-1}_{Y,X})}\\
\omega(Y\otimes X)\ar[d]_{\omega^{Y,X}}\\
\omega(Y)\otimes\omega(X)\ar[r]_{\sigma_{\omega(Y),\omega(X)}}&
\omega(X)\otimes\omega(Y)\ar[r]_{\rho^{-1}_{\omega(X)\otimes\omega(Y)}}&
(\omega(X)\otimes\omega(Y))\otimes k
}
\end{aligned}
\end{equation}
for all $X,Y\in|\sym{C}|$.
\end{theorem}

\begin{lemma}
Under the conditions of Theorem~\ref{thm_coquasi}, the universal
$r$-form and its weak convolution inverse are given by
\begin{eqnarray}
\label{eq_rbasis}
&&r({[\theta|v]}_X\otimes{[\zeta|w]}_Y)\\
&=& \ev_{\omega(X\otimes Y)}\bigl((\zeta\circ(\theta\otimes\id_Y)\circ\alpha^{-1}_{\hat V,X,Y})\otimes
((\id_{\hat V}\otimes\sigma_{Y,X})\circ\alpha_{\hat V,Y,X}\circ(w\otimes\id_X)\circ v)\bigr),\nn\\
\label{eq_rbarbasis}
&&\bar r({[\theta|v]}_X\otimes{[\zeta|w]}_Y)\\
&=& \ev_{\omega(X\otimes Y)}\bigl(
(\theta\circ(\zeta\otimes\id_X)\circ\alpha^{-1}_{\hat V,Y,X}\circ(\id_{\hat V}\otimes\sigma^{-1}_{Y,X}))\otimes
(\alpha_{\hat V,X,Y}\circ(v\otimes\id_Y)\circ w)\bigr).\nn
\end{eqnarray}
\end{lemma}

\begin{proof}
We have to show that~\eqref{eq_rbasis} and~\eqref{eq_rbarbasis} make
the diagrams~\eqref{eq_coquasiuniv1} and~\eqref{eq_coquasiuniv2}
commute. In order to prove this, one needs the definitions of the
morphisms that appear in these diagrams as well as
Proposition~\ref{prop_coevprod}.
\end{proof}

\begin{proof}[Proof of Theorem~\ref{thm_coquasi}]
We verify in a direct computation that $r$ and $\bar r$
of~\eqref{eq_rbasis} and~\eqref{eq_rbarbasis} satisfy the conditions
of Definition~\ref{def_coquasi}. In order to
show~\eqref{eq_coquasidef}, \eqref{eq_coquasiinv1},
\eqref{eq_coquasiinv2} and~\eqref{eq_almostcomm}, one needs
Proposition~\ref{prop_coevprod} and the
relations~\eqref{eq_relations}. In order to
show~\eqref{eq_coquasitensor1} and~\eqref{eq_coquasitensor2}, one
needs in addition the pentagon axiom of the associator of $\sym{C}$,
the hexagon axioms for the braiding of $\sym{C}$, and the cyclic
property of the trace involved in $\ev_{\omega(-)}$.
\end{proof}

\begin{definition}
\label{def_drinfeld}
Let $(H,\mu,\eta,\Delta,\epsilon,S,r)$ be a coquasitriangular
WHA. Then we define
\begin{myenumerate}
\item
the linear form
\begin{equation}
\label{eq_qform}
q\colon H\otimes H\to k,\quad x\otimes y\mapsto r(x^\prime\otimes y^\prime)r(y^\pprime\otimes x^\pprime),
\end{equation}
\item
the linear forms (\emph{dual Drinfel'd elements})
\begin{alignat}{2}
u&\colon H\to k,&&\quad
x\mapsto r(S(x^\pprime)\otimes x^\prime),\\
\label{eq_dualdrinfeldv}
v&\colon H\to k,&&\quad
x\mapsto r(S(x^\prime)\otimes x^\pprime).
\end{alignat}
\end{myenumerate}
\end{definition}

\begin{proposition}
Under the assumptions of Theorem~\ref{thm_coquasi}, the linear form
$q\colon H\otimes H\to k$ is the unique linear map making the diagram
\begin{equation}
\label{eq_quniv}
\begin{aligned}
\xymatrix{
\omega(X)\otimes\omega(Y)\ar[rr]^{\delta^{\omega\otimes\omega}_{X,Y}}
\ar[d]_{\omega_{X,Y}}\ar[ddrr]&&
(\omega(X)\otimes\omega(Y))\otimes(H\otimes H)\ar[dd]^{\id_{\omega(X)}\otimes\id_{\omega(Y)}\otimes q}\\
\omega(X\otimes Y)\ar[d]_{Q_{X,Y}}\\
\omega(X\otimes Y)\ar[r]_{\omega^{X,Y}}&
\omega(X)\otimes\omega(Y)\ar[r]_{\rho^{-1}_{\omega(X)\otimes\omega(Y)}}&
(\omega(X)\otimes\omega(Y))\otimes k
}
\end{aligned}
\end{equation}
commute. Here $Q_{X,Y}=\sigma_{Y,X}\circ\sigma_{X,Y}\colon X\otimes
Y\to X\otimes Y$, $X,Y\in|\sym{C}|$.
\end{proposition}

\begin{proof}
From~\eqref{eq_rbasis}, one can calculate using
Proposition~\ref{prop_coevprod} that
\begin{eqnarray}
&&q({[\theta|v]}_X\otimes{[\zeta|w]}_Y)\\
&=&\ev_{\omega(X\otimes Y)}\bigl((\zeta\circ(\theta\otimes\id_Y)\circ\alpha^{-1}_{\hat V,X,Y})\otimes
((\id_{\hat V}\otimes Q_{X,Y})\circ\alpha_{\hat V,X,Y}\circ(v\otimes\id_Y)\circ w)\bigr).\nn
\end{eqnarray}
Then commutativity of~\eqref{eq_quniv} can be verified in a direct
computation.
\end{proof}

\begin{definition}
\label{def_coribbon}
A \emph{coribbon} WHA $(H,\mu,\eta,\Delta,\epsilon,S,r,\nu)$ over a
field $k$ is a coquasitriangular WHA
$(H,\mu,\eta,\Delta,\epsilon,S,r)$ over $k$ with a convolution
invertible and dual central linear form $\nu\colon H\to k$
(\emph{universal ribbon twist}) that satisfies the following
conditions:
\begin{eqnarray}
\label{eq_coribbon1}
\nu(xy)   &=& \nu(x^\prime)\nu(y^\prime)r(x^\pprime\otimes y^\pprime)r(y^\ppprime\otimes x^\ppprime),\\
\label{eq_coribbon2}
\nu(S(x)) &=& \nu(x),
\end{eqnarray}
for all $x,y\in H$.
\end{definition}

\begin{theorem}
\label{thm_coribbon}
Let $\omega\colon\sym{C}\to\Vect_k$ be the long forgetful functor. Then
$H=\coend(\sym{C},\omega)$ is a coribbon WHA
$(H,\mu,\eta,\Delta,\epsilon,S,r,\nu)$ where $\nu\colon H\to k$ is
determined from the universal property of the coend by commutativity
of
\begin{equation}
\label{eq_coribbonuniv1}
\begin{aligned}
\xymatrix{
\omega(X)\ar[rr]^{\delta^\omega_X}\ar[dd]_{\omega(\nu_X)}\ar[ddrr]&&
\omega(X)\otimes H\ar[dd]^{\id_{\omega(X)}\otimes\nu}\\
\\
\omega(X)\ar[rr]_{\rho^{-1}_{\omega(X)}}&&
\omega(X)\otimes k
}
\end{aligned}
\end{equation}
and its convolution inverse $\bar\nu\colon H\to k$ by commutativity of
\begin{equation}
\label{eq_coribbonuniv2}
\begin{aligned}
\xymatrix{
\omega(X)\ar[rr]^{\delta^\omega_X}\ar[dd]_{\omega(\nu_X^{-1})}\ar[ddrr]&&
\omega(X)\otimes H\ar[dd]^{\id_{\omega(X)}\otimes\bar\nu}\\
\\
\omega(X)\ar[rr]_{\rho^{-1}_{\omega(X)}}&&
\omega(X)\otimes k
}
\end{aligned}
\end{equation}
for all $X\in|\sym{C}|$.
\end{theorem}

\begin{lemma}
Under the conditions of Theorem~\ref{thm_coribbon}, the universal
ribbon twist $\nu$ and its convolution inverse $\bar\nu$ are given by
\begin{eqnarray}
\label{eq_nubasis}
\nu({[\theta|v}]_X)
&=& \ev_{\omega(X)}\bigl(\theta\otimes((\id_{\hat V}\otimes\nu_X)\circ v)\bigr),\\
\bar\nu({[\theta|v]}_X)
\label{eq_nubarbasis}
&=& \ev_{\omega(X)}\bigl(\theta\otimes((\id_{\hat V}\otimes\nu^{-1}_X)\circ v)\bigr),
\end{eqnarray}
for all $v\in\omega(X)$, $\theta\in{\omega(X)}^\ast$, $X\in|\sym{C}|$.
\end{lemma}

\begin{proof}
We have to show that the expressions~\eqref{eq_nubasis}
and~\eqref{eq_nubarbasis} make the diagrams~\eqref{eq_coribbonuniv1}
and~\eqref{eq_coribbonuniv2} commute. This follows immediately from
the definitions.
\end{proof}

\begin{proof}[Proof of Theorem~\ref{thm_coribbon}]
We verify in a direct computation that~\eqref{eq_nubasis}
and~\eqref{eq_nubarbasis} satisfy the conditions of
Definition~\ref{def_coribbon}. In order to see that $\nu$ and
$\bar\nu$ are convolution inverse to each other, one just needs the
dual bases of Proposition~\ref{prop_hayashidual}. Showing that $\nu$
is dual central requires in addition the
relations~\eqref{eq_relations}. Verification of~\eqref{eq_coribbon1}
requires all these and Proposition~\ref{prop_coevprod} as well as the
condition~\eqref{eq_ribbon1}. Finally, in order to
verify~\eqref{eq_coribbon2}, we use~\eqref{eq_antipodebasis}, the left
autonomous structure of $\sym{C}$, the condition~\eqref{eq_ribbon2} as
well as the cyclic property of the trace involved in
$\ev_{\omega(-)}$.
\end{proof}

\subsection{Special properties of the coend}

\begin{proposition}
\label{prop_computehsht}
Let $\omega\colon\sym{C}\to\Vect_k$
be the long forgetful functor and $H=\coend(\sym{C},\omega)$ be the
reconstructed WHA. Then
\begin{eqnarray}
\label{eq_htequal}
H_t &=& \{\,{[\theta|\rho_{\hat V}^{-1}]}_{\one}\mid\quad
\theta\colon\hat V\otimes\one\to\hat V\,\},\\
\label{eq_hsequal}
H_s &=& \{\,{[\rho_{\hat V}|v]}_{\one}\mid\quad
v\colon\hat V\to\hat V\otimes\one\,\},
\end{eqnarray}
and
\begin{equation}
\label{eq_intersect}
H_t\cap H_s\cong k.
\end{equation}
\end{proposition}

\begin{proof}
Let us show~\eqref{eq_htequal}. For each $\theta\in{\omega(X)}^\ast$
and $v\in\omega(X)$, $X\in|\sym{C}|$, we have
$\epsilon_t({[\theta|v]}_X)={[\Phi_X(v)\circ\Psi_X(\theta)\circ\rho_{\hat
V}|\rho_{\hat V}^{-1}]}_{\one}$, \ie\ $H_t$ is included in the set
given. Conversely, for each $\theta\in{\omega(\one)}^\ast$,
$\epsilon_t({[\theta|\rho_{\hat V}^{-1}]}_{\one}) =
{[\Phi_{\one}(\rho_{\hat
V}^{-1})\circ\Psi_{\one}(\theta)\circ\rho_{\hat V}|\rho_{\hat
V}^{-1}]}_{\one} = {[\theta|\rho_{\hat V}^{-1}]}_{\one}$, using the
triangle axiom of $\sym{C}$, and so the given set in contained in
$H_t$. Let us show~\eqref{eq_hsequal}. For each
$\theta\in{\omega(X)}^\ast$ and $v\in\omega(X)$, $X\in|\sym{C}|$, we
have $\epsilon_s({[\theta|v]}_X) ={[\rho_{\hat V}|\rho_{\hat
V}^{-1}\circ\eta\circ v]}_{\one}$, \ie\ $H_s$ is included in the set
given. Conversely, for each $v\in\omega(\one)$, indeed
$\epsilon_s({[\rho_{\hat V}|v]}_{\one})={[\rho_{\hat V}|v]}_{\one}$,
\ie\ the given set in contained in $H_s$. Finally,
\eqref{eq_intersect} follows from~\eqref{eq_htequal}
and~\eqref{eq_hsequal}.
\end{proof}

There is one more condition that is satisfied by every WHA
$H=\coend(\sym{C},\omega)$ for a modular category $\sym{C}$ and the
long forgetful functor $\omega$. This condition is the invertibility
of the $S$-matrix. Since this condition is preserved by equivalences
of semisimple $k$-linear ribbon categories, we take care of this in
Section~\ref{sect_equiv} where we show that $\sym{C}$ is equivalent to
the category $\sym{M}^H$ of finite-dimensional right $H$-comodules.

Most proofs in this section were done by (1) defining the structure
maps in terms of a universal property; (2) expressing these maps in
terms of a convenient basis; (3) verifying their properties using this
basis. Alternatively, it would have been possible to establish their
properties directly from their defining commutative diagrams. In order
to write down or even sketch these proofs, however, one needs extra
large paper, and so we have reverted to the first method involving a
basis.

\section{Corepresentation theory}
\label{sect_corep}

In this section, we consider the category $\sym{M}^H$ of
finite-dimensional right $H$-comodules of some WBA $H$ and show that
it has the structure of a monoidal category. If $H$ is a WHA, then
each object has a specified left-dual, \ie\ $\sym{M}^H$ is
left-autonomous. If $H$ is coquasitriangular, then $\sym{M}^H$ is
braided, and if $H$ is coribbon, then $\sym{M}^H$ is ribbon.

For easier reference, we have collected all definitions relevant to
monoidal categories in Appendix~\ref{app_monoidal}, to left-autonomous
categories in Appendix~\ref{app_dual}, to ribbon categories in
Appendix~\ref{app_ribbon}, and to additive and abelian categories in
Appendix~\ref{app_abelian}.

\subsection{Preparation}

The proofs of the propositions in the section on corepresentations are
all elementary although some of them are rather laborious. They rely
on the following facts about WBAs, WHAs, and coquasitriangular or
ribbon WHAs that we collect this subsection. Some of them are quite
challenging to verify.

\begin{lemma}
Let $(H,\mu,\eta,\Delta,\epsilon)$ be a WBA. Then
\begin{eqnarray}
\label{eq_522f}
\epsilon_s(1^\prime)\otimes 1^\pprime &=& 1^\prime\otimes 1^\pprime,\\
\label{eq_522k}
\epsilon_s(xy) &=& \epsilon_s(\epsilon_s(x)y),\\
\label{eq_522m}
h^\prime\otimes h^\pprime &=& \epsilon_s(h^\prime)\otimes h^\pprime,\\
\label{eq_522p}
x\epsilon_t(y) &=& \epsilon(x^\prime y)x^\pprime,\\
\label{eq_522q}
\epsilon_s(x)y &=& y^\prime\epsilon(xy^\pprime),\\
\label{eq_526c}
1^\prime\otimes(1^\pprime h) &=& h^\prime\otimes h^\pprime,\\
\label{eq_526a}
x^\prime\epsilon_s(x^\pprime) &=& x,\\
\label{eq_526o}
x^\prime\epsilon(x^\pprime h) &=& xh,\\
\label{eq_526p}
{(\epsilon_s(x^\pprime))}^\prime\otimes\bigl(x^\prime{(\epsilon_s(x^\pprime))}^\pprime\bigr)
&=& \epsilon_s(x^\prime)\otimes x^\pprime,\\
\label{eq_526e}
\epsilon(xy^\prime)\epsilon_s(y^\pprime) &=& \epsilon_s(xy),\\
\label{eq_527b}
{(\ell h)}^\prime\otimes{(\ell h)}^\pprime &=& \ell 1^\prime\otimes h 1^\pprime,\\
\label{eq_527c}
{(\ell h)}^\prime\otimes{(\ell h)}^\pprime &=& 1^\prime \ell\otimes 1^\pprime h,
\end{eqnarray}
for all $x,y\in H$ and $h\in H_s$, $\ell\in H_t$. If
$(H,\mu,\eta,\Delta,\epsilon,S)$ is a WHA, then
\begin{eqnarray}
\label{eq_526q}
\epsilon_s(x^\pprime)\otimes(S(x^\prime)x^\ppprime)
&=& {(\epsilon_s(x))}^\prime\otimes{(\epsilon_s(x))}^\pprime,\\
\label{eq_526r}
\epsilon(hy^\pprime)y^\prime S(y^\ppprime)
&=& h\epsilon_t(y),
\end{eqnarray}
for all $x,y\in H$ and $h\in H_s$.
\end{lemma}

\begin{lemma}
Let $H$ be a WBA and $V$ be a right $H$-comodule. Then
\begin{equation}
\label{eq_569}
\epsilon(h{(v_V)}_H){(v_V)}_V\otimes\epsilon_s(V_H)
= \epsilon(hv_H){(v_V)}_V\otimes\epsilon_s({(v_V)}_H)
\end{equation}
for all $h\in H_s$ and $v\in V$.
\end{lemma}

\begin{proof}
From the comodule axioms, \eqref{eq_526e} and~\eqref{eq_522q}.
\end{proof}

\begin{lemma}
Let $(H,\mu,\eta,\Delta,\epsilon,S,r)$ be a coquasitriangular WHA. Then
\begin{equation}
\label{eq_5228g}
r(x^\prime\otimes y^\prime)\epsilon_t(y^\pprime)\epsilon_s(x^\pprime)
= \epsilon_t(x^\prime)\epsilon_s(y^\prime)r(x^\pprime\otimes y^\pprime)
\end{equation}
for all $x,y\in H$.
\end{lemma}

\begin{proof}
From the axioms, \eqref{eq_527b} and~\eqref{eq_527c}.
\end{proof}

\begin{lemma}
\label{lem_pivotal}
Let $(H,\mu,\eta,\Delta,\epsilon,S,r,\nu)$ be a coribbon WHA and $v$ be its
dual Drinfel'd element~\eqref{eq_dualdrinfeldv}. Then the \emph{pivotal form}
$w\colon H\to k$, $x\mapsto v(x^\prime)\nu(x^\pprime)$ is dual group-like and
satisfies
\begin{equation}
S^2(x) = \bar w(x^\prime)x^\pprime w(x^\ppprime)
\end{equation}
for all $x\in H$.
\end{lemma}

\begin{proof}
Analogous to the situation in a ribbon Hopf algebra. Rather tedious.
\end{proof}

\subsection{Monoidal structure}

The definitions in this section follow Nill~\cite{Ni98}, but are here
given in a form that does not assume finite-dimensionality of $H$.
\begin{proposition}
\label{prop_hscomodule}
Let $(H,\mu,\eta,\Delta,\epsilon)$ be a WBA. Then $H_s$ forms a right
$H$-comodule with
\begin{equation}
\label{eq_hscomodule}
\beta_{H_s}\colon H_s\to H_s\otimes H,\qquad
x\mapsto x^\prime\otimes x^\pprime.
\end{equation}
\end{proposition}

\begin{proof}
By~\eqref{eq_rightcoideal}, coassociativity and the counit property.
\end{proof}

\begin{proposition}
Let $H$ be a WBA and $V,W\in|\sym{M}^H|$. Then the $k$-vector space
\begin{equation}
\label{eq_truncspace}
V\hotimes W := \{\,v\otimes w\in V\otimes W\mid\quad v\otimes w= (v_V\otimes w_W)\epsilon(v_Hw_H)\,\}
\end{equation}
forms a right $H$-comodule with
\begin{equation}
\label{eq_trunccomodule}
\beta_{V\hotimes W}\colon V\hotimes W\to (V\hotimes W)\otimes H,\quad
v\otimes w \mapsto (v_V\otimes w_W)\otimes (v_Hw_H).
\end{equation}
\end{proposition}

\begin{proof}
Consequence of the WBA axioms and of the comodule axioms.
\end{proof}

The tensor product $\hotimes$ is often called the \emph{truncated
tensor product}. We note that the $k$-linear map
\begin{equation}
\label{eq_truncidempotent}
P_{V,W}\colon V\otimes W\to V\otimes W,\quad
v\otimes w\mapsto (v_V\otimes w_W)\epsilon(v_Hw_H)
\end{equation}
forms an idempotent, and that $V\hotimes W$ is its image.

\begin{theorem}
Let $H$ be a WBA. Then the category $\sym{M}^H$ forms a $k$-linear
monoidal category $(\sym{M}^H,\hotimes,H_s,\alpha,\lambda,\rho)$ with
\begin{alignat}{3}
\label{eq_comodleftunit}
\lambda_V &\colon H_s\hotimes V\to V,&&\quad h\otimes v\mapsto v_V\epsilon(hv_H),\\
\label{eq_comodrightunit}
\rho_V    &\colon V\hotimes H_s\to V,&&\quad v\otimes h\mapsto v_V\epsilon(v_Hh),
\end{alignat}
and isomorphisms $\alpha_{U,V,W}\colon(U\hotimes V)\hotimes W\to
U\hotimes(V\hotimes W)$ induced from the associator of $\Vect_k$.
\end{theorem}

\begin{proof}
\begin{myenumerate}
\item
We have already seen that $\sym{M}^H$ is $k$-linear as a
category. Since $\hotimes$ is $k$-bilinear on morphisms, $\sym{M}^H$
is also $k$-linear as a monoidal category (\cf\
Definition~\ref{def_preadditivecat}).
\item
We claim that $\lambda_V$ and $\rho_V$ are invertible with inverses
\begin{alignat}{3}
\lambda^{-1}_V &\colon V\to H_s\hotimes V,&&\quad v\mapsto (1^\prime\otimes v_V)\epsilon(1^\pprime v_H),\\
\rho^{-1}_V    &\colon V\to V\hotimes H_s,&&\quad v\mapsto v_V\otimes\epsilon_s(v_H).
\end{alignat}
While $\lambda_V$ and $\rho_V$ are obviously well defined, we have to
verify that $\lambda_V^{-1}$ maps into $H_s\hotimes V$ which follows
from~\eqref{eq_522f} and that $\rho_V^{-1}$ maps into $V\hotimes
H_s$ which follows from~\eqref{eq_526p}. In order to verify that
$\lambda_V\circ\lambda_V^{-1}=\id_V$, one needs the coaction of $H$
on $H_s$, the axioms of the comodule $V$ and the axioms of a
WBA. For $\lambda_V^{-1}\circ\lambda_V=\id_{H_s\hotimes V}$, one
needs in addition~\eqref{eq_522q}, \eqref{eq_522m}
and~\eqref{eq_526c}; for $\rho_V\circ\rho_V^{-1}=\id_V$, one
needs~\eqref{eq_526a} and for $\rho_V^{-1}\circ\rho_V=\id_{V\hotimes
H_S}$~\eqref{eq_526o}, \eqref{eq_522k} and~\eqref{eq_522q}.
\item
Using~\eqref{eq_522m} and~\eqref{eq_522q}, one shows that
$\lambda_V$ is a morphism of right $H$-comodules and
using~\eqref{eq_526p} that $\rho_V^{-1}$ is.
\item
Naturality of $\lambda_V$ and $\rho_V$ follows from the comodule
axioms of $V$ and from the properties of a morphism of comodules.
\item
The pentagon axiom can be proven from the comodule axioms and the
axioms of a WBA, and the triangle axiom from~\eqref{eq_526o}
and~\eqref{eq_522q}.
\end{myenumerate}
\end{proof}

Finally, we show that the forgetful functor $\sym{M}^H\to\Vect_k$ has
a separable Frobenius structure
(Definition~\ref{def_specialfrob}). This result is precisely dual to
that of~\cite{Sz05}. First, we recall the image factorization of an
idempotent in an abelian category.

\begin{proposition}[see, for example~\cite{Bo94}]
\label{prop_splitting} Let $\cal{C}$ be an abelian category and
$p\colon A\to A$ be an idempotent. The image factorization of $p$
yields an object $p(A)$ (the \emph{image} of $p$), which is unique up
to isomorphism, together with morphisms $\coim p\colon A\to p(A)$ (the
\emph{coimage} map) and $\im p\colon p(A)\to A$ (the \emph{image} map)
such that the following diagram commutes:
\begin{equation} \label{eq_factorization}
\begin{aligned}
\xymatrix@!C=2.3pc{
A\ar[r]^-{\coim p}\ar[dr]_{p}&p(A)\ar[d]^{\im p}\\
&A
}
\end{aligned}
\end{equation}
Since $\cal{C}$ is abelian, the idempotent $p$ splits. The splitting
is given precisely by the two morphisms of the image factorization,
and so we have $\id_{p(A)}=\coim p\circ\im p$.
\end{proposition}

\begin{proposition}
\label{prop_forgetful}
Let $(H,\mu,\eta,\Delta,\epsilon)$ be a WBA and
$U\colon\sym{M}^H\to\Vect_k$ be the forgetful functor that assigns to
each finite-dimensional right $H$-comodule $X$ its underlying vector
space $UX$ and to each morphism of right $H$-comodules $f\colon X\to
Y$ the underlying linear map $Uf\colon UX\to UY$. Then
$(U,U_{X,Y},U_0,U^{X,Y},U^0)$ is a $k$-linear faithful functor taking
values in $\fdVect_k$ with a separable Frobenius structure where
\begin{eqnarray}
U_{X,Y}=\coim P_{X,Y}  \colon UX\otimes UY         &\to& P_{X,Y}(UX\otimes UY),\\
U_0    =\eta           \colon k                    &\to& H_s,\\
U^{X,Y}=\im P_{X,Y}    \colon P_{X,Y}(UX\otimes UY)&\to& UX\otimes UY,\\
U^0    =\epsilon|_{H_s}\colon H_s                  &\to& k.
\end{eqnarray}
Here, $P_{X,Y}$ is the idempotent~\eqref{eq_truncidempotent}, and so
its image is the truncated tensor tensor product, \ie\ the vector
space underlying the tensor product in $\sym{M}^H$,
\begin{equation}
P_{X,Y}(UX\otimes UY)=U(X\hotimes Y).
\end{equation}
Furthermore, $H_s=U\one$ is the vector space underlying the monoidal
unit (Proposition~\ref{prop_hscomodule}).
\end{proposition}

\begin{proof}
\begin{myenumerate}
\item
$U$ is $k$-linear and faithful because $UX$ and $Uf$ are just the
underlying vector space and linear map, respectively.
\item
In order to see that $(U,U_{X,Y},U_0)$ is lax monoidal
(Definition~\ref{def_lax}), we have to verify the following.
\begin{myenumerate}
\item
The hexagon axiom $U\alpha_{X,Y,Z}\circ U_{X\otimes
Y,Z}\circ(U_{X,Y}\otimes\id_{UZ}) = U_{X,Y\otimes
Z}\circ(\id_{UX}\otimes U_{Y,Z})\circ\alpha_{UX,UY,UZ}$ follows from
the definitions $U_{X,Y}(x\otimes y)=(x_X\otimes
y_Y)\epsilon(x_Hy_H)$ and $U_{X\otimes Y,Z}(x\otimes y\otimes
z)=(x_X\otimes y_Y\otimes z_Z)\epsilon((x_Hy_H)z_H)$, \etc, and from
the axioms of a WBA.
\item
The first square $U\lambda_X\circ
U_{\one,X}\circ(U_0\otimes\id_{UX})=\lambda_{UX}$ follows from the
definitions $U_0(1)=1\in H_s$; $U_{\one,X}(h\otimes
x)=(h^\prime\otimes x_X)\epsilon(h^\pprime x_H)$ for $h\in H_s$,
$x\in UX$; $U\lambda_X(h\otimes x)=x_X\epsilon(hx_H)$ for $h\in
H_s$, $x\in UX$; and from the axioms of a WBA. Recall that
$\lambda_{UX}$ on the right hand side is the unit constraint of
$\Vect_k$.
\item
The second square $U\rho_X\circ U_{X,\one}\circ(\id_{UX}\otimes
U_0)=\rho_{UX}$ follows from the definitions $U_0(1)=1\in H_s$;
$U_{X,\one}(x\otimes h)=(x_X\otimes h^\prime)\epsilon(x_Hh^\pprime)$
for $x\in UX$, $h\in H_s$; $U\rho_X(x\otimes h)=x_X\epsilon(x_Hh)$
for $x\in UX$, $h\in H_s$; and from the axioms of a WBA. Again,
$\rho_{UX}$ is the unit constraint of $\Vect_k$.
\end{myenumerate}
\item
In order to see that $(U,U^{X,Y},U^0)$ is oplax monoidal, we need to
verify:
\begin{myenumerate}
\item
The hexagon axiom~\eqref{eq_oplaxhexagon} holds because
$U^{X,Y}(x\otimes y)=x\otimes y$ for all $x\otimes y\in U(X\hotimes
Y)\subseteq UX\otimes UY$.
\item
In order to verify the first square
$\lambda_{UX}\circ(U^0\otimes\id_{UX})\circ U^{\one,X}(h\otimes
x)=U\lambda_X(h\otimes x)$ for all $h\otimes x\in U(\one\hotimes
X)$, we use the fact that $h\otimes x=(h^\prime\otimes
x_X)\epsilon(h^\pprime x_H)$ and the definitions
$U^{\one,X}(h\otimes x)=h\otimes x$, $U^0(h)=\epsilon(h)$,
$\lambda_{UX}(1\otimes x)=x$ (unit constraint of $\Vect_k$) and
$U\lambda_X(h\otimes x)=x_X\epsilon(hx_H)$.
\item
In order to verify the second square $\rho_{UX}\circ(\id_{UX}\otimes
U^0)\circ U^{X,\one}(x\otimes h)=U\rho_X(x\otimes h)$ for all
$x\otimes h\in U(X\hotimes\one)$, we use the fact that $x\otimes
h=(x_X\otimes h^\prime)\epsilon(x_Hh^\pprime)$ and the definitions
$U^{X,\one}(x\otimes h)=x\otimes h$, $U^0(h)=\epsilon(h)$,
$\rho_{UX}(x\otimes 1)=x$ and $U\rho_X(x\otimes
h)=x_X\epsilon(x_Hh)$.
\end{myenumerate}
\item
Finally, we verify the compatibility conditions of
Definition~\ref{def_specialfrob}:
\begin{myenumerate}
\item
The splitting of the idempotent $P_{X,Y}$ yields $U_{X,Y}\circ
U^{X,Y}=\id_{U(X\hotimes Y)}$.
\item
In order to show $(\id_{UX}\otimes
U_{Y,Z})\circ\alpha_{UX,UY,UZ}\circ(U^{X,Y}\otimes\id_{UZ})(x\otimes
y\otimes z)=U^{X,Y\otimes Z}\circ U\alpha_{X,Y,Z}\circ U_{X\otimes
Y,Z}(x\otimes y\otimes z)$ for all $x\otimes y\otimes z\in
U(X\hotimes Y)\otimes UZ$, we use the fact that $x\otimes y\otimes
z=(x_X\otimes y_Y\otimes z)\epsilon(x_Hy_H)$ as well as the axioms
of a WBA.
\item
The proof of
$(U_{X,Y}\otimes\id_{UZ})\circ\alpha^{-1}_{UX,UY,UZ}\circ(\id_{UX}\otimes
U^{Y,Z})=U^{X\otimes Y,Z}\circ U\alpha^{-1}_{X,Y,Z}\circ
U_{X,Y\otimes Z}$ is analogous.
\end{myenumerate}
\end{myenumerate}
\end{proof}

\subsection{Duality}

\begin{proposition}
Let $H$ be a WHA and $(V,\beta_V)$ be a finite-dimensional right
$H$-comodule. Then the dual vector space $V^\ast$ forms a right
$H$-comodule with
\begin{equation}
\beta_{V^\ast}\colon V^\ast\to V^\ast\otimes H,\qquad
\theta\mapsto (v\mapsto \theta(v_V)\otimes S(v_H)).
\end{equation}
\end{proposition}

\begin{proof}
Consequence of the WBA axioms and of the comodule axioms.
\end{proof}

\begin{theorem}
Let $H$ be a WHA. Then the category $\sym{M}^H$ is left-autonomous if
the left-dual of each $V\in|\sym{M}^H|$ is chosen as
$(V^\ast,\ev_V,\coev_V)$ where $V^\ast$ is the vector space dual to
$V$ and
\begin{alignat}{3}
\ev_V   &\colon V^\ast\hotimes V\to H_s,&&\quad \theta\otimes v\to\theta(v_V)\epsilon_s(v_H),\\
\coev_V &\colon H_s\to V\hotimes V^\ast,&&\quad x\mapsto\sum_j ({(e_j)}_V\otimes e^j)\epsilon(x{(e_j)}_H).
\end{alignat}
Here we have used the evaluation and coevaluation maps that turn
$V^\ast$ into a left-dual of $V$ in $\Vect_k$:
\begin{alignat}{3}
\ev_V^{(\Vect_k)}   &\colon V^\ast\otimes V\to k,&&\quad \theta\otimes v\mapsto \theta(v),\\
\coev_V^{(\Vect_k)} &\colon k\to V\otimes V^\ast,&&\quad 1\mapsto \sum_j e_j\otimes e^j.
\end{alignat}
\end{theorem}

\begin{proof}
While $\ev_V$ is obviously well defined, we have to show that
$\coev_V$ maps into $V\hotimes V^\ast$. This follows
from~\eqref{eq_522f} and~\eqref{eq_526r}. In order to show that
$\ev_V$ is a morphism of right $H$-comodules, one
needs~\eqref{eq_526q} and for $\coev_V$ one needs~\eqref{eq_522f},
\eqref{eq_526r}, \eqref{eq_522p} and~\eqref{eq_522m}. The triangle
identities can be proven using~\eqref{eq_569}.
\end{proof}

\subsection{Ribbon structure}

In this section, we show that if $H$ is coribbon, the category
$\sym{M}^H$ is ribbon. As soon as we give the braiding and the ribbon
twist, the proofs are straightforward.

\begin{proposition}
\label{prop_comodbraided}
Let $(H,\mu,\eta,\Delta,\epsilon,S,r)$ be a coquasitriangular
WHA. Then $\sym{M}^H$ is a braided monoidal category with braiding
$\sigma_{V,W}\colon V\hotimes W\to W\hotimes V$ given by
\begin{equation}
\sigma_{V,W}(v\otimes w)=(w_W\otimes v_V)r(w_H\otimes v_H)
\end{equation}
for all $V,W\in|\sym{M}^H|$ and $v\in V$, $w\in W$. If $H$ is
cotriangular, then $\sym{M}^H$ is symmetric monoidal.
\end{proposition}

Note that $Q_{V,W}=\sigma_{W,V}\circ\sigma_{V,W}$ can be obtained by
$Q_{V,W}(v\otimes w)=(v_V\otimes w_W)q(v_H\otimes w_H)$, \cf\
Definition~\ref{def_drinfeld}.

\begin{proposition}
\label{prop_comodribbon}
Let $(H,\mu,\eta,\Delta,\epsilon,S,r,\nu)$ be a coribbon WHA. Then
$\sym{M}^H$ is a ribbon category with ribbon twist $\nu_V\colon V\to
V$ given by
\begin{equation}
\nu_V(v) = v_V\nu(v_H)
\end{equation}
for all $V\in|\sym{M}^H|$ and $v\in V$.
\end{proposition}

\begin{remark}
Note that the forgetful functor $U\colon\sym{M}^H\to\Vect_k$ of
Proposition~\ref{prop_forgetful} is in general neither braided nor
ribbon although both $\sym{M}^H$ and $\Vect_k$ are ribbon categories.
\end{remark}

\subsection{Special properties of modular categories}

\begin{proposition}
\label{prop_semisimple}
Let $(C,\Delta,\epsilon)$ be a split cosemisimple coalgebra over a
field $k$. Then $\sym{M}^C$ is semisimple according to
Definition~\ref{def_semisimple}. If $C$ is in addition
finite-dimensional over $k$, then $\sym{M}^C$ is finitely semisimple.
\end{proposition}

\begin{proof}
Let ${\{V_j\}}_{j\in I}$ be a family of objects $V_j\in|\sym{M}^C|$
that contains one and only one representative per isomorphism class of
irreducible right $C$-comodules. We show that the conditions~(3a)
to~(3c) of Definition~\ref{def_semisimple} are satisfied.
\begin{myenumerate}
\item[(3a)]
Each $V_j$, $j\in I$, satisfies $\End(V_j)\cong k$ by
Lemma~\ref{lem_corep}(3).
\item[(3b)]
By Lemma~\ref{lem_corep}(2).
\item[(3c)]
The morphisms $\imath_\ell^{(X)}$ and $\pi_\ell^{(X)}$ are those
that define the finite biproduct~\eqref{eq_cosemisimplebiprod}.
\end{myenumerate}
\end{proof}

\begin{lemma}[see {\cite[Lemma 4.5]{Ni98}}]
\label{lem_nill}
Let $H$ be a WBA. If $H_t\cap H_s\cong k$, then $\End(H_s)\cong k$ in
$\sym{M}^H$ where $H_s$ is the monoidal unit object. In particular,
every morphism $f\colon H_s\to H_s$ in $\sym{M}^H$ is of the form
$f=c\id_{H_s}$ where $c\in k$ can be determined from the condition
$f(1)=c\,1$.
\end{lemma}

\begin{corollary}
\label{cor_copure}
Let $H$ be a split cosemisimple WBA such that $H_t\cap H_s\cong
k$. Then $\sym{M}^H$ is semisimple, and there exists a $0\in I$ such
that $H_s\cong V_0$.
\end{corollary}

\begin{proof}
By the lemma, $H_s$ is simple in $\sym{M}^H$. By
Corollary~\ref{cor_allsimplethere}, this implies that $H_s\cong V_0$
for some $0\in I$.
\end{proof}

\begin{proposition}
Let $H$ be a split cosemisimple WHA. Then $\sym{M}^H$ is semisimple,
and for each $j\in I$, there is some $j^\ast\in I$ such that
${(V_j)}^\ast\cong V_{j^\ast}$.
\end{proposition}

\begin{proof}
For $j\in I$, $V_j$ is an irreducible right $H$-comodule. We show that
every morphism $f\colon V_j^\ast\to V_j^\ast$ is of the form
$f=c\cdot\id_{V_j^\ast}$ with some $c\in k$, and so $V_j^\ast$ is
simple. By Corollary~\ref{cor_allsimplethere}, this implies that
$V_j^\ast\cong V_{j^\ast}$ for some $j^\ast\in I$. This is done as
follows. Given any $f\colon V_j^\ast\to V_j^\ast$, define $g\colon
V_j\to V_j$ by
\begin{equation}
g = \rho_{V_j}\circ(\id_{V_j}\otimes\ev_{V_j})\circ\alpha_{V_j,V_j^\ast,V_j}
\circ((\id_{V_j}\otimes f)\otimes\id_{V_j})\circ(\coev_{V_j}\otimes\id_{V_j})\circ\lambda^{-1}_{V_j}.
\end{equation}
By the triangle identities,
\begin{equation}
f = g^\ast = \lambda_{V_j^\ast}\circ(\ev_{V_j}\otimes\id_{V_j^\ast})\circ\alpha^{-1}_{V_j^\ast,V_j,V_j^\ast}
\circ(\id_{V_j^\ast}\otimes (g\otimes\id_{V_j^\ast}))\circ(\id_{V_j^\ast}\otimes\coev_{V_j})\circ\rho_{V_j^\ast}^{-1},
\end{equation}
but since $V_j$ is simple by assumption, $g=c\cdot\id_{V_j}$ for some
$c\in k$, and we can use $k$-linearity of $\sym{M}^H$ as a monoidal
category and another triangle identity in order to show that
$f=c\cdot\id_{V_j^\ast}$.
\end{proof}

\begin{lemma}
\label{lem_ribbontrace}
Let $H$ be a coribbon WHA, $V\in|\sym{M}^H|$ and $f\colon V\to
V$. Then the trace $\tr_V(f)\colon H_s\to H_s$ is given by
\begin{equation}
\label{eq_tracedetail}
\tr_V(f)(h) = \sum_{j,\ell=1}^n f_{j\ell}\epsilon_s(S({(hc_{\ell j})}^\prime))\,
w({(hc_{\ell j})}^\pprime)
\end{equation}
for all $h\in H_s$. Here $n=\dim_k(V)$ is the $k$-dimension of $V$;
the $f_{j\ell}\in k$ are the matrix elements of $f$, \ie\
$f(v_\ell)=\sum_{j=1}^nv_jf_{j\ell}$; the $c_{\ell j}\in H$ are the
coefficients of $V$, \ie\ $\beta_V(v_j)=\sum_{\ell=1}^nv_\ell\otimes
c_{\ell j}$; and $w\colon H\to k$ is the pivotal form
(Lemma~\ref{lem_pivotal}).
\end{lemma}

\begin{proof}
The trace of $f$ in the ribbon category $\sym{M}^H$ is given by
\begin{equation}
\tr_V(f) = \ev_V\circ\sigma_{V,V^\ast}\circ(\nu_V\otimes\id_{V^\ast})
\circ(f\otimes\id_{V^\ast})\circ\coev_V\colon H_s\to H_s.
\end{equation}
We insert the definitions of these morphisms and use~\eqref{eq_522q}
and~\eqref{eq_5228g}.
\end{proof}

In the following, we restrict ourselves to coribbon WHAs in which
$H_t\cap H_s\cong k$. By Lemma~\ref{lem_nill}, this implies that
$\End(H_s)\cong k$, and so all traces take values in the field
$k$.
\b
In order to express the coefficients of the $S$-matrix~\eqref{eq_smatrix} of
$\sym{M}^H$ in terms of the linear form $q\colon H\otimes H\to k$
of~\eqref{eq_qform}, we proceed as follows.

\begin{definition}
Let $H$ be a coribbon WHA over a field $k$ such that $H_t\cap H_s\cong k$ and
$V\in|\sym{M}^H|$, $n=\dim_k V$. We call
\begin{equation}
\chi_V=\sum_{j=1}^n c_{jj}\in H
\end{equation}
the \emph{dual character} of $V$ and
\begin{equation}
T_V=\sum_{j,\ell=1}^n c_{j\ell}w(c_{\ell j})\in H
\end{equation}
the \emph{dual quantum character} of $V$.
\end{definition}

Note that a dual central linear form $\alpha\colon H\to k$ defines a natural
transformation of the identity functor $f^{(\alpha)}\colon
1_{\sym{M}^H}\Rightarrow 1_{\sym{M}^H}$ via $f^{(\alpha)}_V=(\id_V\otimes
\alpha)\circ\beta_V$ for all $V\in|\sym{M}^H|$. Traces in $\sym{M}^H$ can thus be
expressed in terms of the dual quantum characters.

\begin{proposition}
\label{prop_simpletrace}
Let $H$ be a coribbon WHA over a field $k$ such that $H_t\cap H_s\cong k$,
$V\in|\sym{M}^H|$, and $\alpha\colon H\to k$ be dual central. Then
\begin{equation}
\label{eq_tracef}
\tr_V(f^{(\alpha)}_V) = c^{(\alpha)}_V\,\id_{H_s},
\end{equation}
where the element $c^{(\alpha)}_V\in k$ is determined by
\begin{equation}
\label{eq_tracef2}
\alpha(T_V^\prime)\epsilon_s(S(T_V^\pprime)) = c^{(\alpha)}_V\,\eta(1).
\end{equation}
\end{proposition}

\begin{proof}
By Lemma~\ref{lem_nill}, we can determine $c^{(\alpha)}_V$ in~\eqref{eq_tracef}
by evaluation at $\eta(1)\in H_s$. We use the formula~\eqref{eq_tracedetail}.
\end{proof}

In the special case in which $H$ is a Hopf algebra, we have $H_s\cong k$ and
$\epsilon_s=\epsilon$, and so Proposition~\ref{prop_simpletrace} reduces to
$\tr_V(f^{(\alpha)}_V)=\alpha(T_V)$ as expected. In the case in which $H$ is a
WHA and $\epsilon(\eta(1))\neq 0$, we can apply $\epsilon$
to~\eqref{eq_tracef2} and obtain
$c^{(\alpha)}_V=\alpha(T_V)/\epsilon(\eta(1))$.

In order to deal with the $S$-matrix, we need the analogue of
Proposition~\ref{prop_simpletrace} for endomorphisms of a tensor product of
modules. Note that a linear form $\gamma\colon H\otimes H\to k$ that satisfies
\begin{eqnarray}
\label{eq_lineartensor1}
x^\prime y^\prime\gamma(x^\pprime\otimes y^\pprime)
&=& \gamma(x^\prime\otimes y^\prime)x^\pprime y^\pprime,\\
\label{eq_lineartensor2}
\epsilon(x^\prime y^\prime)\gamma (x^\pprime\otimes y^\pprime)
&=& \gamma(x\otimes y),
\end{eqnarray}
for all $x,y\in H$, defines a morphisms $f^{(\gamma)}_{V,W}\in\End(V\hotimes
W)$ via
\begin{equation}
f^{(\gamma)}_{V,W}=(\id_{V\hotimes W}\otimes\gamma)
\circ(\id_V\otimes\sigma_{H,W}\otimes\id_H)\circ(\beta_V\otimes\beta_W)
\end{equation}
In particular, $q\colon H\otimes H\to k$ of~\eqref{eq_qform} satisfies these
conditions and gives rise to the morphism
$f^{(q)}_{V,W}=Q_{V,W}=\sigma_{W,V}\circ\sigma_{V,W}$ whose trace in $\sym{M}^H$ is
the $S$-matrix.

\begin{proposition}
\label{prop_tensortrace}
Let $H$ be a coribbon WHA over a field $k$ such that $H_t\cap H_s\cong k$,
$V,W\in|\sym{M}^H|$ and $\gamma\colon H\otimes H\to k$ be a linear form that
satisfies~\eqref{eq_lineartensor1} and~\eqref{eq_lineartensor2}. Then
\begin{equation}
\tr_{V\hotimes W}(f^{(\gamma)}_{V,W}) = c^{(\gamma)}_{V,W}\,\id_{H_s}
\end{equation}
where the element $c^{(\gamma)}_{V,W}\in k$ is determined by
\begin{equation}
\label{eq_tracef3}
\gamma(T_V^\prime\otimes T_W^\prime)\epsilon_s(S(T_V^\pprime T_W^\pprime))
= c^{(\gamma)}_{V,W}\,\eta(1).
\end{equation}
\end{proposition}

\begin{proof}
Analogous to Proposition~\ref{prop_simpletrace}. Apply~\eqref{eq_tracedetail}
to $V\hotimes W$ and exploit that $w\colon H\to k$ is dual group-like.
\end{proof}

In the special case in which $H$ is a Hopf algebra,
Proposition~\ref{prop_tensortrace} reduces to $\tr_{V\hotimes
W}(f^{(\gamma)}_{V,W})=\gamma(T_V\otimes T_W)$ as expected. In the case in
which $H$ is a WHA and $\epsilon(\eta(1))\neq 0$, we can apply $\epsilon$
to~\eqref{eq_tracef3} and obtain $c^{(\gamma)}_{V,W}=\gamma(T_V\otimes
T_W)/\epsilon(\eta(1))$.

\begin{definition}
Let $H$ be a coribbon WHA over a field $k$ such that $H_t\cap H_s\cong k$. Let
$T(H)=\Span_k\{\,T_V\mid\,V\in|\sym{M}^H|\,\}\subseteq H$ denote the space of
dual quantum characters of $H$. We define a linear form $\tilde q\colon
T(H)\otimes T(H)\to k$, $T_V\otimes T_W\to \tilde q_{V,W}$ where the $\tilde
q_{V,W}\in k$ are determined by
\begin{equation}
q(T_V^\prime\otimes T_W^\prime)\epsilon_s(S(T_V^\pprime T_W^\pprime))
= \tilde q_{V,W}\,\eta(1).
\end{equation}
$H$ is called \emph{weakly cofactorizable} if every linear form $\phi\colon
T(H)\to k$ is of the form $\phi(-)=\tilde q(-\otimes x)$ for some $x\in T(H)$.
\end{definition}

In the special cases in which $H$ is a Hopf algebra or in which $H$ is a WHA
with $\epsilon(\eta(1))\neq 0$, the condition of weak cofactorizability reduces
to the requirement that the bilinear form $q\colon H\otimes H\to k$ be
non-degenerate if restricted to $T(H)$. The following Corollary finally spells
out the relationship with the $S$-matrix of $\sym{M}^H$.

\begin{corollary}
\label{cor_cofactorize}
Let $H$ be a finite-dimensional, split cosemisimple, coribbon WHA over a field
$k$ such that $H_t\cap H_s\cong k$. Let $\{V_j\}_{j\in I}$ denote a set of
representatives of the isomorphism classes of simple objects of $\sym{M}^H$.
$H$ is weakly cofactorizable if and only if the matrix with coefficients
$S_{j\ell}=\tilde q_{V_j,V_\ell}$ is invertible.
\end{corollary}

\noindent
The results of the present section can be summarized as follows.

\begin{theorem}
\label{thm_comodrep}
Let $H$ be a finite-dimensional, split cosemisimple, coribbon WHA over a field
$k$ such that $H_t\cap H_s \cong k$.
\begin{myenumerate}
\item
$\sym{M}^H$ is a finitely semisimple additive ribbon category.
\item
$\sym{M}^H$ is modular if and only if $H$ is weakly cofactorizable.
\end{myenumerate}
\end{theorem}

\noindent
The following terminology is therefore appropriate.

\begin{definition}
Let $H$ be a WHA over a field $k$. $H$ is called \emph{comodular} if
$H$ is a finite-dimensional, split cosemisimple, weakly
cofactorizable, coribbon WHA such that $H_t\cap H_s\cong k$.
\end{definition}

Theorem~\ref{thm_comodrep}(2) generalizes the result of
Takeuchi~\cite[Theorem~4.6(1)]{Ta01} twofold: (1) from Hopf algebras to WHAs
and (2) by removing the assumptions on the underlying field $k$.

Note that the condition of split cosemisimplicity appears as a consequence of
requiring $\End(V_j)\cong k$ for the simple objects of a modular category,
rather than allowing $\End(V_j)$ to be a finite skew field extension over
$k$. Nevertheless, everything that has been done so far, works for any field
$k$.

In Section~\ref{sect_equiv} we show that if $\sym{C}$ satisfies all
conditions of Definition~\ref{def_modular} except maybe for (3), \ie\
non-degeneracy of the $S$ matrix, and if $H=\coend(\sym{C},\omega)$ is
the WHA reconstructed from $\sym{C}$ with respect to the long
forgetful functor $\omega$, then weak cofactorizability of $H$ is both
necessary and sufficient for the non-degeneracy of the $S$-matrix.

Compared with the sufficient conditions stated in~\cite[Lemma 8.2]{NiTu03}, we
have not only removed assumptions on the underlying field $k$ (algebraic
closure and that the characteristic of $k$ does not divide $\dim_k(H_s)$), but
our condition of weak cofactorizability is indeed weaker than the condition
(dual to) factorizability used in~\cite{NiTu03}. That condition would read in
our context as follows.

A coquasitriangular WHA $H$ is called \emph{cofactorizable} if every
linear form $\phi\colon H\to k$ that satisfies
\begin{equation}
\phi(y^\prime)\epsilon_t(y^\pprime)=\epsilon_t(y^\prime)\phi(y^\pprime)
\end{equation}
for all $y\in H$, is of the form $\phi(-) = q(-\otimes x)$ for some
$x\in H$. The condition of weak cofactorizability, in contrast,
requires non-degeneracy of $q$ only on dual quantum characters.

\subsection{Further properties}

In this section, we collect further results on the reconstructed WHA.

\begin{definition}
Let $H$ be a WBA.
\begin{myenumerate}
\item
$H$ is called \emph{copure} if the monoidal unit object $H_s$ of
$\sym{M}^H$ is irreducible.
\item
$H$ is called \emph{connected} if $Z(H)\cap H_t\cong k$.
\item
If $H$ is finite-dimensional, $H$ is called \emph{coconnected} if
$H^\ast$ is connected.
\item
$H$ is called a \emph{face algebra}~\cite{Ha93} if $H_s$ is a
commutative algebra.
\end{myenumerate}
\end{definition}

\begin{proposition}
Let $\sym{C}$ be a modular category, $\omega\colon\sym{C}\to\Vect_k$
be the long forgetful functor and $H=\coend(\sym{C},\omega)$ be the
reconstructed WHA.
\begin{myenumerate}
\item
$H_t\cong R$ are isomorphic as $k$-algebras.
\item
$H_s\cong R$ are isomorphic as $k$-algebras.
\item
$H_{\mathrm{min}}\cong R^\op\otimes R$ are isomorphic as
$k$-algebras and
\begin{equation}
H_{\mathrm{min}} = {\omega(\one)}^\ast\otimes\omega(\one) =
\{\,{[\theta|v]}_{\one}\mid\quad
v\colon\hat V\to\hat V\otimes\one; \theta\colon\hat V\otimes\one\to\hat V\,\}.
\end{equation}
\item
$H$ is regular.
\item
$H$ is a face algebra.
\item
The dual Drinfel'd elements of $H$ are given by
\begin{eqnarray}
u({[\theta|v]}_X) &=& \ev_{\omega(X)}(\theta\otimes((D_{\hat V}^{-1}\otimes\nu_X)\circ v\circ D_{\hat V})),\\
v({[\theta|v]}_X) &=& \ev_{\omega(X)}(\theta\otimes((D_{\hat V}\otimes\nu_X^{-1})\circ v\circ D_{\hat V}^{-1})),
\end{eqnarray}
for all $v\in\omega(X)$, $\theta\in{\omega(X)}^\ast$ and
$X\in|\sym{C}|$.
\item
The pivotal form $w\colon H\to k$ (Lemma~\ref{lem_pivotal}) is given
by
\begin{equation}
w({[\theta|v]}_X) = \ev_{\omega(X)}\bigl((D_{\hat V}^{-1}\circ\theta\circ(D_{\hat V}\otimes\id_X))
\otimes v\bigr),
\end{equation}
for all $v\in\omega(X)$, $\theta\in{\omega(X)}^\ast$ and
$X\in|\sym{C}|$.
\item
The dual character associated with $X\in|\sym{C}|$ is the element
\begin{equation}
\sum_j{[e^j_{(X)}|e^{(X)}_j]}_X\in H.
\end{equation}
\item
The dual quantum character associated with $X\in|\sym{C}|$ is the
element
\begin{equation}
\sum_j{[D_{\hat V}^{-1}\circ e^j_{(X)}\circ (D_{\hat V}\otimes\id_X)|e_j^{(X)}]}_X\in H.
\end{equation}
\end{myenumerate}
\end{proposition}

\begin{proof}
\begin{myenumerate}
\item
From the idempotent basis ${(\lambda_j)}_j$ of
Proposition~\ref{prop_verlinde}, one obtains a pair of dual bases
${(e^{(\one)}_j)}_j$ and ${(e^j_{(\one)})}_j$ of $\omega(\one)$ and
${\omega(\one)}^\ast$ with respect to $\ev_{\omega(\one)}$ by
\begin{equation}
e^{(\one)}_j = \rho_{\hat V}^{-1}\circ\lambda_j,
\qquad
e_{(\one)}^j = \lambda_j\circ\rho_{\hat V}.
\end{equation}
It follows from the triangle axiom of $\sym{C}$ that these from a
basis of orthogonal idempotents of $H_t$.
\item
Analogous.
\item
According to~\cite[Section 3]{Ni02}, $H_\mathrm{min}$ is generated
as an algebra by $H_t\cup H_s$. From
Proposition~\ref{prop_computehsht}, we see that
\begin{equation}
{\omega(\one)}^\ast\otimes\omega(\one) = \Span_k
\{\,{[\theta|v]}_{\one}\mid\quad v\colon\hat V\to\hat V\otimes\one;\,
\theta\colon\hat V\otimes\one\to\hat V\,\}
\end{equation}
is the vector space generated by $H_t\cup H_s$. We verify in a
direct calculation that it is an algebra with unit ${[\rho_{\hat
V}|\rho^{-1}_{\hat V}]}_{\one}$ and multiplication
\begin{equation}
\mu({[e^j_{(\one)}|e_\ell^{(\one)}]}_{\one}\otimes
{[e^m_{(\one)}|e_n^{(\one)}]}_{\one})
= \delta_{jm}\delta_{\ell n}{[e^j_{(\one)}|e_\ell^{(\one)}]}_{\one},
\end{equation}
using again the triangle axiom. This equation also shows the
isomorphism of algebras $H_\mathrm{min}\cong R\otimes R\cong
R^\op\otimes R$.
\item
Using~\eqref{eq_ssquare},
\begin{eqnarray}
S^2({[e^j_{(\one)}|e_\ell^{(\one)}]}_{\one})
&=& {[D^{-1}_{\hat V}\circ e^j_{(\one)}\circ(D_{\hat V}\otimes\id_\one)|
(D^{-1}_{\hat V}\otimes\id_\one)\circ e_\ell^{(\one)}\circ D_{\hat V}]}_{\one}\nn\\
&=& {[e^j_{(\one)}|e_\ell^{(\one)}]}_{\one},
\end{eqnarray}
for all $j,\ell\in I$ by naturality of $\rho_{\hat V}$.
\item
By~(1).
\end{myenumerate}

\noindent
The remaining claims are proven by direct computation.
\end{proof}

\begin{proposition}
Every comodular WHA $H$ is coconnected and copure.
\end{proposition}

\begin{proof}
Since $H$ is finite-dimensional and $H_t\cap H_s\cong k$, $H$ is
coconnected by~\cite[Proposition 3.11]{Ni01}. By
Corollary~\ref{cor_copure}, $H_s\cong V_0$ for some $0\in I$. But by
Proposition~\ref{prop_semisimple}, each $V_j$, $j\in I$, is an
irreducible right $H$-comodule.
\end{proof}

\section{Equivalence of categories}
\label{sect_equiv}

Let us compare the original modular category $\sym{C}$ with the
category $\sym{M}^H$ of finite-dimensional right $H$-comodules of the
reconstructed WHA $H=\coend(\sym{C},\omega)$. We show that
$\sym{C}\simeq\sym{M}^H$ are equivalent as ribbon categories. In this
section, most commutative diagrams are in $\sym{M}^H$. We highlight
this fact by putting the hat on the truncated tensor product
$\hotimes$ and by writing $H_s$ rather than $\one$ for its monoidal
unit object.

\subsection{Equivalence of monoidal categories}

In order to see that $\sym{C}\simeq\sym{M}^H$ are equivalent as
monoidal categories, we show the following.

\begin{theorem}
\label{thm_monequiv}
Let $\sym{C}$ be a modular category, $\omega\colon\sym{C}\to\Vect_k$
be the long forgetful functor and $H=\coend(\sym{C},\omega)$ be the
reconstructed WHA.
\begin{myenumerate}
\item
The long forgetful functor factors through $\sym{M}^H$, \ie\ the diagram
\begin{equation}
\begin{aligned}
\xymatrix{
\sym{C}\ar[rr]^F\ar[ddrr]_{\omega}&&\sym{M}^H\ar[dd]^U\\
\\
&&\Vect_k
}
\end{aligned}
\end{equation}
commutes. Here $U\colon\sym{M}^H\to\Vect_k$ is the forgetful
functor of Proposition~\ref{prop_forgetful}.
\item
The functor $F$ is $k$-linear, essentially surjective and fully
faithful.
\item
$(F,F_{X,Y},F_0)$ forms a strong monoidal functor with
\begin{alignat}{2}
F_{X,Y}&\colon FX\hotimes FY\to F(X\otimes Y),&&\quad
f\otimes g\mapsto \alpha_{\hat V,X,Y}\circ(f\otimes\id_Y)\circ g,\\
F_0&\colon H_s\mapsto F\one,&&\quad
{[\rho_{\hat V}|v]}_{\one}\mapsto v.
\end{alignat}
\end{myenumerate}
\end{theorem}

\begin{proof}
\begin{myenumerate}
\item
The functor $F$ sends the objects and morphisms of $\sym{C}$ to the
same vector spaces and linear maps as $\omega$ does, \ie\
$FX=\omega(X)$ and $Ff=\omega(f)$ for all $X,Y\in|\sym{C}|$ and
$f\colon X\to Y$. In order to show that $F$ is well defined as a
functor to $\sym{M}^H$, we have to verify the following:
\begin{myenumerate}
\item
For each $X\in|\sym{C}|$, $FX=\omega(X)$ forms a right
$H$-comodule, \cf~\eqref{eq_coaction}.
\item
For each morphism $f\colon X\to Y$ of $\sym{C}$,
$Ff=\omega(f)=(\id_{\hat V}\otimes f)\circ-$ is a morphism of right
$H$-comodules because
\begin{equation}
\delta^\omega_Y\circ\omega(f)[v] = (\omega(f)\otimes\id_H)\circ\delta^\omega_X[v]
\end{equation}
for all $v\in\omega(X)$. This can be verified by using the
coaction~\eqref{eq_coaction}, the relations~\eqref{eq_relations},
the form of the dual morphism~\eqref{eq_hayashimor} and the
properties of the coevaluation~\eqref{eq_hayashicoev}.
\end{myenumerate}
\item
$F$ is obviously $k$-linear since $\omega$ is. $F$ is essentially
surjective because $H$ is split cosemisimple and therefore every
finite-dimensional right $H$-comodule $M$ is of the form
\begin{equation}
\label{eq_hayashisurj}
M\cong\bigoplus_{\ell=1}^n\omega(V_{j_\ell})\cong\omega\left(\bigoplus_{\ell=1}^nV_{j_\ell}\right)
\end{equation}
for some $j_1,\ldots,j_n\in I$. Here we have used that $\sym{C}$ is
additive and thus has all finite biproducts. $F$ is faithful because
$\omega$ is (Proposition~\ref{prop_faithful}). In order to see that
$F$ is full, consider some morphism $f\colon M\to N$ of $\sym{M}^H$,
decompose both source and target as in~\eqref{eq_hayashisurj}. Since
the $\omega(V_j)$, $j\in I$, are simple in $\sym{M}^H$, morphisms
between these are either null or multiples of the identity. The
latter are of the form $(\id_{\hat
V}\otimes\id_{V_j})\circ-=\omega(\id_{V_j})$.
\item
In order to see that $F_{X,Y}$ and $F_0$ are well defined, we have
to show that these linear maps are morphisms of right $H$-comodules,
\ie\
\begin{myenumerate}
\item
For all ${[\rho_{\hat V}|v]}_{\one}\in H_s$, we have
\begin{equation}
\delta^\omega_{\one}\circ F_0 ({[\rho_{\hat V}|v]}_{\one})
= (F_0\otimes\id_H)\circ\beta_{H_s}({[\rho_{\hat V}|v]}_{\one}).
\end{equation}
In order to show this, we need the coaction~\eqref{eq_coaction} on
$F\one=\omega(\one)$ and the coaction~\eqref{eq_hscomodule} on
$H_s$. Recall that the subalgebra $H_s$ of
Proposition~\ref{prop_hsht} was computed for the reconstructed WHA
in Proposition~\ref{prop_computehsht}(2).
\item
For all $X,Y\in|\sym{C}|$, $v\in FX=\omega(X)$ and $w\in
FY=\omega(Y)$, we find
\begin{eqnarray}
&&(F_{X,Y}\otimes\id_H)\circ(\id_{\omega(X)}\otimes\id_{\omega(Y)}\otimes\mu)
\circ(\id_{\omega(X)}\otimes\sigma_{H,\omega(Y)}\otimes\id_H)\nn\\
&&\circ(\delta^\omega_X\otimes\delta^\omega_Y)[v\otimes w]\nn\\
&=& \delta^\omega_{X\otimes Y}\circ F_{X,Y}[v\otimes w],
\end{eqnarray}
where the left hand side is the comodule structure of the tensor
product $FX\hotimes FY$ from~\eqref{eq_trunccomodule}. In order to
verify the equation, we use the coactions of~\eqref{eq_coaction} and
Proposition~\ref{prop_coevprod}.
\end{myenumerate}
$F_0$ is an isomorphism with inverse
\begin{equation}
F_0^{-1}\colon F\one\to H_s,\quad
v\mapsto {[\rho_{\hat V}|v]}_{\one}.
\end{equation}
In order to see that the $F_{X,Y}$ are isomorphisms, we note that
$F_{X,Y}$ is the restriction of $\omega_{X,Y}$
of~\eqref{eq_longfrob1} to $FX\hotimes FY$, \ie\ the restriction to
the truncated tensor product of $\sym{M}^H$,
\cf~\eqref{eq_truncspace}.

We see that
$\omega^{X,Y}\circ\omega_{X,Y}\colon\omega(X)\otimes\omega(Y)\to\omega(X)\otimes\omega(Y)$
agrees with the idempotent $P_{\omega(X),\omega(Y)}$
of~\eqref{eq_truncidempotent} because for all $v\in\omega(X)$ and $w\in\omega(Y)$,
\begin{eqnarray}
\label{eq_leftinverse}
&& \omega^{X,Y}\circ\omega_{X,Y}[v\otimes w]\nn\\
&=& (\id_{\omega(X)}\otimes\id_{\omega(Y)}\otimes(\epsilon\circ\mu))
\circ(\id_{\omega(X)}\otimes\sigma_{H,\omega(Y)}\otimes\id_H)
\circ(\delta^\omega_X\otimes\delta^\omega_Y)[v\otimes w]\nn\\
&=& P_{\omega(X),\omega(Y)}(v\otimes w),
\end{eqnarray}
using the coactions~\eqref{eq_coaction} and
Proposition~\ref{prop_coevprod}.

Therefore $\omega^{X,Y}$ maps into the truncated tensor product
$FX\hotimes FZ$. It is a left-inverse of $F_{X,Y}$ because
of~\eqref{eq_leftinverse} and a right-inverse because
of~\eqref{eq_splitepic}.

Finally $F_{X,Y}$ satisfies the hexagon axiom of a strong monoidal
functor because $\omega_{X,Y}$ does. The two square axioms have to
be verified explicitly:
\begin{myenumerate}
\item
First, for all ${[\rho_{\hat V}|v]}_{\one}\in H_s$ and
$w\in\omega(X)$, we have
\begin{eqnarray}
&& F\lambda_X\circ F_{\one,X}\circ(F_0\otimes\id_{FX}) ({[\rho_{\hat V}|v]}_{\one}\otimes w)\nn\\
&=& ((\epsilon\circ\mu)\otimes\id_{FX})\circ(\id_{H_s}\otimes\sigma_{FX,H})
\circ(\id_{H_s}\otimes\delta^\omega_X)({[\rho_{\hat V}|v]}_{\one}\otimes w)\nn\\
&=& \lambda_{FX}({[\rho_{\hat V}|v]}_{\one}\otimes w),
\end{eqnarray}
using the triangle axiom in $\sym{C}$, the
relations~\eqref{eq_relations}, the left-duals of
Proposition~\ref{prop_hayashidual} and~\eqref{eq_526o}.
\item
Similarly, we verify
\begin{eqnarray}
&& F\rho_X\circ F_{X,\one}\circ(\id_{FX}\otimes F_0) (w\otimes {[\rho_{\hat V}|v]}_{\one})\nn\\
&=& (\id_{FX}\otimes(\epsilon\circ\mu))\circ(\delta^\omega_X\otimes\id_H)(w\otimes {[\rho_{\hat V}|v]}_{\one})\nn\\
&=& \rho_{FX}(w\otimes{[\rho_{\hat V}|v]}_{\one}).
\end{eqnarray}
\end{myenumerate}
\end{myenumerate}
\end{proof}

\begin{corollary}
Under the conditions of Theorem~\ref{thm_monequiv}, the categories
\begin{equation}
\sym{C}\simeq_\otimes\sym{M}^{\coend(\sym{C},\omega)}
\end{equation}
are equivalent as monoidal categories.
\end{corollary}

\begin{proof}
Since $F$ is essentially surjective and fully faithful,
by~\cite[Theorem IV.4.1]{Ma73}, $F$ is part of an adjoint equivalence
$F\dashv G$. By Proposition~\ref{prop_monadjunction}, the fact that
$F$ is strong monoidal implies that $G$ is lax monoidal and both unit
and counit of the adjoint equivalence are monoidal natural
transformations.
\end{proof}

\subsection{Equivalence of ribbon categories}

In this section, we show that the original modular category $\sym{C}$
is equivalent to the category $\sym{M}^H$, $H=\coend(\sym{C},\omega)$,
as a ribbon category.

\begin{proposition}
Let $\sym{C}$ be a modular category, $\omega\colon\sym{C}\to\Vect_k$
be the long forgetful functor and $H=\coend(\sym{C},\omega)$ the
reconstructed coquasitriangular WHA. Then the functor
$F\colon\sym{C}\to\sym{M}^H$ of Theorem~\ref{thm_monequiv} is braided.
\end{proposition}

\begin{proof}
We have to show that the condition~\eqref{eq_braidedfunctor} holds for
$F$, \ie\ that
\begin{equation}
\sigma_{FX,FY}(v\otimes w) = F_{Y,X}^{-1}\circ F(\sigma_{X,Y})\circ F_{X,Y}(v\otimes w)
\end{equation}
for all $v\in FX=\omega(X)$, $w\in FY=\omega(Y)$,
$X,Y\in|\sym{C}|$. Here, $\sigma_{FX,FY}$ is the braiding obtained in
Proposition~\ref{prop_comodbraided} from the coquasitriangular
structure of Theorem~\ref{thm_coquasi}. The claim is an immediate
consequence of the definitions.
\end{proof}

\begin{proposition}
Let $\sym{C}$ be a modular category, $\omega\colon\sym{C}\to\Vect_k$
be the long forgetful functor and $H=\coend(\sym{C},\omega)$ be the
reconstructed coribbon WHA. Then the functor
$F\colon\sym{C}\to\sym{M}^H$ of Theorem~\ref{thm_monequiv} is ribbon.
\end{proposition}

\begin{proof}
We have to show that the condition~\eqref{eq_ribbonfunctor} holds for
$F$, \ie\ that
\begin{equation}
\nu_{FX}(v) = F(\nu_X)(v)
\end{equation}
for all $v\in FX=\omega(X)$, $X\in|\sym{C}|$. Here, $\nu_{FX}$ is the
ribbon twist obtained in Proposition~\ref{prop_comodribbon} from the
coribbon structure of Theorem~\ref{thm_coribbon}. The claim follows
immediately from the definitions.
\end{proof}

\subsection{Equivalence of modular categories}

\begin{definition}
Let $\sym{C}$ and $\sym{C}^\prime$ be modular categories with the same
$k=\End(\one)=\End(\one^\prime)$. We say that $\sym{C}$ and
$\sym{C}^\prime$ are \emph{equivalent as modular categories} if there
is a functor $F\colon\sym{C}\to\sym{C}^\prime$ that is $k$-linear,
essentially surjective, fully faithful, strong monoidal and ribbon.
\end{definition}

\begin{theorem}
Let $\sym{C}$ be a modular category and $\omega$ be the long forgetful
functor. Then
\begin{equation}
\sym{C}\simeq\sym{M}^{\coend(\sym{C},\omega)}
\end{equation}
are equivalent as modular categories.
\end{theorem}

\begin{proof}
The functor $F$ of Theorem~\ref{thm_monequiv} has these properties.
\end{proof}

\begin{corollary}
\label{cor_abelian}
Each modular category $\sym{C}$ is abelian, and its long forgetful
functor $\omega$ is exact.
\end{corollary}

\begin{proof}
Since the functor $F\colon\sym{C}\to\sym{M}^{\coend(\sym{C},\omega)}$
is part of an equivalence of categories and the category
$\sym{M}^{\coend(\sym{C},\omega)}$ is abelian, so is $\sym{C}$. As
part of an equivalence, $F$ is exact, and since $U$ of
Theorem~\ref{thm_monequiv} is exact, too, so is $\omega$.
\end{proof}

Finally, we can complete the characterization of the WHA
$H=\coend(\sym{C},\omega)$ reconstructed from a modular category
$\sym{C}$ and the long forgetful functor
$\omega\colon\sym{C}\to\Vect_k$ (Section~\ref{sect_reconstruct}). The
following result complements Theorem~\ref{thm_comodrep}(2).

\begin{theorem}
Let $\sym{C}$ be a modular category with $k=\End(\one)$ and
$\omega\colon\sym{C}\to\Vect_k$ be the long forgetful functor. Then
$H=\coend(\sym{C},\omega)$ is a comodular WHA over $k$.
\end{theorem}

\begin{proof}
$\sym{C}\simeq\sym{M}^H$ are equivalent as modular categories, and
$F\colon\sym{C}\to\sym{M}^H$ of Theorem~\ref{thm_monequiv} is a
$k$-linear, essentially surjective, fully faithful, strong monoidal
ribbon functor. Such a functor preserves simple objects up to
isomorphism and preserves traces
(Proposition~\ref{prop_ribbonfunctortrace}), and so it preserves the
non-degeneracy of the $S$-matrix of Definition~\ref{def_modular}(3) as
well. Then Corollary~\ref{cor_cofactorize} implies weak
cofactorizability of $H$, and so $H$ is comodular.
\end{proof}

\subsection{Morita equivalence and the choice of the forgetful functor}

When we start with a modular category $\sym{C}$ and reconstruct a
comodular Weak Hopf Algebra $H=\coend(\sym{C},\omega)$, we always work
with the canonical choice of the long forgetful functor
$\omega\colon\sym{C}\to\Vect_k$. When we start with a comodular WHA
$H$ over $k$, however, the category $\sym{M}^H$ always comes with a
forgetful functor $U\colon\sym{M}^H\to\Vect_k$
(Proposition~\ref{prop_forgetful}) which may or may not agree with the
canonical choice of the long forgetful functor. In order to better
understand the situation, let us recall the following results
from~\cite[Theorem 2.1.12 and Lemma 2.2.1]{Sc92}.

There is a category $\mathfrak{C}_k$ whose objects $(\sym{C},\omega)$
are \emph{categories over} $\Vect_k$, \ie\ pairs of a small category
$\sym{C}$ with a functor $\omega\colon\sym{C}\to\Vect_k$ that takes
values in $\fdVect_k$. Its morphisms are \emph{functors over}
$\Vect_k$. A functor
$[F,\xi]\colon(\sym{C},\omega)\to(\sym{C}^\prime,\omega^\prime)$ over
$\Vect_k$ is an equivalence class of pairs $(F,\xi)$ where
$F\colon\sym{C}\to\sym{C}^\prime$ is a functor and
$\xi\colon\omega\Rightarrow\omega^\prime\circ F$ a natural
equivalence. The equivalence relation is such that $[F,\xi]$ is an
isomorphism in $\mathfrak{C}_k$ if and only if
$F\colon\sym{C}\to\sym{C}^\prime$ is an equivalence of categories.

Then Tannaka--Kre\v\i n duality between coalgebras and their
categories of finite-dimensional comodules is an adjunction
\begin{equation}
\begin{aligned}
\xymatrix{
\mathfrak{C}_k\ar@/^2ex/[rr]^{\mathrm{Coend}(-)}="p"&&
\CoAlg_k\ar@/^2ex/[ll]^{\sym{M}^-}="q"
\ar@{-|}"p"+<0ex,-2.5ex>;"q"+<0ex,2.5ex>
}
\end{aligned}
\end{equation}
Here $\CoAlg_k$ is the category of coalgebras over $k$ and their
homomorphisms. The functor $\coend$ applied to a category
$(\sym{C},\omega)$ over $\Vect_k$ gives the universal coend, using the
functor $\omega$ supplied, and the functor $\sym{M}^-$ applied to a
coalgebra $H$ gives the category $\sym{M}^H$ of finite-dimensional
right $H$-comodules with the forgetful functor $U^H$ of
Proposition~\ref{prop_forgetful}, viewed as a category
$(\sym{M}^H,U^H)$ over $\Vect_k$.

The counit of the adjunction is always a natural equivalence, \ie\ for
each coalgebra $H$ over $k$, $H\cong\coend(\sym{M}^H,U)$ are
isomorphic as coalgebras. If $(\sym{C},\omega)$ is a category over
$\Vect_k$ for which $\sym{C}$ is $k$-linear abelian and $\omega$ is
$k$-linear, faithful and exact, then the unit of the adjunction is an
isomorphism as well, \ie\ $(\sym{C},\omega)$ is isomorphic in
$\mathfrak{C}_k$ to
$(\sym{M}^{\coend(\sym{C},\omega)},U^{\coend(\sym{C},\omega)})$. This
means that $\sym{C}\simeq\sym{M}^{\coend(\sym{C},\omega)}$ are
equivalent as categories. We have shown this by hand for the case in
which $\sym{C}$ is modular and $\omega$ the long forgetful functor.

If one starts with a comodular WHA $H$ whose underlying functor
$U^H\colon\sym{M}^H\to\Vect_k$ is not naturally isomorphic to the long
forgetful functor $\omega\colon\sym{M}^H\to\Vect_k$, the above
adjunction yields a coalgebra $\coend(\sym{M}^H,U^H)$ that is
isomorphic to $H$, but the coalgebra we have reconstructed in
Section~\ref{sect_reconstruct}, is $\coend(\sym{M}^H,\omega)$ which
need not be isomorphic to $H$ as a coalgebra.

Our $\coend(\sym{M}^H,\omega)$ is in general only Morita equivalent to
$H$. It is simply a canonical choice in the class of all comodular
WHAs whose categories of finite-dimensional comodules are equivalent
to $\sym{M}^H$ as modular categories.

\section{Example}
\label{sect_example}

In this section, we present the reconstructed Weak Hopf Algebra $H$
for the modular category $\sym{C}$ associated with the quantum group
$U_q(\ssl_2)$, $q$ a root of unity. We use the diagrammatic
description of~\cite{KaLi94} and precisely follow their notation.

Let $r\in\{2,3,4,\ldots\}$ and $q=\exp\frac{\pi}{r}$. For simplicity,
we work over the complex numbers $k=\C$. The isomorphism classes of
simple objects of $\sym{C}$ are indexed by the set
$I=\{0,1,\ldots,r-2\}$. By $V_j$, we denote a specific representative
of the class $j\in I$. Its identity morphism is visualized by a
straight line, labeled by $j\in I$,
\begin{equation}
\begin{aligned}
\begin{pspicture}(2,2)
\psline(1,0)(1,2)
\rput(1.5,1.8){$\scriptstyle j$}
\end{pspicture}
\end{aligned}
\end{equation}
All our diagrams are plane projections of oriented framed tangles,
drawn in blackboard framing. The coherence theorem for ribbon
categories~\cite{ReTu90} makes sure that each diagram defines a
morphism of $\sym{C}$. Since $\sym{C}$ is $k$-linear, we can take
formal linear combinations of diagrams with coefficients in $k$. All
our diagrams are read from top to bottom.

Two special features of $U_q(\ssl_2)$ are exploited. First, the simple
objects are isomorphic to their duals, and the choice of
representatives $V_j$, $j\in I$, of the simple objects is such that
${(V_j)}^\ast=V_j$ are equal rather than merely isomorphic. This
allows us to omit any arrows from the diagrams that would indicate the
orientation of the ribbon tangle.

Second, there are no higher multiplicities, \ie\ for all $a,b,c\in I$,
we have $\dim_k\Hom(V_a\otimes V_b,V_c)\in\{0,1\}$. More precisely,
$\Hom(V_a\otimes V_b,V_c)\cong k$ if and only if the triple $(a,b,c)$
is \emph{admissible}. Otherwise, $\Hom(V_a\otimes V_b,V_c)=\{0\}$.

\begin{definition}
A triple $(a,b,c)\in I^3$ is called \emph{admissible} if the following
conditions hold.
\begin{myenumerate}
\item
$a+b+c\equiv 0$ mod $2$ (\emph{parity}),
\item
$a+b-c\geq 0$ and $b+c-a\geq 0$ and $c+a-b\geq 0$ (\emph{quantum triangle inequality}),
\item
$a+b+c\leq 2r-4$ (\emph{non-negligibility}).
\end{myenumerate}
\end{definition}

\noindent
A special choice of basis vector of $\Hom(V_a\otimes V_b,V_c)$ is
denoted by a trivalent vertex,
\begin{equation}
\begin{aligned}
\begin{pspicture}(4,2)
\rput(2,1){\trivalent{a}{b}{c}}
\end{pspicture}
\end{aligned}
\end{equation}
If we draw such a diagram for a triple $(a,b,c)\in I^3$ that is not
admissible, by convention, we multiply the entire diagram by zero. We
denote by $\Delta_j$ the categorical dimension of $V_j$ and by
$\theta(a,b,c)$ the evaluation of the theta graph,
\begin{equation}
\Delta_j =
\begin{aligned}
\begin{pspicture}(4,2)
\pscircle(2,1){0.5}
\rput(1,1){$\scriptstyle j$}
\end{pspicture}
\end{aligned}
\qquad\qquad
\theta(a,b,c) =
\begin{aligned}
\begin{pspicture}(6,2)
\pscircle[fillstyle=solid,fillcolor=black](1,1){0.15}
\pscircle[fillcolor=black](5,1){0.15}
\psbezier(1,1)(1,2.5)(5,2.5)(5,1)
\psbezier(1,1)(1,-0.5)(5,-0.5)(5,1)
\psline(1,1)(5,1)
\rput(2,2.5){$\scriptstyle a$}
\rput(3,1.5){$\scriptstyle b$}
\rput(2,-0.5){$\scriptstyle c$}
\end{pspicture}
\end{aligned}
\end{equation}
Note that $\Delta_j\neq 0$ for all $j\in I$ and $\theta(a,b,c)\neq 0$
for all admissible triples $(a,b,c)\in I^3$. For each $j\in I$, we
define the vector spaces
\begin{equation}
\omega(V_j) = \Hom(\hat V,\hat V\otimes V_j) = \Span_k\Biggl\{\quad
\begin{aligned}
\begin{pspicture}(4,2)
\rput(2,1){\trivalentdown{p}{q}{j}}
\end{pspicture}
\end{aligned}
\quad\Biggl|\Biggr.\quad
p,q\in I\quad\Biggr\}
\end{equation}
and
\begin{equation}
{\omega(V_j)}^\ast = \Hom (\hat V\otimes V_j,\hat V) = \Span_k\Biggl\{\quad
\begin{aligned}
\begin{pspicture}(4,2)
\rput(2,1){\trivalent{p}{j}{q}}
\end{pspicture}
\end{aligned}
\quad\Biggl|\Biggr.\quad
p,q\in I\quad\Biggr\}
\end{equation}
where $\hat V=\oplus_{j\in I}V_j$ denotes the \emph{universal
object}. This notation is compatible with the remainder of the present
article, but \emph{not} with~\cite{KaLi94}. There, the universal
object is called $\omega$ whereas our $\omega$ is the long forgetful
functor. In the following, we prefer the bases
${(e^{(V_j)}_{pq})}_{pq}$ and ${(e^{pq}_{(V_j)})}_{pq}$ with
\begin{equation}
e_{pq}^{(V_j)} =
\begin{aligned}
\begin{pspicture}(4,2)
\rput(2,1){\trivalentdown{p}{q}{j}}
\end{pspicture}
\end{aligned}
\qquad\mbox{and}\qquad
e^{pq}_{(V_j)} = \frac{\Delta_q}{\theta(p,q,j)}
\begin{aligned}
\begin{pspicture}(4,2)
\rput(2,1){\trivalent{p}{j}{q}}
\end{pspicture}
\end{aligned}
\end{equation}
where $p,q\in I$ such that $(p,q,j)$ is admissible. The reconstructed
WHA is the vector space
\begin{equation}
H=\bigoplus_{j\in I}{\omega(V_j)}^\ast\otimes\omega(V_j).
\end{equation}
A convenient basis of $H$ is given by the vectors of the form
\begin{equation}
{[e^{pq}_{(V_j)}|e^{(V_j)}_{rs}]}_{V_j} := e^{pq}_{(V_j)}\otimes e^{(V_j)}_{rs}
\end{equation}
for $j\in I$ and $p,q,r,s\in I$ such that $(p,q,j)$ and $(r,s,j)$ are
admissible. We can now give the coalgebra structure
$(H,\Delta,\epsilon)$:
\begin{eqnarray}
\Delta( {[e^{pq}_{(V_j)}|e^{(V_j)}_{rs}]}_{V_j} ) &=&
\sum_{t,u\in I} {[e^{pq}_{(V_j)}|e^{(V_j)}_{tu}]}_{V_j}\otimes
{[e^{tu}_{(V_j)}|e^{(V_j)}_{rs}]}_{V_j},\\
\epsilon( {[e^{pq}_{(V_j)}|e^{(V_j)}_{rs}]}_{V_j} )&=&
\delta_{ps}\delta_{qr}.
\end{eqnarray}
More generally, the counit is the evaluation of the following trace,
\begin{equation}
\epsilon\Bigl(\Bigl[
\begin{aligned}
\begin{pspicture}(3,2)
\rput(1.5,1){\boxbasisdual{A}}
\end{pspicture}
\end{aligned}
\Bigr|\Bigl.
\begin{aligned}
\begin{pspicture}(3,2)
\rput(1.5,1){\boxbasis{B}}
\end{pspicture}
\end{aligned}
\Bigr]\Bigr) =
\begin{aligned}
\begin{pspicture}(5,4)
\rput(2,3){\boxbasis{B}}
\rput(2,1){\boxbasisdual{A}}
\rput(4,2){\circlebox{D}}
\psbezier(2,0)(2,-1)(4,-1)(4,0)
\psbezier(4,4)(4,5)(2,5)(2,4)
\psline(4,0)(4,1)
\psline(4,3)(4,4)
\end{pspicture}
\end{aligned}
\end{equation}
By this notation, we mean that one takes whatever ribbon tangles $A$
and $B$ occur in the argument of $\epsilon$ and pastes them into the
diagram on the right. All open ends of the tangles are labeled by
simple objects, and when one connects two of them, say labeled by
$p\in I$ and $q\in I$, the composition of morphisms is zero unless
$p=q$, \ie\ one has to write down a prefactor of $\delta_{pq}$. For
example, putting
\begin{equation}
\begin{aligned}
\begin{pspicture}(4,2)
\rput(2,1){\trivalent{r}{j}{s}}
\end{pspicture}
\end{aligned}
\qquad\mbox{below}\qquad
\begin{aligned}
\begin{pspicture}(4,2)
\rput(2,1){\trivalentdown{p}{q}{k}}
\end{pspicture}
\end{aligned}
\qquad\mbox{gives}\qquad\delta_{qr}\delta_{kj}
\begin{aligned}
\begin{pspicture}(4,4)
\pscircle[fillstyle=solid,fillcolor=black](2,3){0.15}
\pscircle[fillstyle=solid,fillcolor=black](2,1){0.15}
\psbezier(2,1)(1,1)(1,3)(2,3)
\psbezier(2,1)(3,1)(3,3)(2,3)
\psline(2,3)(2,4)
\psline(2,1)(2,0)
\rput(2.5,0.2){$\scriptstyle s$}
\rput(0.8,2){$\scriptstyle q$}
\rput(3.2,2){$\scriptstyle k$}
\rput(2.5,3.8){$\scriptstyle p$}
\end{pspicture}
\end{aligned}
\end{equation}
The morphism $D$ is defined as
\begin{equation}
\begin{aligned}
\begin{pspicture}(2,4)
\rput(1,2){\circlebox{D}}
\psline(1,3)(1,3.5)
\psline(1,1)(1,0.5)
\end{pspicture}
\end{aligned}
\quad=\sum_{j\in I}\Delta_j^{-1}
\begin{aligned}
\begin{pspicture}(1,3)
\psline(0.5,0)(0.5,3)
\rput(1,2.5){$\scriptstyle j$}
\end{pspicture}
\end{aligned}
\end{equation}
We extend our notation for the elements of $H$ to $\omega(X)=\Hom(\hat
V,\hat V\otimes X)$ and ${\omega(X)}^\ast=\Hom(\hat V\otimes X,\hat
V)$ for any object $X$ of $\sym{C}$. Such elements are denoted by
\begin{equation}
{\Bigl[
\begin{aligned}
\begin{pspicture}(3,2)
\rput(1.5,1){\boxbasisdual{}}
\rput(2.5,2){$\scriptstyle X$}
\end{pspicture}
\end{aligned}
\Bigr|\Bigl.
\begin{aligned}
\begin{pspicture}(3,2)
\rput(1.5,1){\boxbasis{}}
\rput(2.5,0){$\scriptstyle X$}
\end{pspicture}
\end{aligned}
\Bigr]}_X
\end{equation}
and they are indeed elements of $H$ if we impose for each morphism
$f\colon X\to Y$ the relations
\begin{equation}
\Bigl[
\begin{aligned}
\begin{pspicture}(3,5)
\rput(1.5,1){\boxbasisdual{A}}
\psline(1,2)(1,4.5)
\rput(2,3){\circlebox{f}}
\psline(2,4)(2,4.5)
\rput(2.5,1.9){$\scriptstyle Y$}
\rput(2.5,4.1){$\scriptstyle X$}
\end{pspicture}
\end{aligned}
\Bigr|\Bigl.
\begin{aligned}
\begin{pspicture}(3,5)
\rput(1.5,3.5){\boxbasis{B}}
\psline(1,0)(1,2.5)
\psline(2,0)(2,2.5)
\rput(2.5,0){$\scriptstyle X$}
\end{pspicture}
\end{aligned}
\Bigr]_X
\qquad = \qquad
\Bigl[
\begin{aligned}
\begin{pspicture}(3,5)
\rput(1.5,1){\boxbasisdual{A}}
\psline(1,2)(1,4.5)
\psline(2,2)(2,4.5)
\rput(2.5,4){$\scriptstyle Y$}
\end{pspicture}
\end{aligned}
\Bigr|\Bigl.
\begin{aligned}
\begin{pspicture}(3,5)
\rput(1.5,3.5){\boxbasis{B}}
\psline(1,0)(1,2.5)
\rput(2,1.5){\circlebox{f}}
\psline(2,0)(2,0.5)
\rput(2.5,0.4){$\scriptstyle Y$}
\rput(2.5,2.6){$\scriptstyle X$}
\end{pspicture}
\end{aligned}
\Bigr]_Y
\end{equation}
Before we present the algebra structure of $H$, we recall the
recoupling identity
\begin{equation}
\begin{aligned}
\begin{pspicture}(5,4)
\pscircle[fillstyle=solid,fillcolor=black](2,2){0.15}
\pscircle[fillstyle=solid,fillcolor=black](3,2){0.15}
\psline(1,1)(2,2)
\psline(1,3)(2,2)
\psline(2,2)(3,2)
\psline(3,2)(4,3)
\psline(3,2)(4,1)
\rput(0.5,3){$\scriptstyle b$}
\rput(0.5,1){$\scriptstyle a$}
\rput(4.5,3){$\scriptstyle c$}
\rput(4.5,1){$\scriptstyle d$}
\rput(2.5,1.5){$\scriptstyle j$}
\end{pspicture}
\end{aligned}
\qquad=\sum_{i\in I}
\left\{\begin{matrix}a&b&i\\ c&d&j\end{matrix}\right\}_q
\begin{aligned}
\begin{pspicture}(4,5)
\pscircle[fillstyle=solid,fillcolor=black](2,2){0.15}
\pscircle[fillstyle=solid,fillcolor=black](2,3){0.15}
\psline(1,1)(2,2)
\psline(3,1)(2,2)
\psline(2,2)(2,3)
\psline(2,3)(1,4)
\psline(2,3)(3,4)
\rput(0.5,1){$\scriptstyle a$}
\rput(3.5,1){$\scriptstyle d$}
\rput(0.5,4){$\scriptstyle b$}
\rput(3.5,4){$\scriptstyle c$}
\rput(2.5,2.5){$\scriptstyle i$}
\end{pspicture}
\end{aligned}
\end{equation}
which holds whenever the triples $(a,b,j)$ and $(c,d,j)$ on the left
hand side are admissible. Here the \emph{quantum $6j$-symbols} can be
computed as follows.
\begin{equation}
\left\{\begin{matrix}a&b&i\\ c&d&j\end{matrix}\right\}_q =
\frac{\Delta_i}{\theta(a,d,i)\theta(b,c,i)}
\begin{aligned}
\begin{pspicture}(5,4)
\pscircle[fillstyle=solid,fillcolor=black](1,2){0.15}
\pscircle[fillstyle=solid,fillcolor=black](3,2){0.15}
\pscircle[fillstyle=solid,fillcolor=black](2,1){0.15}
\pscircle[fillstyle=solid,fillcolor=black](2,3){0.15}
\psline(2,1)(1,2)
\psline(2,1)(3,2)
\psline(1,2)(3,2)
\psline(1,2)(2,3)
\psline(3,2)(2,3)
\psbezier(2,1)(2,0)(4,0)(4,1)
\psbezier(2,3)(2,4)(4,4)(4,3)
\psline(4,1)(4,3)
\rput(1.1,1.3){$\scriptstyle a$}
\rput(2.9,1.3){$\scriptstyle d$}
\rput(1.1,2.7){$\scriptstyle b$}
\rput(2.9,2.7){$\scriptstyle c$}
\rput(2,1.6){$\scriptstyle j$}
\rput(4.5,2){$\scriptstyle i$}
\end{pspicture}
\end{aligned}
\end{equation}
The algebra structure $(H,\mu,\eta)$ is given by
\begin{eqnarray}
\eta(1) &=& \sum_{j,\ell\in I}
{[e^{jj}_{(V_0)}|e_{\ell\ell}^{(V_0)}]}_{V_0},\\
\mu\Bigl(\Bigl[
\begin{aligned}
\begin{pspicture}(3,2)
\rput(1.5,1){\boxbasisdual{A}}
\rput(2.5,2){$\scriptstyle X$}
\end{pspicture}
\end{aligned}
\Bigr|\Bigl.
\begin{aligned}
\begin{pspicture}(3,2)
\rput(1.5,1){\boxbasis{B}}
\rput(2.5,0){$\scriptstyle X$}
\end{pspicture}
\end{aligned}
\Bigr]_X\otimes\Bigl[
\begin{aligned}
\begin{pspicture}(3,2)
\rput(1.5,1){\boxbasisdual{C}}
\rput(2.5,2){$\scriptstyle Y$}
\end{pspicture}
\end{aligned}
\Bigr|\Bigl.
\begin{aligned}
\begin{pspicture}(3,2)
\rput(1.5,1){\boxbasis{D}}
\rput(2.5,0){$\scriptstyle Y$}
\end{pspicture}
\end{aligned}
\Bigr]_Y\Bigr) &=& \Bigl[
\begin{aligned}
\begin{pspicture}(4,4)
\rput(1,2.5){\boxbasisdual{A}}
\rput(1.5,1){\boxbasisdualwide{C}}
\psline(2.5,2)(2.5,3.5)
\rput(2,3.5){$\scriptstyle X$}
\rput(3,3.5){$\scriptstyle Y$}
\end{pspicture}
\end{aligned}
\Bigr|\Bigl.
\begin{aligned}
\begin{pspicture}(4,4)
\rput(1.5,2.5){\boxbasiswide{D}}
\rput(1,1){\boxbasis{B}}
\psline(2.5,1.5)(2.5,0)
\rput(2,0){$\scriptstyle X$}
\rput(3,0){$\scriptstyle Y$}
\end{pspicture}
\end{aligned}
\Bigr]_{X\otimes Y}
\end{eqnarray}
In terms of our favourite bases, the multiplication reads
\begin{eqnarray}
&&\mu({[e^{pq}_{(V_j)}|e_{rs}^{(V_j)}]}_{V_j}\otimes{[e^{ab}_{(V_\ell)}|e_{cd}^{(V_\ell)}]}_{V_\ell})\\
&=& \delta_{qa}\delta_{rd}\sum_{u\in I}
{[e^{pb}_{(V_u)}|e_{cs}^{(V_u)}]}_{V_u}\,
\left\{\begin{matrix}p&j&u\\ \ell&b&a\end{matrix}\right\}_q\,
\left\{\begin{matrix}c&\ell&u\\j&s&d\end{matrix}\right\}_q\,
\frac{\Delta_a\,\theta(p,b,u)\theta(j,\ell,u)}{\Delta_u\,\theta(p,a,j)\theta(a,b,\ell)}.\nn
\end{eqnarray}
The antipode of $H$ is given by
\begin{equation}
S\Bigl(\Bigl[
\begin{aligned}
\begin{pspicture}(3,2)
\rput(1.5,1){\boxbasisdual{A}}
\rput(2.5,2){$\scriptstyle X$}
\end{pspicture}
\end{aligned}
\Bigr|\Bigl.
\begin{aligned}
\begin{pspicture}(3,2)
\rput(1.5,1){\boxbasis{B}}
\rput(2.5,0){$\scriptstyle X$}
\end{pspicture}
\end{aligned}
\bigr]_X\Bigr) = \Bigl[
\begin{aligned}
\begin{pspicture}(4,4)
\rput(1.5,2){\boxbasis{B}}
\psbezier(2,1)(2,0)(3,0)(3,1)
\psline(3,1)(3,5)
\rput(1.5,4){\circlebox{D}}
\rput(1,0){\circlebig{D^{-1}}}
\psline(1,-1)(1,-1.5)
\rput(3.5,5){$\scriptstyle X$}
\end{pspicture}
\end{aligned}
\Bigr|\Bigl.
\begin{aligned}
\begin{pspicture}(4,2)
\rput(1.5,1){\boxbasisdual{A}}
\psbezier(2,2)(2,3)(3,3)(3,2)
\psline(3,2)(3,0)
\rput(3.5,0){$\scriptstyle X$}
\end{pspicture}
\end{aligned}
\Bigr]_X
\end{equation}
which reads in our basis
\begin{equation}
S({[e^{pq}_{(V_j)}|e_{rs}^{(V_j)}]}_{V_j} =
{[e^{rs}_{(V_j)}|e_{pq}^{(V_j)}]}_{V_j}\,
\frac{\Delta_q\,\theta(r,s,j)}{\Delta_r\,\theta(p,q,j)}.
\end{equation}
We finally list the coquasitriangular structure
\begin{equation}
r\Bigl(\Bigl[
\begin{aligned}
\begin{pspicture}(3,2)
\rput(1.5,1){\boxbasisdual{A}}
\rput(2.5,2){$\scriptstyle X$}
\end{pspicture}
\end{aligned}
\Bigr|\Bigl.
\begin{aligned}
\begin{pspicture}(3,2)
\rput(1.5,1){\boxbasis{B}}
\rput(2.5,0){$\scriptstyle X$}
\end{pspicture}
\end{aligned}
\Bigr]_X\otimes\Bigl[
\begin{aligned}
\begin{pspicture}(3,2)
\rput(1.5,1){\boxbasisdual{C}}
\rput(2.5,2){$\scriptstyle Y$}
\end{pspicture}
\end{aligned}
\Bigr|\Bigl.
\begin{aligned}
\begin{pspicture}(3,2)
\rput(1.5,1){\boxbasis{D}}
\rput(2.5,0){$\scriptstyle Y$}
\end{pspicture}
\end{aligned}
\Bigr]_Y\Bigr) =
\begin{aligned}
\begin{pspicture}(6,10)
\rput(2,1){\boxbasisdualwide{C}}
\rput(1.5,3){\boxbasisdual{A}}
\rput(1.5,7){\boxbasis{D}}
\rput(2,9){\boxbasiswide{B}}
\psbezier(2.25,0)(2.25,-1)(5,-1)(5,0)
\psbezier(2.25,10)(2.25,11)(5,11)(5,10)
\rput(5,5){\circlebox{D}}
\psline(5,0)(5,4)
\psline(5,10)(5,6)
\psline(1,4)(1,6)
\psbezier(2,6)(2,5.5)(2.3,5.2)(2.5,5)
\psbezier(2.5,5)(2.7,4.8)(3,4)(3,2)
\psbezier(2,4)(2,4.2)(2.2,4.6)(2.3,4.8)
\psbezier(3,8)(3,7)(2.9,5.4)(2.7,5.2)
\end{pspicture}
\end{aligned}
\end{equation}
and the universal ribbon form
\begin{equation}
\nu\Bigl(\Bigl[
\begin{aligned}
\begin{pspicture}(3,2)
\rput(1.5,1){\boxbasisdual{A}}
\rput(2.5,2){$\scriptstyle X$}
\end{pspicture}
\end{aligned}
\Bigr|\Bigl.
\begin{aligned}
\begin{pspicture}(3,2)
\rput(1.5,1){\boxbasis{B}}
\rput(2.5,0){$\scriptstyle X$}
\end{pspicture}
\end{aligned}
\Bigr]_X\Bigr) =
\begin{aligned}
\begin{pspicture}(5,5)
\rput(1.5,1){\boxbasisdual{A}}
\rput(1.5,4){\boxbasis{B}}
\psbezier(1.5,0)(1.5,-1)(4,-1)(4,0)
\psbezier(1.5,5)(1.5,6)(4,6)(4,5)
\rput(4,2.5){\circlebox{D}}
\psline(4,0)(4,1.5)
\psline(4,3.5)(4,5)
\psline(1,2)(1,3)
\psbezier(2,3)(2,2.5)(2,2)(2.5,2)
\psbezier(2.5,2)(2.8,2)(3,2.3)(3,2.5)
\psbezier(3,2.5)(3,3)(2.5,2.8)(2.3,2.3)
\end{pspicture}
\end{aligned}
\end{equation}

\appendix
\section{Background on tensor categories}
\label{app_tensor}

In this appendix, we collect the relevant definitions and properties
of monoidal, autonomous, braided monoidal, ribbon and abelian
categories, following Freyd--Yetter~\cite{FrYe92},
Schauenburg~\cite{Sc92}, Turaev~\cite{Tu94} and
MacLane~\cite{Ma73}. In order to keep the appendix short, we write
down identities involving morphisms rather than the more familiar
commutative diagrams.

\subsection{Monoidal categories}
\label{app_monoidal}

\begin{definition}
A \emph{monoidal category}
$(\sym{C},\otimes,\one,\alpha,\lambda,\rho)$ is a category $\sym{C}$
with a bifunctor $\otimes\colon\sym{C}\times\sym{C}\to\sym{C}$
(\emph{tensor product}), an object $\one\in|\sym{C}|$ (\emph{monoidal
unit}) and natural isomorphisms $\alpha_{X,Y,Z}\colon(X\otimes
Y)\otimes Z\to X\otimes(Y\otimes Z)$ (\emph{associator}),
$\lambda_X\colon \one\otimes X\to X$ (\emph{left-unit constraint}) and
$\rho_X\colon X\otimes\one\to X$ (\emph{right-unit constraint}) for
all $X,Y,Z\in|\sym{C}|$, subject to the pentagon axiom
\begin{equation}
\alpha_{X,Y,Z\otimes W}\circ\alpha_{X\otimes Y,Z,W}
= (\id_X\otimes\alpha_{Y,Z,W})\circ\alpha_{X,Y\otimes Z,W}\circ(\alpha_{X,Y,Z}\otimes\id_W)
\end{equation}
and the triangle axiom
\begin{equation}
\rho_X\otimes\id_Y=(\id_X\otimes\lambda_Y)\circ\alpha_{X,\one,Y}
\end{equation}
for all $X,Y,Z,W\in|\sym{C}|$.
\end{definition}

\begin{definition}
\label{def_lax}
Let $(\sym{C},\otimes,\one,\alpha,\lambda,\rho)$ and
$(\sym{C}^\prime,\otimes^\prime,\one^\prime,\alpha^\prime,\lambda^\prime,\rho^\prime)$
be monoidal categories.
\begin{myenumerate}
\item
A \emph{lax monoidal functor}
$(F,F_{X,Y},F_0)\colon\sym{C}\to\sym{C}^\prime$ consists of a
functor $F\colon\sym{C}\to\sym{C}^\prime$, morphisms $F_{X,Y}\colon
FX\otimes^\prime FY\to F(X\otimes Y)$ that are natural in
$X,Y\in|\sym{C}|$, and of a morphism $F_0\colon\one^\prime\to
F\one$, subject to the hexagon axiom
\begin{equation}
F_{X,Y\otimes Z}\circ(\id_{FX}\otimes^\prime F_{Y,Z})\circ\alpha^\prime_{FX,FY,FZ}
= F\alpha_{X,Y,Z}\circ F_{X\otimes Y,Z}\circ(F_{X,Y}\otimes^\prime\id_{FZ})
\end{equation}
and the two squares
\begin{eqnarray}
\lambda^\prime_{FX} &=& F\lambda_X\circ F_{\one,X}\circ(F_0\otimes^\prime\id_{FX}),\\
\rho^\prime_{FX}    &=& F\rho_X\circ F_{X,\one}\circ(\id_{FX}\otimes^\prime F_0)
\end{eqnarray}
for all $X,Y,Z\in|\sym{C}|$.
\item
An \emph{oplax monoidal functor}
$(F,F^{X,Y},F^0)\colon\sym{C}\to\sym{C}^\prime$ consists of a
functor $F\colon\sym{C}\to\sym{C}^\prime$, morphisms $F^{X,Y}\colon
F(X\otimes Y)\to FX\otimes^\prime FY$ that are natural in
$X,Y\in|\sym{C}|$, and of a morphism $F^0\colon
F\one\to\one^\prime$, subject to the hexagon axiom
\begin{equation}
\label{eq_oplaxhexagon}
(\id_{FX}\otimes^\prime F^{Y,Z})\circ F^{X,Y\otimes Z}\circ F\alpha_{X,Y,Z}
= \alpha^\prime_{FX,FY,FZ}\circ(F^{X,Y}\otimes^\prime\id_{FZ})\circ F^{X\otimes Y,Z}
\end{equation}
and the two squares
\begin{eqnarray}
F\lambda_X &=& \lambda^\prime_{FX}\circ(F^0\otimes^\prime\id_{FX})\circ F^{\one,X},\\
F\rho_X    &=& \rho^\prime_{FX}\circ(\id_{FX}\otimes^\prime F^0)\circ F^{X,\one}
\end{eqnarray}
for all $X,Y,Z\in|\sym{C}|$.
\item
A \emph{strong monoidal functor}
$(F,F_{X,Y},F_0)\colon\sym{C}\to\sym{C}^\prime$ is a lax monoidal
functor such that all $F_{X,Y}$, $X,Y\in|\sym{C}|$ and $F_0$ are
isomorphisms.
\end{myenumerate}
\end{definition}

\begin{definition}
Let $(F,F_{X,Y},F_0)\colon\sym{C}\to\sym{C}^\prime$ and
$(G,G_{X,Y},G_0)\colon\sym{C}\to\sym{C}^\prime$ be lax monoidal
functors between monoidal categories $\sym{C}$ and $\sym{C}^\prime$. A
\emph{monoidal natural transformation} $\eta\colon F\Rightarrow G$ is
a natural transformation such that
\begin{equation}
\eta_{X\otimes Y}\circ F_{X,Y} = G_{X,Y}\circ(\eta_X\otimes^\prime\eta_Y)
\end{equation}
for all $X,Y\in\sym{C}$.
\end{definition}

There is a similar notion of monoidal natural transformation if the
functors are oplax rather than lax. Compositions of [lax, oplax,
strong] monoidal functors are again [lax, oplax, strong] monoidal. The
following result is well known, but quite laborious to verify.

\begin{proposition}
\label{prop_monadjunction}
Let $\sym{C}$ and $\sym{C}^\prime$ be monoidal categories and $F\dashv
G\colon\sym{C}^\prime\to\sym{C}$ be an adjunction with unit
$\eta\colon 1_{\sym{C}}\Rightarrow G\circ F$ and counit
$\epsilon\colon F\circ G\Rightarrow 1_{\sym{C}^\prime}$.
\begin{myenumerate}
\item
If $F$ has an oplax monoidal structure $(F,F^{C_1,C_2},F^0)$, then
$G$ has a lax monoidal structure $(G,G_{D_1,D_2},G_0)$ as follows,
\begin{eqnarray}
G_{D_1,D_2} &=& G(\epsilon_{D_1}\otimes\epsilon_{D_2})
\circ G(F^{G(D_1),G(D_2)})\circ\eta_{G(D_1)\otimes G(D_2)},\\
G_0 &=& G(F^0)\circ\eta_\one.
\end{eqnarray}
\item
If $F$ is strong monoidal, then both $\eta$ and $\epsilon$ are monoidal
natural transformations.
\item
If $F$ is strong monoidal and the adjunction is an equivalence, then
$G$ is strong monoidal.
\end{myenumerate}
\end{proposition}

By an equivalence of monoidal categories, we mean an equivalence of
categories such that one of the functors is strong monoidal. One can
then chose the other functor in such a way that one has an adjoint
equivalence and apply Proposition~\ref{prop_monadjunction}, items~(2)
and~(3). We denote such an equivalence by
$\sym{C}\simeq_\otimes\sym{D}$.

\subsection{Duality}
\label{app_dual}

\begin{definition}
Let $(\sym{C},\otimes,\one,\alpha,\lambda,\rho)$ be a monoidal
category.
\begin{myenumerate}
\item
A \emph{left-dual} $(X^\ast,\ev_X,\coev_X)$ of an object
$X\in|\sym{C}|$ consists of an object $X^\ast\in|\sym{C}|$ and
morphisms $\ev_X\colon X^\ast\otimes X\to\one$ (\emph{left
evaluation}) and $\coev_X\colon\one\to X\otimes X^\ast$ (\emph{left
coevaluation}) that satisfy the triangle identities
\begin{eqnarray}
\label{eq_zigzag1}
\rho_X\circ(\id_X\otimes\ev_X)\circ\alpha_{X,X^\ast,X}
\circ(\coev_X\otimes\id_X)\circ\lambda_X^{-1} &=& \id_X,\\
\label{eq_zigzag2}
\lambda_{X^\ast}\circ(\ev_X\otimes\id_{X^\ast})\circ\alpha^{-1}_{X^\ast,X,X^\ast}
\circ(\id_{X^\ast}\otimes\coev_X)\circ\rho_{X^\ast}^{-1} &=& \id_{X^\ast}.
\end{eqnarray}
\item
A \emph{right-dual} $(\bar X,\bar\ev_X,\bar\coev_X)$ of
$X\in|\sym{C}|$ consists of an object $\bar X\in|\sym{C}|$ and
morphisms $\bar\ev_X\colon X\otimes\bar X\to\one$ (\emph{right
evaluation}) and $\bar\coev_X\colon\one\to\bar X\otimes X$
(\emph{right coevaluation}) that satisfy the triangle identities
\begin{eqnarray}
\label{eq_zigzag3}
\lambda_X\circ(\bar\ev_X\otimes\id_X)\circ\alpha^{-1}_{X,\bar X,X}
\circ(\id_X\otimes\bar\coev_X)\circ\rho_X^{-1} &=& \id_X,\\
\label{eq_zigzag4}
\rho_{\bar X}\circ(\id_{\bar X}\otimes\bar\ev_X)\circ\alpha_{\bar X,X,\bar X}
\circ(\bar\coev_X\otimes\id_{\bar X})\circ\lambda^{-1}_{\bar X} &=& \id_{\bar X}.
\end{eqnarray}
\end{myenumerate}
\end{definition}

\begin{definition}
Let $(\sym{C},\otimes,\one,\alpha,\lambda,\rho)$ be a monoidal category
and $f\colon X\to Y$ be a morphism of $\sym{C}$.
\begin{myenumerate}
\item
If both $X$ and $Y$ have left-duals, the \emph{left-dual} of $f$ is defined as
\begin{equation}
f^\ast := \lambda_{X^\ast}\circ(\ev_Y\otimes\id_{X^\ast})\circ\alpha^{-1}_{Y^\ast,Y,X^\ast}
\circ(\id_{Y^\ast}\otimes(f\otimes\id_{X^\ast}))\circ(\id_{Y^\ast}\otimes\coev_X)\circ\rho^{-1}_{Y^\ast}.
\end{equation}
\item
If both $X$ and $Y$ have right-duals, the \emph{right-dual} of $f$ is defined as
\begin{equation}
\bar f:=\rho_{\bar X}\circ(\id_{\bar X}\otimes\bar\ev_Y)\circ\alpha_{\bar X,Y,\bar Y}
\circ((\id_{\bar X}\otimes f)\otimes\id_{\bar Y})\circ(\bar\coev_X\otimes\id_{\bar Y})\circ\lambda^{-1}_{\bar Y}.
\end{equation}
\end{myenumerate}
\end{definition}

\begin{definition}
A \emph{[left-, right-]autonomous category} is a monoidal category in
which each object is equipped with a specified [left-, right-]dual. An
\emph{autonomous category} is a monoidal category that is both left-
and right-autonomous.
\end{definition}

Note that every autonomous category is monoidally closed because the
functor $-\otimes X^\ast$ is a right adjoint of $-\otimes X$ and $\bar
X\otimes -$ is a right-adjoint of $X\otimes -$ for all
$X\in|\sym{C}|$. In particular, the tensor product in an autonomous
category preserves colimits in both arguments.

\subsection{Ribbon categories}
\label{app_ribbon}

\begin{definition}
A \emph{braided monoidal category}
$(\sym{C},\otimes,\one,\alpha,\lambda,\rho,\sigma)$ is a monoidal
category $(\sym{C},\otimes,\one,\alpha,\lambda,\rho)$ with natural
isomorphisms $\sigma_{X,Y}\colon X\otimes Y\to Y\otimes X$ for all
$X,Y\in|\sym{C}|$ that satisfy the two hexagon axioms
\begin{eqnarray}
\sigma_{X\otimes Y,Z} &=& \alpha_{Z,X,Y}\circ(\sigma_{X,Z}\otimes\id_Y)
\circ\alpha^{-1}_{X,Z,Y}\circ(\id_X\otimes\sigma_{Y,Z})\circ\alpha_{X,Y,Z},\\
\sigma_{X,Y\otimes Z} &=& \alpha^{-1}_{Y,Z,X}\circ(\id_Y\otimes\sigma_{X,Z})
\circ\alpha_{Y,X,Z}\circ(\sigma_{X,Y}\otimes\id_Z)\circ\alpha^{-1}_{X,Y,Z}
\end{eqnarray}
for all $X,Y,Z\in|\sym{C}|$. The category is called \emph{symmetric
monoidal} if in addition
\begin{equation}
\sigma_{Y,X}\circ\sigma_{X,Y} = \id_{X\otimes Y}
\end{equation}
for all $X,Y\in|\sym{C}|$.
\end{definition}

\begin{definition}
Let $(\sym{C},\otimes,\one,\alpha,\lambda,\rho,\sigma)$ and
$(\sym{C}^\prime,\otimes^\prime,\one^\prime,\alpha^\prime,\lambda^\prime,\rho^\prime,\sigma^\prime)$
be braided monoidal categories. A lax monoidal functor
$(F,F_{X,Y},F_0)\colon\sym{C}\to\sym{C}^\prime$ is called
\emph{braided} if
\begin{equation}
\label{eq_braidedfunctor}
F\sigma_{X,Y}\circ F_{X,Y} = F_{Y,X}\circ\sigma^\prime_{FX,FY}
\end{equation}
for all $X,Y\in|\sym{C}|$.
\end{definition}

\begin{proposition}
Let $\sym{C}$ and $\sym{C}^\prime$ be braided monoidal categories and
$F\dashv G\colon\sym{C}^\prime\to\sym{C}$ be an adjoint
equivalence. If $F$ is strong monoidal and braided, then so is $G$.
\end{proposition}

By an equivalence of braided monoidal categories, we therefore mean an
equivalence of categories one functor of which is strong monoidal and
braided.

\begin{definition}
\label{def_ribboncat}
A \emph{ribbon category}
$(\sym{C},\otimes,\one,\alpha,\lambda,\rho,{(-)}^\ast,\ev,\coev,\sigma,\nu)$
is a left-autonomous category
$(\sym{C},\otimes,\one,\alpha,\lambda,\rho,{(-)}^\ast,\ev,\coev)$ that
is braided monoidal as
$(\sym{C},\otimes,\one,\alpha,\lambda,\rho,\sigma)$ with natural
isomorphisms (\emph{ribbon twist}) $\nu_X\colon X\to X$ such that
\begin{equation}
\label{eq_ribbon1}
\nu_{X\otimes Y} = \sigma_{Y,X}\circ\sigma_{X,Y}\circ(\nu_X\otimes\nu_Y)
\end{equation}
and
\begin{equation}
\label{eq_ribbon2}
(\nu_X\otimes\id_{X^\ast})\circ\coev_X = (\id_X\otimes\nu_{X^\ast})\circ\coev_X
\end{equation}
for all $X,Y\in|\sym{C}|$.
\end{definition}

Note that in every ribbon category $\sym{C}$, there are natural
isomorphisms $\tau_X\colon X\to {X^\ast}^\ast$ for all
$X\in|\sym{C}|$, given by
\begin{equation}
\label{eq_pivotal}
\tau_X = \lambda_{{X^\ast}^\ast}\circ(\ev_X\otimes\id_{{X^\ast}^\ast})
\circ(\sigma_{X,X^\ast}\otimes\id_{{X^\ast}^\ast})
\circ(\nu_X\otimes\coev_{X^\ast})\circ\rho_X^{-1},
\end{equation}
that satisfy $(\tau_X)\ast = \tau_{X^\ast}^{-1}$.

Every ribbon category $\sym{C}$ is not only left-autonomous, but also
right-autonomous with $(\bar X,\bar\ev_x,\bar\coev_X)$ where $\bar
X=X^\ast$ and
\begin{eqnarray}
\label{eq_rightdual1}
\bar\ev_X   &=& \ev_X\circ\sigma_{X,X^\ast}\circ(\nu_X\otimes\id_{X^\ast}),\\
\label{eq_rightdual2}
\bar\coev_X &=& (\id_{X^\ast}\otimes\nu_X)\circ\sigma_{X,X^\ast}\circ\coev_X
\end{eqnarray}
for all $X\in|\sym{C}|$. The left- and the right-dual of any morphism
$f\colon X\to Y$ agree, $f^\ast=\bar f$.

\begin{definition}
Let
$(\sym{C},\otimes,\one,\alpha,\lambda,\rho,{(-)}^\ast,\ev,\coev,\sigma,\nu)$
be a ribbon category, $X\in|\sym{C}|$, and $f\colon X\to X$ be a
morphism of $\sym{C}$. Then we define
\begin{myenumerate}
\item
the \emph{trace} of $f$ by
\begin{equation}
\tr_X(f) := \bar\ev_X\circ(f\otimes\id_{X^\ast})\circ\coev_X\colon\one\to\one,
\end{equation}
\item
the \emph{dimension} of $X$ by
\begin{equation}
\dim(X) := \tr_X(\id_X).
\end{equation}
\end{myenumerate}
\end{definition}

\begin{proposition}
Let
$(\sym{C},\otimes,\one,\alpha,\lambda,\rho,{(-)}^\ast,\ev,\coev,\sigma,\nu)$
be a ribbon category. Then
\begin{myenumerate}
\item
$\tr_X(f) = \tr_{X^\ast}(f^\ast)$ for all $f\colon X\to X$.
\item
$\tr_X(g\circ f) = \tr_Y(f\circ g)$ for all $f\colon X\to Y$ and
$g\colon Y\to X$.
\item
$\tr_{X_1\otimes X_2}(h_1\otimes h_2)=\tr_{X_1}(h_1)\tr_{X_2}(h_2)$
for all $h_j\colon X_j\to X_j$, $j\in\{1,2\}$.
\end{myenumerate}
\end{proposition}

\begin{definition}
Let $\sym{C}$ and $\sym{C}^\prime$ be ribbon categories. A
\emph{ribbon functor} $(F,F_{X,Y},F_0)\colon\sym{C}\to\sym{C}^\prime$
is a lax monoidal functor that is braided and satisfies
\begin{equation}
\label{eq_ribbonfunctor}
F\nu_X=\nu^\prime_{FX}
\end{equation}
for all $X\in|\sym{C}|$.
\end{definition}

\begin{proposition}
Let $\sym{C}$ and $\sym{C}^\prime$ be ribbon categories and $F\dashv
G\colon\sym{C}^\prime\to\sym{C}$ be an adjoint equivalence. If $F$ is
strong monoidal and ribbon, then so is $G$.
\end{proposition}

By an equivalence of ribbon categories we therefore mean an
equivalence of categories one functor of which is strong monoidal and
ribbon. The following proposition states what strong monoidal ribbon
functors do to traces.

\begin{proposition}
Let $\sym{C}$ and $\sym{C}^\prime$ be ribbon categories and
$(F,F_{X,Y},F_0)\colon\sym{C}\to\sym{C}^\prime$ be a strong monoidal
ribbon functor. Then for each morphism $f\colon X\to X$ of $\sym{C}$,
the diagram
\begin{equation}
\label{eq_ribbonfunctortrace}
\begin{aligned}
\xymatrix{
\one^\prime\ar[dd]_{\tr_{FX}(Ff)}\ar[rr]^{F_0}&&
F\one\ar[dd]^{F\tr_X(f)}\\
\\
\one^\prime\ar[rr]_{F_0}&&
F\one
}
\end{aligned}
\end{equation}
commutes.
\end{proposition}

\subsection{Abelian and semisimple categories}
\label{app_abelian}

\begin{definition}
A category $\sym{C}$ is called \emph{$\mathbf{Ab}$-enriched} if it is enriched
in the category $\mathbf{Ab}$ of abelian groups, \ie\ if $\Hom(X,Y)$
is an abelian group for all objects $X,Y\in|\sym{C}|$ and if the
composition of morphisms is $\Z$-bilinear.

Let $k$ be a commutative ring. A category $\sym{C}$ is called
\emph{$k$-linear} if it is enriched in ${}_k\sym{M}$, the category of
$k$-modules, \ie\ if $\Hom(X,Y)$ is a $k$-module for all
$X,Y\in|\sym{C}|$ and if the composition of morphisms is $k$-bilinear.

A functor $F\colon\sym{C}\to\sym{C}^\prime$ between [$\mathbf{Ab}$-enriched,
$k$-linear] categories is called [\emph{additive}, $k$-\emph{linear}]
if it induces homomorphisms of [additive groups, $k$-modules]
\begin{equation}
\Hom(X,Y)\to\Hom(FX,FY)
\end{equation}
for all $X,Y\in|\sym{C}|$.
\end{definition}

\begin{definition}
\label{def_preadditivecat}
A monoidal category $(\sym{C},\otimes,\one,\alpha,\lambda,\rho)$ is
called [\emph{$\mathbf{Ab}$-enriched}, $k$-\emph{linear}] if $\sym{C}$ is
[$\mathbf{Ab}$-enriched, $k$-linear] and if the tensor product of morphisms is
[$\Z$-bilinear, $k$-bilinear].
\end{definition}

\begin{definition}
An \emph{additive category} is an $\mathbf{Ab}$-enriched category that
has a terminal object and all binary products. A \emph{preabelian
category} is an $\mathbf{Ab}$-enriched category that has all finite
limits. An \emph{abelian category} is a preabelian category in which
every monomorphism is a kernel and in which every epimorphism is a
cokernel.

A functor $F\colon\sym{C}\to\sym{C}^\prime$ between preabelian
categories is called \emph{exact} if it preserves all finite limits.
\end{definition}

Recall that in an $\mathbf{Ab}$-enriched category, an object is
terminal if and only if it is initial and if and only if it is
null. Every additive category has all finite biproducts. An
equivalence of [$\mathbf{Ab}$-enriched, $k$-linear] categories is an
equivalence of categories one functor of which is [additive,
$k$-linear].

\begin{definition}
\label{def_semisimple}
Let $\sym{C}$ be a $k$-linear category, $k$ a commutative ring.
\begin{myenumerate}
\item
An object $X\in|\sym{C}|$ is called \emph{simple} if $\End(X)\cong
k$ are isomorphic as $k$-modules.
\item
An object $X\in|\sym{C}|$ is called \emph{null} if
$\End(X)\cong\{0\}$.
\item
The category $\sym{C}$ is called \emph{semisimple} if there exists a
family ${\{V_j\}}_{j\in I}$ of objects $V_j\in|\sym{C}|$, $I$ some
index set, such that
\begin{myenumerate}
\item
$V_j$ is simple for all $j\in I$.
\item
$\Hom(V_j,V_\ell)=\{0\}$ for all $j,\ell\in I$ for which $j\neq\ell$.
\item
For each object $X\in|\sym{C}|$, there is a finite sequence
$j_1^{(X)},\ldots,j_{n^X}^{(X)}\in I$, $n^X\in\N_0$, and morphisms
$\imath_\ell^{(X)}\colon V_{j_\ell}\to X$ and $\pi_\ell^{(X)}\colon
X\to V_{j_\ell}$ such that
\begin{equation}
\id_X = \sum_{\ell=1}^{n^X}\imath^X_\ell\circ\pi^X_\ell.
\end{equation}
and
\begin{equation}
\label{eq_semisimpledualbases}
\pi^X_\ell\circ{\imath}^X_m = \left\{
\begin{matrix}
\id_{V_{j^X_\ell}},&\mbox{if}\quad \ell=m,\\
0,                 &\mbox{else}
\end{matrix}
\right.
\end{equation}
\end{myenumerate}
\item
The category is called \emph{finitely semisimple} (also
\emph{Artinian semisimple}) if it is semisimple with a finite index
set $I$ in condition~(3).
\end{myenumerate}
\end{definition}

\begin{proposition}[see {\cite[Lemma II.4.2.2]{Tu94}}]
\label{prop_semisimplehom}
Let $\sym{C}$ be a $k$-linear category, $k$ a commutative ring. If
$\sym{C}$ is [finitely] semisimple, then there is a [finite] set
$\sym{J}\subseteq|\sym{C}|$ of objects each of which is non-null such
that
\begin{eqnarray}
\label{eq_semisimplehom}
\Phi\colon\bigoplus_{J\in\sym{J}}\Hom(X,J)\otimes\Hom(J,Y) &\to& \Hom(X,Y),\nn\\
f\otimes g &\mapsto& g\circ f,
\end{eqnarray}
is an isomorphism for all $X,Y\in|\sym{C}|$.
\end{proposition}

\begin{lemma}
Let $\sym{C}$ be a $k$-linear category, $k$ a field, and $\Hom(X,Y)$
be a finite-dimensional vector space over $k$ for all
$X,Y\in|\sym{C}|$ and let $\sym{J}$ be a set of objects that satisfies
the conditions of Proposition~\ref{prop_semisimplehom}.
\begin{myenumerate}
\item
Each $J\in\sym{J}$ is simple.
\item
If $X\in|\sym{C}|$ is simple, then there exists some $J_X\in\sym{J}$
such that $X\cong J_X$. For all other $J\in\sym{J}$, $J\not\cong
J_X$, we have $\Hom(X,J)=\{0\}=\Hom(J,X)$.
\item
If $X,Y\in|\sym{C}|$ are both simple, then either $X\cong Y$ or
$\Hom(X,Y)=\{0\}$.
\end{myenumerate}
\end{lemma}

\begin{proof}
The idea for this proof is that both source and target of the
isomorphism~\eqref{eq_semisimplehom} are finite-dimensional vector
spaces over $k$. We can therefore count dimensions.
\end{proof}

\begin{corollary}
\label{cor_allsimplethere}
Let $\sym{C}$ be a semisimple $k$-linear category with family
${\{V_j\}}_{j\in I}$ of simple objects, $k$ a field, and $\Hom(X,Y)$
be a finite-dimensional vector space over $k$ for all
$X,Y\in|\sym{C}|$. If $X\in\sym{C}$ is simple, then there exists some
$j\in V_j$ such that $X\cong V_j$.
\end{corollary}

\subsection{$\mathbf{Ab}$-enriched and non-degenerate ribbon categories}

\begin{proposition}
Let
$(\sym{C},\otimes,\one,\alpha,\lambda,\rho,{(-)}^\ast,\ev,\coev,\sigma,\nu)$
be an $\mathbf{Ab}$-enriched ribbon category.
\begin{myenumerate}
\item
The abelian group $k:=\End(\one)$ is a unital commutative ring
with respect to the composition of morphisms.
\item
The category $\sym{C}$ is $k$-linear as a monoidal category.
\item
For all objects $X\in|\sym{C}|$, the trace
\begin{equation}
\tr_X\colon\Hom(X,X)\to k,
\end{equation}
is $k$-linear.
\item
For all objects $X,Y\in|\sym{C}|$, the map
\begin{equation}
\label{eq_traceform}
\Hom(Y,X)\otimes\Hom(X,Y)\to k,\quad f\otimes g\mapsto\tr_X(f\circ g)
\end{equation}
is $k$-bilinear.
\end{myenumerate}
\end{proposition}

\begin{definition}
\label{def_nondegenerate}
An $\mathbf{Ab}$-enriched ribbon category
$(\sym{C},\otimes,\one,\alpha,\lambda,\rho,{(-)}^\ast,\ev,\coev,\sigma,\nu)$
is called \emph{non-degenerate} if the bilinear
forms~\eqref{eq_traceform} are non-degenerate for all objects
$X,Y\in|\sym{C}|$, \ie\ if $\tr_X(f\circ g)=0$ for all $g\colon X\to
Y$ implies $f=0$.
\end{definition}

If we work with semisimple ribbon categories, we also require the set
of representatives of the simple objects to contain the monoidal unit
and to be closed under duality.

\begin{definition}
\label{def_semisimpleribbon}
An $\mathbf{Ab}$-enriched ribbon category
$(\sym{C},\otimes,\one,\alpha,\lambda,\rho,{(-)}^\ast,\ev,\coev,\sigma,\nu)$,
$k=\End(\one)$, is called [\emph{finitely}] \emph{semisimple} if the
underlying $k$-linear category is [finitely] semisimple and the family
${\{V_j\}}_j$ of Definition~\ref{def_semisimple}(3) satisfies the
following conditions.
\begin{myenumerate}
\item
There is an element $0\in I$ such that $V_0\cong\one$.
\item
For each $j\in I$, there is some $j^\ast\in I$ such that $V_{j^\ast}\cong V_j^\ast$.
\end{myenumerate}
\end{definition}

\begin{proposition}
\label{prop_invertibledimension}
Let $\sym{C}$ be a semisimple $\mathbf{Ab}$-enriched ribbon category
with family ${\{V_j\}}_{j\in I}$ as in
Definition~\ref{def_semisimple}(3). Then for all $j\in I$, $\dim V_j$
is invertible in $k$~\cite[Lemma II.4.2.4]{Tu94}.
\end{proposition}

\noindent
The following proposition gives conditions under which ribbon functors
preserve traces.

\begin{proposition}
\label{prop_ribbonfunctortrace}
Let $\sym{C}$ and $\sym{C}^\prime$ be semisimple $k$-linear ribbon
categories, $k$ a field, and
$(F,F_{X,Y},F_0)\colon\sym{C}\to\sym{C}^\prime$ be a strong monoidal
$k$-linear ribbon functor. Then for each morphism $f\colon X\to X$ of
$\sym{C}$,
\begin{equation}
\tr_X(f) = \tr_{FX}(Ff).
\end{equation}
\end{proposition}

\begin{proof}
Since the monoidal units of $\sym{C}$ and $\sym{C}^\prime$ are simple,
$\tr_X(f)=\lambda_f\id_{\one}$ and
$\tr_{Ff}(Ff)=\lambda^\prime_f\id_{\one^\prime}$ for some
$\lambda_f,\lambda_f^\prime\in k$. By $k$-linearity of $F$,
$F\tr_X(f)=\lambda_f\id_{F\one}$, and so~\eqref{eq_ribbonfunctortrace}
implies that $\lambda_f\id_{F\one}=\lambda_f^\prime\id_{F\one}$ and
therefore $\lambda_f=\lambda_f^\prime$.
\end{proof}

\begin{proposition}
\label{prop_nondegenerate}
Let $\sym{C}$ be an $\mathbf{Ab}$-enriched non-degenerate ribbon
category, $k=\End(\one)$ be a field and $\Hom(X,Y)$ be a
finite-dimensional vector space over $k$ for all $X,Y\in|\sym{C}|$. If
$\sym{C}$ satisfies all conditions of a finitely semisimple category
of Definition~\ref{def_semisimple} except maybe
for~\eqref{eq_semisimpledualbases}, then the $\imath_\ell^{(X)}$ and
$\pi_\ell^{(X)}$ can be chosen in such a way
that~\eqref{eq_semisimpledualbases} holds as well.
\end{proposition}

\begin{proof}
Consider the bilinear form
\begin{eqnarray}
\Psi\colon\Hom(X,\hat V)\otimes\Hom(\hat V,X) &\to& k,\nn\\
f\otimes g &\mapsto& \sum_{\ell=1}^{n^{(X)}}\tr_{\hat V}
(f\circ\imath_\ell^{(X)}\circ\pi_\ell^{(X)}\circ g){(\dim V_{j_\ell^{(X)}})}^{-1}
\end{eqnarray}
which is non-degenerate by
Proposition~\ref{prop_invertibledimension}. The
${(\imath_\ell^{(X)})}_{\ell}$ form a basis of $\Hom(\hat V,X)$, and
since $k$ is a field and the $\Hom$ spaces are finite-dimensional
vector spaces, we can choose a dual basis
${(\tilde\pi_\ell^{(X)})}_\ell$ of $\Hom(X,\hat V)$. Then for any
$1\leq p,q\leq n^{(X)}$ and $p\neq q$,
$0=\Psi(\tilde\pi_p^{(X)}\otimes\imath_q^{(X)})$ implies
$\tilde\pi_p^{(X)}\circ\imath_q^{(X)}=0$ by non-degeneracy. Finally,
if $p=q$, $1=\Psi(\tilde\pi_p^{(X)}\otimes\imath_p^{(X)})$ implies
that $\tilde\pi_p^{(X)}\circ\imath_p^{(X)}=\id_{V_{j_p}^{(X)}}$
because $V_{j_p^{(X)}}$ is simple.
\end{proof}

\subsubsection*{Acknowledgements}

The author would like to thank Gabriella B{\"o}hm, Shahn Majid and
Korn{\'e}l Szlach{\'a}nyi for valuable discussions and everyone at
RMKI Budapest for their hospitality.

\newenvironment{hpabstract}{%
\renewcommand{\baselinestretch}{0.2}
\begin{footnotesize}%
}{\end{footnotesize}}%
\newcommand{\hpeprint}[2]{%
\href{http://www.arxiv.org/abs/#1}{\texttt{arxiv:#1#2}}}%
\newcommand{\hpspires}[1]{%
\href{http://www.slac.stanford.edu/spires/find/hep/www?#1}{SPIRES Link}}%
\newcommand{\hpmathsci}[1]{%
\href{http://www.ams.org/mathscinet-getitem?mr=#1}{\texttt{MR #1}}}%
\newcommand{\hpdoi}[1]{%
\href{http://dx.doi.org/#1}{\texttt{DOI #1}}}%
\newcommand{\hpjournal}[2]{%
\href{http://dx.doi.org/#2}{\textsl{#1\/}}}%

\end{document}